\documentclass[twoside,11pt]{article}

\usepackage{jmlr2e,color}
\usepackage{amsmath}\begin{normalsize}
\end{normalsize}
\usepackage{algorithm}
\usepackage[noend]{algpseudocode}
\usepackage{graphicx}
\usepackage{subcaption}
\usepackage{mathtools}
\usepackage{amsthm} 
\usepackage{amsmath}
\usepackage{amssymb}
\usepackage{enumitem}
\usepackage{physics}
\usepackage{setspace}

\DeclareCaptionLabelSeparator{periodquad}{. \;}
\DeclareMathOperator*{\argmax}{arg\,max}
\newcommand\numberthis{\addtocounter{equation}{1}\tag{\theequation}}

\newcommand{\verts}{\; | \;}
\def\wh{\widehat}

\newcommand{\bA}{\boldsymbol{A}}
\newcommand{\ba}{\boldsymbol{a}}
\newcommand{\bB}{\boldsymbol{B}}
\newcommand{\bb}{\boldsymbol{b}}
\newcommand{\bC}{\boldsymbol{C}}

\newcommand{\bD}{\boldsymbol{D}}
\newcommand{\bE}{\boldsymbol{E}}
\newcommand{\be}{\boldsymbol{e}}
\newcommand{\bF}{\boldsymbol{F}}
\newcommand{\bG}{\boldsymbol{G}}
\newcommand{\bH}{\boldsymbol{H}}
\newcommand{\bI}{\boldsymbol{I}}
\newcommand{\bL}{\boldsymbol{L}}
\newcommand{\bM}{\boldsymbol{M}}
\newcommand{\bN}{\boldsymbol{N}}
\newcommand{\bO}{\boldsymbol{O}}
\newcommand{\bP}{\boldsymbol{P}}
\newcommand{\bQ}{\boldsymbol{Q}}

\newcommand{\bR}{\boldsymbol{R}}
\newcommand{\br}{\boldsymbol{r}}

\newcommand{\bT}{\boldsymbol{T}}
\newcommand{\bV}{\boldsymbol{V}}
\newcommand{\bv}{\boldsymbol{v}}
\newcommand{\bW}{\boldsymbol{W}}
\newcommand{\bw}{\boldsymbol{w}}
\newcommand{\bX}{\boldsymbol{X}}
\newcommand{\bx}{\boldsymbol{x}}
\newcommand{\bY}{\boldsymbol{Y}}
\newcommand{\bZ}{\boldsymbol{Z}}

\newcommand{\bxi}{\boldsymbol{\xi}}
\newcommand{\bXi}{\boldsymbol{\Xi}}
\newcommand{\bTheta}{\boldsymbol{\Theta}}
\newcommand{\btheta}{\hbox{$\boldsymbol{\theta}$}}
\newcommand{\bPi}{\boldsymbol{\Pi}}
\newcommand{\bpi}{\boldsymbol{\pi}}
\newcommand{\boldeta}{\boldsymbol{\eta}}
\newcommand{\bDelta}{\boldsymbol{\Delta}}
\newcommand{\bphi}{\boldsymbol{\phi}}
\newcommand{\bLambda}{\boldsymbol{\Lambda}}

\newcommand{\otheta}{\overline{\theta}}
\newcommand{\diag}{\hbox{\textbf{diag}}}
\newcommand{\trueregdegree}{\hbox{$D_0$}}
\newcommand{\obsregdegree}{\hbox{$\hat{D}$}}
\newcommand{\trueregdegreehalf}{\hbox{$D_0^{-\frac{1}{2}}$}}
\newcommand{\obsregdegreehalf}{\hbox{$\hat{D}^{-\frac{1}{2}}$}}
\newcommand{\maxthetai}{\otheta \vee \theta_i}
\newcommand{\minthetai}{\otheta \wedge \theta_i}
\newcommand{\minthetaa}{\otheta \wedge \theta_a}
\newcommand{\minthetax}{\otheta \wedge \theta_x}
\newcommand{\minthetaab}{\theta_a \wedge \theta_b}

\newtheorem{assumption}{Assumption}
\newtheorem{thm}{Theorem}[section]
\newtheorem{lem}{Lemma}[section]
\newtheorem{cor}{Corollary}[section]
\newtheorem{prop}{Proposition}[section]

\newtheorem{rem}{Remark}

\usepackage{lastpage}

\ShortHeadings{Optimal Parameter Estimation in DCMM Models}{Jiang and Fan}
\firstpageno{1}

\begin{document}

\title{Optimal Parameter Estimation in Degree Corrected Mixed Membership Models}

\author{\name Stephen Jiang \email stephenjiang@princeton.edu \\
       \addr Department of Operations Research and Financial Engineering \\
       Princeton University\\
       Princeton, NJ 08544, USA
       \AND
       \name Jianqing Fan \email jqfan@princeton.edu \\
       \addr Department of Operations Research and Financial Engineering \\
       Princeton University\\
      Princeton, NJ 08544, USA}

\maketitle

\begin{abstract}
With the rise of big data, networks have pervaded many aspects of our daily lives, with applications ranging from the social to natural sciences. Understanding the latent structure of the network is thus an important question. In this paper, we model the network using a Degree-Corrected Mixed Membership (DCMM) model, in which every node $i$ has an affinity parameter $\theta_i$, measuring the degree of connectivity, and an intrinsic membership probability vector $\pi_i = (\pi_1, \cdots \pi_K)$, measuring its belonging to one of $K$ communities, and a probability matrix $P$ that describes the average connectivity between two communities. Our central question is to determine the optimal estimation rates for the probability matrix and degree parameters $P$ and $\Theta$ of the DCMM, an often overlooked question in the literature. By providing new lower bounds, we show that simple extensions of existing estimators in the literature indeed achieve the optimal rate.  Simulations lend further support to our theoretical results.
\end{abstract}


\section{Introduction}
A common problem in network science is to cluster the vertices of a graph, based on the adjacency matrix and other graph structures. For example, advertising teams are often interested in identifying individuals with similar preferences, so as to reach the consumers most likely to purchase their product while minimizing commercial costs \citep{Yang2013}. Likewise, determining similar proteins in protein-protein interaction networks enriches our knowledge of biological systems, allowing labs and pharmaceutical companies to innovate quicker \citep{guimera2005functional, zhang2020detecting}, and financial statement fraud behaviors are affected by peer's culture and pressures and the network based on business similarity provides useful information for defining relevant peers \citep{fan2023unearthing}.

Popular approaches to network clustering encompass both algorithmic-based methods, which attempt to optimize a criterion \citep{Newman_2013, Zhang_2014}, and probabilistic-based methods, which attempt to fit a model to the data \citep{goldenberg2010survey}. Due to its simple formulation and roots in the classic Erd\"os-Renyi graph, one particular model named ``stochastic block models'' (SBM) has garnered significant attention over the past decade {\citep{holland1983stochastic, wang1987stochastic, abbe2023communitydetectionstochasticblock}}. In this period, various extensions have been proposed to the SBM, such as the mixed-membership block model \citep{airoldi2008mixed}, which allows for individuals to belong to several communities, and the degree-corrected block model \citep{karrer2011stochastic}, which allows individuals to have different expected degrees and thus connections in the graph. In this paper, we study the degree-corrected mixed-membership block model (DCMM), which combines the novel aspects of both aforementioned models.

More specifically, consider an undirected graph $G=(V,E)$ with $n$ nodes, where $V=[n]:=\{1,2,\ldots,n\}$ is the set of nodes and $E \subseteq [n]\times[n]$ denotes the set of edges or links between the nodes. Given such a graph $G$, consider its symmetric adjacency matrix $\bX \in \mathbb{R}^{n \times n}$ which captures the connectivity structure of $\bX$, namely $x_{ij}=1$ if there exists a link or edge between the nodes $i$ and $j$, i.e., $(i,j) \in E$ and $x_{ij}=0$ otherwise. For the sake of convenience, we assume that the matrix $\bX$ has self-loops, to avoid the tedious computations inherent to the non self-loop case. However,  the results can be readily extended to the non-self-loop case.

In the degree-corrected mixed membership model (DCMM), we assume the existence of a latent community structure, such that the graph $G$ can be decomposed into $K$ disjoint communities $\{\mathcal{C}_1, \cdots \mathcal{C}_k \}$, where the probability that any arbitrary node $i \in [n]$ belongs to a given community is recorded in a membership vector $\bpi_i \in \mathbb{R}^K$ \citep{jin2022mixed}. In particular,
$$
\mathbb{P}(\hbox{node $i$ belongs to community $\mathcal{C}_k$}) = \pi_i(k) \quad \forall \; k \in [K]
$$
A node is called a ``pure node'' if $\bpi_i(k) = 1$ for some community $k \in [K]$. The probability that any two individuals $i, j$ are connected follows a Bernoulli distribution:
$$
\mathbb{P}(\hbox{edge appearing between individuals $i$ and $j$}) = \theta_i \theta_j \bpi_i^T \bP \bpi_j
$$
Here, $\bP \in \mathbb{R}^{K \times K}$ is a nonnegative matrix modeling the connectivity among distinct communities, and $\btheta = (\theta_1, \cdots \theta_n)$ is the degree parameters of each node, with average degree $\otheta = \frac{1}{n} \sum_{i = 1}^n \theta_i$.  We will impose the unit diagonals of $\bP$ for identifiability.  Defining $\bTheta = \diag(\btheta) \in \mathbb{R}^{N \times N}$ to be the diagonal matrix of degrees, we can express the DCMM in a matrix form,
\begin{align*}
\bH = \bTheta \bPi \bP \bPi^T \bTheta, \quad X_{ij} = \hbox{Bernoulli}(H_{ij}) \quad \forall \; 1 \leq i, j \leq n \numberthis \label{DCMMdef}
\end{align*}

This paper aims to establish the optimal estimation rates for the entries $P_{kl}$ and $\theta_i$ of the probability and degree matrices $\bP$ and $\bTheta$. In the following section, we define and impose several reasonable assumptions that will enable us to construct an estimator that achieves the optimal rate. Our estimators are direct generalizations of those in \citet{jin2022mixed}, which address the setting of low degree heterogeneity, to the severe degree heterogeneity setting, so we impose many of the same assumptions. In Section 3, we formally define and establish the estimation rates of our estimators. More importantly, we construct a novel information lower bound for $\bP$ and $\bTheta$. These lower bounds constitute the main result of our paper; the key idea underlying them is to reparameterize the DCMM in terms of a random dot product graph, another popular family of statistical models for network clustering \citep{athreya2017statisticalinferencerandomdot}. In Section 4, we provide simulation results verifying the rates of our estimators. Finally, we conclude in Section 5 by discussing some implications of our work and future directions.

While many papers in the community detection literature seek to estimate the membership matrix $\bPi$, few have focused on estimating $\bP$ and $\bTheta$, obtaining rates instead as a byproduct of their main results. Indeed, the spectral clustering algorithms SPACL and Mixed-SCORE \citep{mao2021estimating, jin2022mixed} differ in their methodology for estimating $\bPi$. However, their estimators for $\bP$ and $\bTheta$ are identical, attaining dependencies on $n$ and $\bTheta$ of $|\hat{P}_{ab} - P_{ab}|\lesssim \sqrt{\frac{\log(n)}{n \overline{\theta}^2}}$ and $|\hat{\theta}_i - \theta_i| \lesssim \sqrt{\frac{1}{n}}(\sqrt{\frac{\theta_i}\otheta} + \frac{\theta_i}{\otheta})$ in the case where $\btheta$ is nearly homogeneous, i.e. $\frac{\theta_{max}}{\theta_{min}} \leq C$ for some constant $C$. In contrast, the best lower bound known for $|\hat{P}_{ab} - P_{ab}|$ until this point had a dependence of $\frac{1}n$ \citep{markshov2018}; in addition, it assumed that $\theta = \Omega(1)$. Our main contributions in this paper are to therefore i) extend the lower bound for $|\hat{P}_{ab} - P_{ab}|$ to the highly heterogeneous degree setting and improve its dependence on $n$, thus obtaining a rate of $\sqrt{\frac{1}{n\otheta}}$ and ii) derive lower bounds $|\theta_i - \theta|$ for the degree parameters, with dependencies $\sqrt{\frac{1}{n}}(\sqrt{\frac{\theta_i}\otheta} + \frac{\theta_i}{\otheta})$ on $n$ and $\btheta$.
 
\section{Notation and Assumptions}
To motivate our estimators, we first note that \citet{jin2022mixed} provides estimators for the DCMM parameters in the low degree heterogeneity setting using the Mixed-SCORE algorithm; \citet{ke2022optimal} later adapts the algorithm to the severe degree heterogeneity situation through the Mixed-SCORE-Laplacian, but only provides an estimator for $\bPi$, even though generalizing the estimators for $\bP$ and $\bTheta$ in \citet{jin2022mixed} would be straightforward. Thus, we simply formally construct and analyze the generalizations of the estimators, thereby imposing many of the same assumptions.

Given a DCMM model \eqref{DCMMdef}, let $\bD_{0}\in\mathbb{R}^{n\times n}$ be a positive diagonal matrix with $\bD_{0}(i,i)= (\be_i+\frac{1}{n}\mathbf{1}_n)^T \bH \mathbf{1}_n$ for all $1\leq i\leq n$. Define
\[
\bG:= K \cdot \bPi^T\bTheta \bD_{0}^{-1}\bTheta \bPi \quad \in\mathbb{R}^{K\times K}. 
\]
For all $k \in [K]$, let $\lambda_k(\bP\bG)$ be the $k$th largest eigenvalue in magnitude of $\bP\bG$, and denote its first right eigenvector by $\boldeta_1 \in \mathbb{R}^K$. As in Ke and Wang (2022), we impose the following conditions.

\begin{assumption} \label{regconds1} (Assumption (3.1) - (3.3) of \citet{ke2022optimal}) The DCMM parameters $(\bTheta,\bPi,\bP)$ satisfy the following requirements:
\begin{itemize} \itemsep -5pt
\item[(a)] $\Vert \bG\Vert\leq c_1$, $\Vert \bG^{-1}\Vert\leq c_1$, and $\min_{1\leq k\leq K}\bigl\{\sum_{i=1}^n\theta_i\bpi_i(k)\bigr\}\geq c_1\|\btheta\|_1$. 
\item[(b)] $\max_{k\neq 1}\{\lambda_k(\bP\bG)\}\leq \min\{(1-c_2)\cdot \lambda_1(\bP\bG), c_2^{-1}\sqrt{K}\}$.
\item[(c)] If $\boldeta_1 > 0$ is the leading right eigenvector of $\bP\bG$, it satisfies $\min_{1\leq k\leq K} \boldeta_1(k)>0$ and $\min_{1\leq k\leq K}\{\boldeta_1(k)\}\geq c_3\max_{1\leq k\leq K}\{\boldeta_1(k)\}$. 
\item[(d)] $P$ is non-singular and has unit diagonals.
\item[(e)] Each community $\mathcal{C}_k$ contains at least one pure node, i.e. for all $k \in [K]$, there exists an $i \in [n]$ for $\bpi_i = \mathbf{e_k}$.
\end{itemize} 
\end{assumption}
As elaborated in \citet{ke2022optimal}, $\bG$ is a measure of the balance among communities. Thus, Assumption \ref{regconds1}(a) ensures the communities are well-balanced. Assumption \ref{regconds1}(b) is a mild eigengap condition. Assumption \ref{regconds1}(c) may appear difficult to satisfy, but actually holds for a variety of settings (see \citet{jin2022mixed}).

Now, we introduce the additional assumptions imposed by this paper. For clarity's sake, however, we first establish some relevant notation. Let $\obsregdegree(i,i)= (\be_i+\frac{1}{n}\mathbf{1}_n)^T \bX \mathbf{1}_n \quad \forall \; i \in [n]$. Define the \textit{regularized graph Laplacian matrices} $\bL_0$ and $\hat{\bL}$ by
\begin{align*}
    \bL_0 & := \boldsymbol{\trueregdegreehalf} \bH \boldsymbol{\trueregdegreehalf} \\
    \hat{\bL} & := \boldsymbol{\obsregdegreehalf} \bX \boldsymbol{\obsregdegreehalf}
\end{align*}
For $1\leq k\leq K$, let $\lambda_k$ be the kth largest eigenvalue in magnitude of $\bL_0$, with corresponding eigenvector ${\xi}_k$; define $\hat{\lambda}_k$ and $\hat{\xi}_k$ analogously for $\hat{\bL}$. Denote 
\begin{align*}
    \bLambda & = \diag(\lambda_1, \cdots, \lambda_K), \qquad \wh \bLambda = \diag(\wh \lambda_1, \cdots, \wh \lambda_K) \in \mathbb{R}^{K \times K} 
\end{align*} 
and for $k \in [K]$, define 
\begin{align*}
    \bLambda_{-k} & = \diag(\lambda_1, \cdots, \lambda_{k - 1}, \lambda_{k + 1}, \cdots, \lambda_K),  \qquad \wh \bLambda_{-k} = \diag(\wh \lambda_1, \cdots, \wh \lambda_{k - 1}, \wh \lambda_{k + 1}, \cdots, \wh \lambda_K) \in \mathbb{R}^{(K - 1) \times (K - 1)}
\end{align*}
i.e. the eigenvalues excluding the $k$th one. Similarly, denote the corresponding eigenvectors by 
\begin{align*}
    \bXi =[\xi_1,\ldots,\xi_K], \qquad \hat{\bXi} =[\hat{\xi}_1,\ldots,\hat{\xi}_K] \in \mathbb{R}^{n \times K}
\end{align*} and for $k \leq K$, define 
\begin{align*}
    \bXi_{-k} & = [\xi_1,\ldots,\xi_K],  \qquad \hat{\bXi}_{-k}=[\hat{\xi}_1,\ldots,\hat{\xi}_K] \in \mathbb{R}^{n \times (K - 1)}
\end{align*}
to be the set of eigenvectors of $\bL_0$ and $\hat{\bL}$ respectively, excluding the $k$th eigenvector. Lastly, let
$$
F:= K \norm{\btheta}^{-2} \cdot \bPi^T\bTheta^2 \bPi \quad \in\mathbb{R}^{K\times K}
$$
\begin{assumption} \label{regconds2} ~
\begin{itemize} \itemsep -5pt
    \item[(a)] There exists a constant $c_4 > 0$ such that for all $k \in [K]$, $\max_{l \neq k}|\lambda_k(\bL_0) - \lambda_l(\bL_0)| \geq c_4 \cdot |\lambda_k(\bL_0)|$.
    \item[(b)] $\lambda_K(\bL_0) \gg \sqrt{\frac{\log n}{n \otheta^2}}$
    \item[(c)] $\norm{F} \leq c_5, \norm{F^{-1}} \leq c_6$.
    \item[(d)] For all $i \in [n]$, $c_7\sqrt{\frac{\log n}{n}} \ll \theta_i \leq c_8$.
    \item[(e)] $\norm{\bP}_{\text{max}} \leq c_9$
    \item[(f)] The number of communities $K$ is fixed.
\end{itemize}
\end{assumption}
Assumption \ref{regconds2}(a) ensures that the eigenvalues of the Laplacian are well-separated. We make this assumption to obtain better deviation bounds on the eigenvectors; more specifically, if certain conditions on $\lambda_K$ and $\otheta$ hold (to be defined in Theorem \ref{estimationP2}), we can obtain better bounds than those in \citet{ke2022optimal}. We acknowledge that in certain models, results allowing a multiplicity of eigenvalues have been established. Since $\norm{\hat{\bL} - \bL_0} \lesssim \sqrt{\frac{\log n}{n \otheta^2}}$ with high probability, Assumption \ref{regconds2}(b) is mild; it ensures the minimum signal is always at least as large as the noise. As elaborated in Assumption 2 of \citet{jin2022mixed}, $F$ is a measure of the balance among communities; thus Assumption \ref{regconds2}(c) is not difficult to satisfy in practice. Assumption \ref{regconds2}(d) and Assumption \ref{regconds2}(e) ensure the model is not too heterogeneous. Assumption \ref{regconds2} is a mild condition imposed to simplify the analysis; if desired, the constant order of $K$ in Assumption \ref{regconds2}(f) can be replaced with growth order $K \asymp o(n^{\alpha})$ at the cost of significantly more tedious computations.

Lastly, we introduce some auxillary notation that will be used throughout this paper. For a matrix $\bA = (A_{ij}) \in \mathbb{R}^{m \times n}$ and any $i \in [m], j \in [n]$, denote the $i$th row of entries by $\bA_{i, :} = [A_{i1}, A_{i2}, \cdots, A_{in}] \in \mathbb{R}^{1 \times n}$, and the $j$th column by $\bA_{:, j} = [A_{1j}, A_{2j}, \cdots, A_{mj}]^T \in \mathbb{R}^{m \times 1}$. Furthermore, we use $a(k)$ to denote the kth component of a vector $\ba$. For each $k \in [K]$, denote $\mathcal{C}_k = \{i \in [n] \verts \bpi_i = \mathbf{e_k}\}$ to be the set of pure nodes in the kth community. Also, for any two sequences $a_n$ and $b_n$, we write $a_n \lesssim b_n$ if there exists a constant $C$ such that $a_n \leq Cb_n$ for any $n$. Denote $a \vee b = \max\{a, b\}$ and $a \wedge b = \min\{a, b\}$.

\section{Main Results} 
In this section, we identify the optimal estimation rates of the entries of the $\bP$ and $\bTheta$ under Assumptions 1 and 2. The main result of this section will be new minimax lower bounds for the estimation of $\bP$ and $\bTheta$. 

\subsection{Estimators for $\bP$ and $\bTheta$}
First, we extend the estimators in \citet{jin2022mixed} to the regularized Laplacian setting. Below, we describe our estimation procedure. Our algorithm essentially follows the normalization and vertex hunting steps of the Mixed-SCORE-Laplacian algorithm \citep{ke2022optimal}. However, there are a couple differences. First, \citet{ke2022optimal} employs a vanilla version of the successive projection algorithm (SPA) for vertex hunting \citep{araujo2001successive}, whereas we employ the vertex hunting algorithm in \cite{bhattacharya2023inferences}, which is a slightly modified version of SPA. Secondly, \citet{ke2022optimal} excludes any node whose degree is below a certain threshold before applying successive projection, whereas we do not prefilter any nodes before applying our vertex hunting algorithm.

\begin{figure}
\begin{algorithm}[H]
\caption{Estimation of $\bP$ and $\bTheta$}
	\begin{algorithmic}[1]		
        \State \textbf{Input} observed Laplacian matrix $\wh D$ and a constant $\phi$.
        \State Compute the row-normalized eigenvectors
    $$
    \wh \br_i := \left[\frac{\wh \xi_2(i)}{\wh \xi_1(i)},\frac{\wh \xi_3(i)}{\wh \xi_1(i)},\dots, \frac{\wh \xi_K(i)}{\wh \xi_1(i)}\right]^T \in \mathbb{R}^{K-1},\quad\forall \; i\in [n];
    $$ 
        \State Input $\wh \br_1, \wh \br_2, \cdots \wh \br_n$ and a radius $\phi$ into a vertex hunting algorithm \ref{algVH}. Denote the estimated vertices it outputs by $\wh \bv_k \; \forall \; k \in [K]$.
    \State Define $$\wh \bQ = \begin{bmatrix}
    1 & \wh\bv_1^T \\
    1 & \wh\bv_2^T \\
    \vdots & \vdots \\
    1 & \wh \bv_K^T
    \end{bmatrix} \in \mathbb{R}^{K \times K}$$ and compute a vector $\wh b_1(k) = (\wh \lambda_1 + \wh \bv_k^T \bLambda \bv_k)^{-\frac{1}{2}} \quad \forall \; k \in [K]$, where $\wh \bb_1 \in \mathbb{R}^{K}$.
    \State For each $i\in [n]$, solve $\hat{\bw}_i\in\mathbb{R}^K$ from the linear equations: $\sum_{k=1}^K \hat{w}_i(k)\hat{\bv}_k=\hat{\br}_i$ and $\sum_{k=1}^K \hat{w}_i(k)=1$. Compute the estimated memberships $\hat{\bpi}_i^*(k)=\max\{\hat{w}_i(k)/\hat{b}_1(k),\, 0\}, \;  \hat{\bpi}_i=\hat{\bpi}_i^*/\|\hat{\bpi}_i^*\|_1  \quad \forall \; i \in [n]$.
    \State Calculate $\wh \bP = \wh \bb_1^T \wh \bQ \wh \bLambda \wh \bQ^T \wh \bb_1$ and $\hat{\bTheta}(i, i) = \hat{\xi}_1(i)\trueregdegree(i, i)^{\frac{1}2}(\hat{\bpi}_i^T\hat{\bb}_1)^{-1}$.
    \State \Return $\wh \bP$ and $\hat{\bTheta}$.
	\end{algorithmic}  \label{algestimationP}
\end{algorithm}

\vspace{\floatsep}

\begin{algorithm}[H]
\caption{Modified Version of Successive Projection (\cite{bhattacharya2023inferences}; Algorithm 1)}
	\begin{algorithmic}[1]		
        \State \textbf{Input} $\wh \br_1, \wh \br_2, \cdots, \wh \br_n$ and a radius $\phi > 0$.
        \State \textbf{Initialize} $Z_i = [1, \wh \br_i^T]^T$, for $i \in [n]$. 
        \State \textbf{for $k \in [K]$ do}
     \State \hspace{\algorithmicindent} Let $i_k = \argmax_{1\leq i\leq n}\left\|\bZ_i\right\|_2$ and $\wh\bv_k'^T = \wh \br_{i_k}$  
        \State \hspace{\algorithmicindent} Update $\bZ_i \leftarrow \bZ_i - \wh \br_{i_k}\wh \br_{i_k}^\top\bZ_i /\left\|\wh \br_{i_k}\right\|_2^2$, for $i\in [n]$
        \State \textbf{end for}
        \State Let $$\hat{\mathcal{C}}_k = \left\{i\in [n]:\left\|\wh\br_i - \wh\bv_k'^T\right\|_2\leq \phi\right\} \text{ and } \wh\bv_k = \frac{1}{|\hat{\mathcal{C}}_k|}\sum_{i\in \hat{\mathcal{C}}_k}\wh\br_i$$
        for $k\in [K]$
        \State \textbf{return} $\hat{\mathcal{C}}_1, \hat{\mathcal{C}}_2, \dots, \hat{\mathcal{C}}_K$ and $\wh\bv_1,\wh\bv_2,\dots, \wh\bv_K$
	\end{algorithmic}  \label{algVH}
\end{algorithm}
\end{figure}

To motivate our choice of estimators in Algorithm \ref{algestimationP}, we note the following proposition. First, we define some relevant notation. \begin{itemize}
\item Denote $$
    \br_i := \left[\frac{\xi_2(i)}{\xi_1(i)},\frac{\xi_3(i)}{\xi_1(i)},\dots, \frac{\xi_K(i)}{\xi_1(i)}\right]^T \in \mathbb{R}^{K-1},\quad\forall \; i\in [n].
    $$ Likewise, for each $k \in [K]$, define the true vertices to be $\bv_k := \br_i$ for any $i \in \mathcal{C}_k$.
    \item Define $Q = \begin{bmatrix}
        1 & \bv_1^T \\
        1 & \bv_2^T \\
        \vdots & \vdots \\
        1 & \bv_K^T
    \end{bmatrix} \in \mathbb{R}^{K \times K}$.
    \item Denote the $1 \times K$ vector $b_1(k) = ( \lambda_1 +  \bv_k^T \bLambda \bv_k)^{-\frac{1}{2}} \quad \forall \; k \in [K]$.
\end{itemize}
\begin{prop} \label{recoverP}
       The relations 
       \begin{align*}
           \bP & = \diag(\bb_1)Q\bLambda Q^T\diag(\bb_1) \\
           \theta(i, i) & = \xi_{1}(i)\trueregdegree(i, i)^{\frac{1}2}(\bpi_i^T\bb_1)^{-1} \quad \forall \; i \in [n]
       \end{align*}
       hold, i.e. the true counterparts of the estimators $\wh \bP$ and $\hat{\bTheta}$ defined in Algorithm \ref{algestimationP} are the plug-in estimators of the true parameter values.
   \end{prop}
   We now establish the estimation rate of $\wh \bP$ under the most general conditions. For all $k \in [K]$, define $\theta_{k, \text{min}} = \min_{i \in \mathcal{C}_k} \theta_i$, i.e. the smallest degree among all pure nodes in community k. Also, let
\begin{align*}
\nu_n = \min\bigl\{K^{-1/2}\lambda_1, \;  \lambda_K\bigr\} \numberthis \label{delta}
\end{align*}
\begin{thm}\label{estimationP1}
   Let $\hat{\bP}$ be the estimator as defined in Algorithm \ref{algestimationP}. Under Assumptions \ref{regconds1} and \ref{regconds2}, for any fixed $a, b \in [K]$, there exists a permutation $T$ on $\{1, 2, \cdot K\}$ such that 
$$
 |\hat{P}_{t(a)t(b)} - P_{ab}| \lesssim \sqrt{\frac{K^4\lambda_1^2\lambda_2^2}{n\otheta(\otheta \wedge (\minthetaab))\nu_n^2}} + \sqrt{\frac{K^3 \lambda^2_1 \log n}{n \otheta^2 \lambda_K^2}}
$$
with probability $1 - o(n^{-10})$.
\end{thm}
In Theorem 3.2 below, however, we establish a lower bound that is independent of the smallest eigenvalue $\lambda_K$. Since our current estimator $\hat{P}$ depends on $\lambda_K$, it may be suboptimal. In order to achieve the optimal dependence, we additionally impose Assumption \ref{regconds3}.

\vspace{2 mm}

\begin{assumption} \label{regconds3} ~
\begin{itemize}
    \item $\otheta \gg \frac{1}{\lambda_K}\sqrt[4]{\frac{\log^2 n}{n}}, \sqrt[4]{\frac{K\log^2 n}{n^2\lambda_K^2}}$
    \item $\frac{\theta_{\text{max}}}{\otheta} \ll \frac{n\otheta^2\lambda_K^2}{K \log n}$
    \item $\lambda_K \gg \sqrt{\frac{K \theta_{\text{max}} \log n}{n\otheta^3}}$
    \item $\min_{i \in (\mathcal{C}_a \cup \mathcal{C}_b)} \theta_i \geq \otheta$
\end{itemize}
\end{assumption}

\vspace{2 mm}

\begin{thm}\label{estimationP2}
   If Assumption \ref{regconds3} holds in addition to the assumptions of Theorem \ref{estimationP1}, the estimation rate of $\hat{\bP}$ can be improved to
$$
 |\hat{P}_{t(a)t(b)} - P_{ab}| \lesssim \sqrt{\frac{K^4 \log n}{n \lambda_K \min_{i \in (\mathcal{C}_a \cup \mathcal{C}_b)}(\theta_i(\minthetai))}} + \sqrt{\frac{K^2\log n}{n \otheta^2}} + \sqrt{\frac{K^4 \lambda_1^2 \log n}{n\otheta^2}}
$$
for any fixed $a, b \in [k]$ and sufficiently large $n$, with probability $1 - o(n^{-10})$.
\end{thm}

The key assumption enabling this bound is the well-separated eigenvalue condition. Indeed, recent literature such as \citet{jin2022mixed} and \citet{ke2022optimal} only assume that the first eigenvalue is well-separated from the others and, as a result, they are only able to obtain $l_{\infty}$ deviation bounds $|\hat{\xi}_1(i) - \xi_1(i)|$ for the first eigenvector. In contrast, assuming all eigenvalues are well-separated enables us to obtain $l_{\infty}$ deviation bounds $|\hat{\xi}_k(i) - \xi_k(i)|$ for all $k \in [k]$ (and, by extension, $l_{\infty}$ deviation bounds on the normalized eigenvector rows $|\hat{\br}_{k - 1}(i) - r_{k - 1}(i)|$).

\begin{rem}
That being said, the bound in Theorem \ref{estimationP2} actually fails to improve upon our bound in Theorem \ref{estimationP1} if $\min_{i \in (\mathcal{C}_a \cup \mathcal{C}_b)} \theta_i \leq \otheta$; that is, there exists some pure node among communities $a, b$ with degree less than the mean degree. Extending the bounds to the case of low-degree pure nodes is an interesting topic for future research.
\end{rem}

Regarding the estimation rate of $\bTheta$, we first impose an additional assumption.

\vspace{2 mm}

\begin{assumption} \textit{(Assumption 3.4 of Ke and Wang, 2022) \label{regconds4}
    There exists a constant $c$ such that $$
    \bigl\{ 1\leq i\leq n: \bpi_i(k)=1, \;\; \theta_i\geq c\bar{\theta}\bigr\}\neq \emptyset, \qquad \mbox{for each }1\leq k\leq K. 
    $$
i.e. each community has a pure node with degree at least a multiple of the mean.}
\end{assumption}

\vspace{2 mm}

We impose this condition because unlike $\hat{P}_{ab}$, $\hat{\theta}_i$ uses an estimate of the membership vector $\bpi_i$ of the ith node. In particular, our estimate of $\bpi_i$ involves all $K$ vertices $\bv_k$, as opposed to just the $a$th and $b$th vertices $\bv_a, \bv_b$; consequently, its error depends on the estimation errors of the vertices. Assumption \ref{regconds4} thus ensures that these vertex estimation errors do not dominate the error $|\hat{\br}_i - \br_i|$ of the node itself.

\begin{thm} \label{estimationTheta}
    Under Assumptions \ref{regconds1} and \ref{regconds2}, and the additional Assumption \ref{regconds4}, there exists a constant $c$ such that
    \begin{align*}
        |\hat{\bTheta}(i, i) - \bTheta(i, i)| & \leq \sqrt{\frac{K\theta_i^2\log n}{n\bar{\theta} (\bar{\theta}\wedge \theta_i)\nu_n^2}} + \sqrt{\frac{K\theta_i^2\log n}{n\otheta^2 \lambda_K^2}}
    \end{align*}
    for all $i \in [n]$ with probability $1 - o(n^{-10})$.
\end{thm}

\subsection{Lower Bounds for $\bP$ and $\bTheta$}

Now, we provide lower bounds for $\bP$ and $\bTheta$, which match up to factors of $K$ and a logarithmic factor. As is standard in minimax analysis, we first define the family of models for which we will establish a lower bound. Given a vector of scalars $\sigma = (n, K, \lambda_1, \lambda_K, \otheta, \tilde{\theta})$, let ${\cal Q}_{n, ab}(\sigma)$ denote the collection of $K$-community DCMM models $(\bTheta, \bPi, \bP)$ satisfying Assumptions \ref{regconds1} and \ref{regconds2} for which:
\begin{align}
    \lambda_1(L) \geq \lambda_1, \lambda_K(L) \geq \lambda_K, \sum_{i = 1}^n \Theta(i, i) = n\otheta, \tilde{\theta} \geq \min_{i \in \mathcal{C}_a \cup \mathcal{C}_b} \theta_i \label{lowerboundsPconds}
\end{align}

\begin{thm} \label{lowerboundsP}
For any vector of scalars $\sigma$, pairs $a$ and $b$, and constant $c_{10} > 0$, define $\tilde {\cal Q}_{n, ab, c_{10}}(\sigma) = {\cal Q}_{n, ab}(n, K, c_{10}\lambda_1, c_{10}\lambda_K, c_{10}\otheta, c_{10}\tilde{\theta})$, i.e. models with eigenvalues and degrees that are $c_{10}$ times those of ${\cal Q}_{n, ab}(\sigma)$. Then, there exist sufficiently small constants $C, C' > 0$ such that, for all fixed $a, b \in [K]$ and sufficiently large $n$, 
$$
\inf_{\hat{\bP}}\sup_{\bP \in \tilde {\cal Q}_{n, ab, C}(\sigma)} |\hat{P}_{ab} - P_{ab}| \geq \frac{C' \lambda_1}{\sqrt{n\otheta \tilde{\theta}}}
$$

up to factors of $K$.
\end{thm}
We define ${\cal Q}_{n, ab}$ in this manner because the estimation rate of $P_{ab}$ in Theorem \ref{estimationP2} depends solely on the parameters in $\sigma$. When $\lambda_1$ and $\lambda_K$ are of the same order, the upper bound in Theorem \ref{estimationP1} matches our lower bound, up to factors of $K$ and a logarithm in $n$. Most importantly, though, when $\lambda_1$ and $\lambda_K$ are not of the same order, Theorem \ref{estimationP2} enables us to still obtain a matching upper bound, if Assumption \ref{regconds3} holds and $\tilde{\theta} \leq \otheta$. We leave the question of when either Assumption \ref{regconds3} does not hold or the degree of the vertices $\tilde{\theta}$ exceeds $\otheta$ to future research; we suspect that the eigenvector bounds developed in Theorem \ref{improveddevbounds} will need to be improved, perhaps via new techniques.

Next, we provide a lower bound for $\bTheta$. 
\begin{thm} \label{lowerboundsTheta}
  Given a vector of scalars $\tau = (n, K, \otheta, \theta_i)$, such that $\otheta, \theta_i$ satisfy Assumption \ref{regconds2}(d), let ${\cal R}_{n, i}(\tau)$ denote the collection of $K$-community DCMM models $(\bTheta, \bPi, \bP)$ satisfying Assumptions \ref{regconds1}, \ref{regconds2}, and Assumption \ref{regconds4} for which:
  \begin{align*}
      \sum_{i = 1}^n \Theta(i, i) = n\otheta \hbox{ and the ith node has degree $\theta_i$} \numberthis \label{lowerboundsThetaconds}
  \end{align*}
  Furthermore,  define $\tilde {\cal R}_{n, ab, c_{11}}(\sigma) = {\cal R}_{n, ab}(n, K, c_{11}\otheta, c_{11}\theta_i)$ for any vector $\tau$ and constant $c_{11} > 0$, i.e. models with eigenvalues and degrees $c_{11}$ times those of ${\cal R}_{n, i}(\tau)$. Then, there exist sufficiently small constants $c_{11}, C$ such that for all sufficiently large $n$, 
\begin{equation}
    \label{LB-rateTheta}
\inf_{\hat{\bTheta}}\sup_{\bTheta \in {\cal R}_{n, i, c_{11}}(\tau)} |\hat{\Theta}(i, i) - \Theta(i, i)| \geq C\sqrt{\frac{\log n}{n\otheta^2}}(\sqrt{\frac{\theta_i}{\otheta}} + \frac{\theta_i}{\otheta})
\end{equation}
\end{thm}
Again, we define ${\cal R}_{n, i}(\tau)$ in this manner because the estimation rate of $\theta_i$ in Theorem \ref{estimationTheta} depends solely on the parameters in $\tau$. When the condition number is small and Assumption \ref{regconds4} holds, the upper bound in \ref{estimationTheta} matches our lower bound, up to factors of $K$ and a logarithm in $n$. We leave the question of improving the lower bound when either Assumption \ref{regconds4} does not hold or $\frac{\lambda_1}{\lambda_K} \gg 1$ to future research.

\section{Simulations}
\begin{figure}
    \centering
    \includegraphics[scale=0.9]{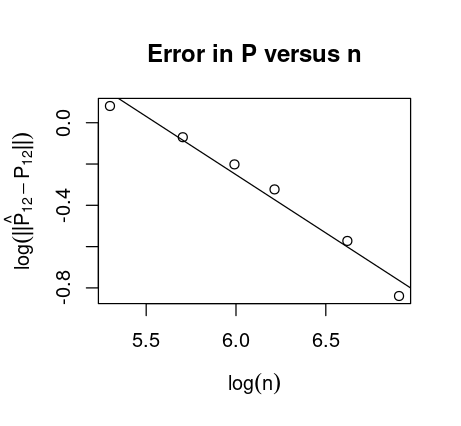}
            \label{fig:SRl}
    \caption{Convergence rate of the error $|\hat{P}_{12} - P_{12}|$ for our algorithm.}\label{fig1}
\end{figure}

We conduct two experiments, testing the theoretical rates of our estimators for $P$ and $\Theta$ respectively. However, we do not use Algorithm \ref{algestimationP} in our experiments. Instead, we use a slightly modified version of Algorithm \ref{algestimationP}; specifically, we use a different vertex hunting algorithm called ``sketched vertex search" (SVS), which is defined in \cite{jin2022mixed}. The reason for this change is that the successive projection algorithm is highly sensitive to node-wise errors. Since our estimator for $\bP$ depends on the individual vertex estimation errors, and $\btheta_i$ on the node-wise embedding errors $|\hat{\br}_i - \br_i|$, it is better to use a more robust vertex hunting algorithm. SVS accomplishes this by first applying a k-means `denoising' step, before applying the successive projection algorithm. The exact details of SVS are listed in Algorithm \ref{algVHmodifiedSVS} below. For more details on the SVS and other vertex hunting algorithms, we refer the reader to \cite{ke2022optimal} and \cite{ke2022scorenormalizationespeciallyhighly}.

\begin{algorithm}[t]
\caption{Modified Version of Sketched Vertex Search}
	\begin{algorithmic}[1]		
        \State \textbf{Input} $\wh \br_1, \wh \br_2, \cdots, \wh \br_n$, a radius $\phi > 0$, and a tuning integer $L \geq K$.
        \State \textbf{Initialize} $Z_i = [1, \wh \br_i^T]^T$, for $i \in [n]$. 
        \State Run k-means clustering on $Z_1, Z_2, \cdots , Z_n$, with $\bL_0$ clusters. Denote the outputs by $x_1, x_2, \cdots, x_L$.
        \State Input $x_1, x_2, \cdots, x_L$ and $\phi$ into Algorithm \ref{algVH}, outputting $\hat{\bv}_1, \hat{\bv}_2, \cdots, \hat{\bv}_K$.
        \State \textbf{return} $\wh\bb_1,\wh\bb_2,\dots, \wh\bb_K$
	\end{algorithmic}  \label{algVHmodifiedSVS}
\end{algorithm}

In the first experiment, we fix $K = 2$ and our probability matrix $P = 0.5I_2 + 0.5\mathbf{1}\mathbf{1}^T$. For each $n \in \{200, 300, 400, 500, 750, 1000\}$, we generate 100 random pairs of parameters $(\bTheta, \bPi)$: to generate $\bTheta$, we draw $\theta_1, \theta_2, \cdots \theta_n \stackrel{iid}{\sim} \hbox{Uniform}([0.05, 0.8])$, and to generate $\bPi$, we set $\bpi_i = [1, 0]^T$ and $\bpi_i = [0, 1]^T$ each for 10$\%$ of the nodes, and draw $\bpi_i = [t_i, 1 - t_i]^T$, where $t_i \stackrel{iid}{\sim} \hbox{Uniform}([0.15, 0.85])$. For each random pair of parameters, we generate a random realization of the adjacency matrix $A$. The left plot in Figure \ref{fig1} depicts the log of the average of the errors $|\hat{P}_{12} - P_{12}|$ of each of the 100 realizations $(\bTheta, \bPi)$, for each $n$. Since $\theta_i = \Omega(1)$, Theorem \ref{estimationP1} predicts the error rate to have dependence $n^{-\frac{1}{2}}$ on $n$. In our plot, the line of best fit has slope -0.56, therefore nearly matching the theoretical rate.

\begin{figure}[!h]     
            \includegraphics[scale = 0.5, width=\textwidth]{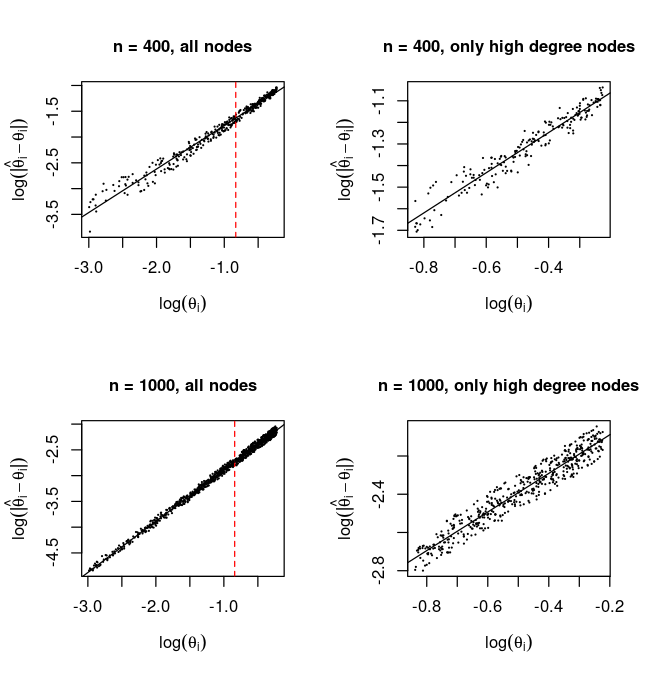}
    \caption{$\log(|\hat{\theta}_i - \theta_i|)$ versus $\log(\theta_i)$ in Experiment 2. Slope of the lines of best fit (from top-left and proceeding clockwise): 0.84, 0.93, 1.01, 0.99. The red dashed lines in the leftmost plots indicate the average degree $\otheta$.}
    \label{fig2} 
\end{figure}

In the second experiment, we seek to investigate the dependence of the degree error $|\hat{\theta}_i - \theta_i|$ versus the degree $\theta_i$; for this reason, we consider a single DCMM parameter set $S = (\bTheta, \bPi, \bP)$, thereby removing any possible confounding effects that may arise from varying $\Pi$ between realizations. In more detail, we first fix the number of datapoints $n$ at a value in the set $\{200, 300, 400, 500, 750, 1000\}$. Then, we fix $K = 2$ and our probability matrix $P = 0.5I_2 + 0.5\mathbf{1}\mathbf{1}^T$. Then, we generate $(\bTheta, \bPi)$ as in the first experiment. Fixing this pair $(\bTheta, \bPi)$, we thus obtain a fixed parameter set $S = (\bTheta, \bPi, \bP)$. Subsequently, we generate 100 realizations of $\mathcal{Q}$ and compute the average of $|\hat{\theta}_i - \theta_i|_1$ for each of these 100 realizations, for all $1 \leq i \leq n$. For the sake of space, we only plot the results for $n = 400, 1000$. Our theory predicts $|\hat{\theta}_i - \theta_i|$ to have a dependence of $\theta_i^1$ for $\theta_i > \otheta$ and a dependence of $\theta_i^{\frac{1}{2}}$ for $\theta_i < \otheta$. We only plot the error for ``high-degree" nodes, i.e. those satisfying $\theta_i > \otheta$, in the two rightmost plots of \ref{fig2}; the lines of best fit in the top-right and bottom-right plots of Figure \ref{fig2} have slopes 0.93 and 1.01 respectively, thus verifying our estimation rate. In the two leftmost plots, we plot all nodes. Since the theoretical rate varies for $\theta_i$ smaller and larger than $\otheta$, we expect a phase transition to occur. However, the top-left and bottom-left plots are nearly linear, with lines of best fit with slopes 0.84 and 0.99 respectively. We conjecture that low-degree nodes may suffer from additional noise, making the phase transition unobservable in our plots.

\section{Discussion}
In this paper, we constructed novel lower bounds for the $\bP$ and $\bTheta$ matrices in the DCMM model, alongside estimators. As a result, we showed that simple extensions of existing estimators can achieve optimal estimation rates for $\bP$ and $\bTheta$, and, additionally, if further assumptions hold on the well-separateness of the eigenvalues and average degree $\otheta$ (i.e. Assumption \ref{regconds3}), we can show a matching upper bound for $P_{ab}$ in the settings where the SNR is low and the minimum degree of a pure vertex $\min_{\theta_i \in (\mathcal{C}_a \cup \mathcal{C}_b)} \theta_i$ is smaller than the average degree $\otheta$. Attaining a optimal rate for $\bP$ in these regimes without Assumption \ref{regconds3} is an interesting question. We conjecture that an entirely new estimator for $\bP$ may be required, as improving the current one would require better entrywise bounds on the eigenvectors and first singular vector $\boldeta_1$ of $\bP\bG$, both challenges that have not been addressed in the literature. Likewise, establishing the optimal rate for $\bTheta$ in the low-SNR regime would be a natural extension. We leave these compelling directions to future research and efforts. 


\acks{The authors do not have competing interests for this work.  The research is in part supported by the NSF grants  DMS-2053832, DMS-2210833 and
ONR N00014-22-1-2340.}

\newpage
\appendix

\section{Basic Results} \label{basicresultssection}
In this section, we state and prove various basic results, which will be frequently referred to in the proofs. First, we characterize the order of the eigenvector entries of $\bL_0$. For the ease of notation, we adapt the ``incoherence" framework, commonly used in denoising problems [e.g. \citep{yan2024inferenceheteroskedasticpcamissing}]. 
\begin{lem} \label{incoherence}
    There exist constants $C, C'$ such that
    \begin{align*}
        \sqrt{\frac{C'\theta_i^2}{n\otheta(\maxthetai)}} \leq \xi_1(i) \leq \sqrt{\frac{C\theta_i^2}{n\otheta(\maxthetai)}}
    \end{align*}
    holds simultaneously for all $i \in [n]$. Define
    \begin{align*}
        \mu_i = \frac{C\theta_i^2}{\otheta(\maxthetai)}
    \end{align*} Then $\xi_1(i) \leq \sqrt{\frac{\mu_i}{n}}$.
\end{lem}
Next, we borrow some fundamental lemmas from \citet{ke2022optimal}.
\begin{lem} \label{lem:order-eigenval} (Extension of Lemma B.1 of \cite{ke2022optimal})
Under Assumptions \ref{regconds1}(a) and \ref{regconds2}(c),  
\begin{align*} 
\lambda_1 \asymp K^{-1}\lambda_1(\bP\bG) > 0, \quad |\lambda_K|\asymp K^{-1}\lambda_K(\bP\bG), \quad \max_{l \neq k} |\lambda_k - \lambda_l|\geq c \lambda_k 
\end{align*}
\end{lem}

\begin{lem} \label{lem:order} (Extension of Lemma B.2 of \cite{ke2022optimal})
Denote $S_1 = \{1\leq i\leq n: \theta_i\geq \bar{\theta}\}, S_2 = \{1\leq i\leq n: \theta_i<\bar{\theta}\}$. There exists a fixed constant $C$ such that 
\begin{align} \label{order-xi}
\xi_1(i)\asymp \sqrt{\frac{\mu_i}{n}} \asymp \frac{1}{\sqrt{n}}\begin{cases}
\sqrt{\theta_i/\bar{\theta}}, & i\in S_1,\\
\theta_i/\bar{\theta}, & i\in S_2,
\end{cases}\qquad 
 \|\bXi_{-1}(i)\|  \leq \sqrt{\frac{K\mu_i}{n}} \leq \frac{C\sqrt K}{\sqrt{n}} 
\begin{cases}
\sqrt{\theta_i/\bar{\theta}}, & i\in S_1,\\
\theta_i/\bar{\theta}, & i\in S_2,
\end{cases}
\end{align}
\end{lem}
Next, we show several results concerning the $\bH$ and $\bW$ matrices, alongside deducing upper bounds on the degree expansion terms defined in the subsection \ref{sectionexpansiondegrees} (which occurs later on in this paper). Define $\mathcal{A}_1$ to be the event that Lemma \ref{basicresults} holds. 
\\[1.5 pt]
\begin{lem}\label{basicresults}
There exist fixed constants $C, C'$ such that for any $i \in [n]$, the following holds with probability $1 - o(n^{-11})$,
\begin{align*}
   Cn\otheta\theta_i \leq \sum_{j = 1}^n H_{ij} & \leq C'n\otheta\theta_i \numberthis \label{basicresultseq1}\\
    D_0(i, i), \hat{D}(i, i) & \geq C'\max \{n\otheta\theta_i, n\otheta^2\}\numberthis \label{basicresultseq2}\\
    \norm{\bH_{i, \cdot}} & \leq C\otheta\theta_i\sqrt{n} \numberthis \label{basicresultseq3}\\
    |[\hat{D} - D_0](i, i)| & \leq 
    \numberthis \label{basicresultseq8} \\
     |\sum_{j = 1}^n W_{ij}| & \leq C\sqrt{n\otheta\theta_i \log n} \numberthis \label{basicresultseq4}\\
     |\sum_{i = 1}^n \sum_{j = 1}^n W_{ij}| & \leq Cn\otheta\sqrt{\log n} \numberthis \label{basicresultseq6} \\
     \norm{\bW_{i, \cdot}} & \leq C\sqrt{n\otheta\theta_i \log n} \numberthis \label{basicresultseq5} \\
     \norm{\bH} \asymp K^{-1}\norm{\btheta}^2 \; , \; & \norm{\bW} \leq C\sqrt{n\otheta\theta_{\text{max}}} \numberthis \label{basicresultseq7}
\end{align*}
\end{lem}
We also show bounds on various expressions involving $\mu_i, \theta_i$ and $\bxi_k$ for later use.
\\[1.5 pt]
\begin{lem}\label{basicsumresults}
Under event $\mathcal{A}_1$, the following bounds hold:
\begin{align}
|\trueregdegree^{-\frac{3}{2}}(a, a) \xi_k(a)| & \lesssim \frac{1}{n^{2}\otheta^2(\otheta \vee \theta_a)} \label{sum1} \\
\sum_{a \neq i} \mu_a \theta_a \trueregdegree^{-1}(a, a) & \lesssim \frac{\sqrt{\log n}}{\otheta} \label{sum2}\\
\sum_{a = 1}^n \xi_k(a) \trueregdegree^{-\frac{3}{2}}(a, a) & \lesssim \frac{1}{n\otheta^3} \label{sum3}\\
\sum_{a = 1}^n \sqrt{\frac{\mu_a \theta_a^2}{\minthetaa}} & \lesssim n\sqrt{\frac{\theta_a^3}{\otheta^2}} \label{sum4}\\
\sum_{a = 1}^n \theta_a\sqrt{\frac{\mu_a}{\otheta \vee \theta_a}} & \lesssim n\sqrt{\otheta} \label{sum5}
\end{align}
\end{lem}

\subsection{Proofs of Basic Lemmas}
\subsubsection{Proof of Lemma \ref{incoherence}}
    As shown in the proof of Lemma B.2 of \citet{ke2022optimal}, $\bB\bB^T = \bG^{-1}$. Since the entries of $\bG^{-1}$ are less than $c_1$, where $c_1$ is the constant in Assumption \ref{regconds1}(a), it follows that the entries of $\bB$ must be smaller than $\sqrt{c_1}$. By construction, 
    \begin{align*}
        |\xi_k(i)| & = |\bD_0^{-\frac{1}{2}}\theta_i \bpi_i^T \bB_{:, k}| \\
        & \leq |\bD_0^{-\frac{1}{2}}\theta_i|\max_{1 \leq l \leq k} |b_{l, k}| \\
        & \leq \sqrt{c_1}|\bD_0^{-\frac{1}{2}}||\theta_i| \leq \sqrt{C} \cdot \frac{\theta_i}{\sqrt{n\otheta(\maxthetai)}}
    \end{align*}
    where the last line is obtained from plugging in Equation \ref{basicresultseq2}, and $C'$ is a fixed constant. Setting $\mu_i = \frac{C\theta_i^2}{\otheta(\maxthetai)}$ suffices.

The lower bound of $\xi_1(i) \geq C' \cdot \frac{\theta_i}{\sqrt{n\otheta(\maxthetai)}}$ follows from Lemma B.2 of \citet{ke2022optimal}. As a result, $\xi_1(i) \asymp \sqrt{\frac{\mu_i}n}$.

\subsubsection{Proof of Lemma \ref{lem:order-eigenval}}
The first two inequalities in Lemma \ref{lem:order-eigenval} are shown in Lemma B.1 of \citet{ke2022optimal}. The third inequality follows by definition from Assumption \ref{regconds2}(a).

\subsubsection{Proof of Lemma \ref{lem:order}}
     As defined in Lemma \ref{incoherence}, $\mu_i \asymp \begin{cases}
\theta_i/\bar{\theta}, & i\in S_1,\\
\theta^2_i/\bar{\theta}^2, & i\in S_2,
\end{cases}$. As a result, Lemma \ref{lem:order} follows immediately from Lemma B.2 of \citet{ke2022optimal}.

\subsubsection{Proof of Lemma \ref{basicresults}} \label{basicresultspf} \renewcommand\labelitemi{-}
 \begin{itemize}
    \item \textbf{Proofs of \eqref{basicresultseq1}, \eqref{basicresultseq2}, \eqref{basicresultseq3},
    \eqref{basicresultseq8}:} 
   For any $i \in [n]$, the following bound always holds,
    \begin{align*}
        \sum_{j = 1}^n H_{ij} & = \theta_i \bpi_i \bP \bPi^T \bTheta \\
        & \geq \theta_i \sum_{k = 1}^K \bpi_i(k) P_{k, k} \sum_{j = 1}^n \theta_j\bpi_j(k) \\
        & \geq  \theta_i \sum_{k = 1}^K \bpi_i(k) \cdot 1 \cdot (c_1 \norm{\btheta}_1) \\
        & \geq c_1n\otheta\theta_i
    \end{align*}
    where the third line follows from Assumption \ref{regconds1}(a) and the fact that $\bP$ has unit diagonals. Furthermore,
    \begin{align*}
        \sum_{j = 1}^n H_{ij} & = \theta_i \sum_{j = 1}^n  \theta_j (\bpi_i \bP \bpi_j^T) \\
        & \leq \theta_i \sum_{j = 1}^n  \theta_j \norm{\bpi_i}_1\norm{\bP}_{\text{max}}\norm{\bpi_j}_1 \\ & \leq \theta_i \sum_{j = 1}^n  C'\theta_j = C'n\otheta\theta_i
    \end{align*}
    where the last line follows from Assumption \ref{regconds2}(e).

    In turn,
    \begin{align*}
        \trueregdegree(i, i) & = \sum_{j = 1}^n H_{ij} + \frac{\delta}{n} \sum_{i = 1}^n \sum_{j = 1}^n H_{ij} \\
        & \geq c_1n\otheta\theta_i + c_1n\otheta^2 \numberthis \label{basicresultspfeq1}
    \end{align*}

    Assume that equations \eqref{basicresultseq4} and \eqref{basicresultseq6} also hold (which will be shown to occur with probability $1 - o(n^{-11})$ later in this proof). By plugging in their bounds, it follows that
     \begin{align*}
        |\hat{D}(i, i) - \trueregdegree(i, i)| & = \sum_{j = 1}^n W_{ij} + \frac{\delta}{n} \sum_{i = 1}^n \sum_{j = 1}^n W_{ij} \\
        & \leq C\sqrt{n\otheta\theta_i \log n} + C'\delta\otheta\sqrt{ \log n} \leq C''\sqrt{n\otheta\theta_i \log n}
        \end{align*}
        where the last relation follows from Assumption \ref{regconds2}(d). Thus, \eqref{basicresultseq8} holds. To show \eqref{basicresultseq2}, 
    \begin{align*}
        \hat{D}(i, i) & = \trueregdegree(i, i) + (\hat{D}(i, i) - \trueregdegree(i, i))\\
        & \geq (c_1n\otheta\theta_i + c_1n\otheta^2) - C''\delta\otheta\sqrt{ \log n} \\
        & \geq C'''(\otheta\theta_i + n\otheta^2)
    \end{align*}
    where the last line follows from Assumption \ref{regconds2}(d).

On the other hand,
    \begin{align*}
         \norm{\bH_{i, \cdot}}^2 & = \sum_{j = 1}^n H_{ij}^2 \leq C\theta_i^2 \sum_{j = 1}^n \theta_j^2 \leq Cn\theta_{\text{max}}^2 \theta_i^2
     \end{align*}

    \item \textbf{Proofs of \eqref{basicresultseq4}, \eqref{basicresultseq6}, and \eqref{basicresultseq5}:} To show \eqref{basicresultseq4}, we use Bernstein's inequality. Since $\max_{i, j \leq n} |W_{ij}| \leq 1$, the inequality
     \begin{align*}
         \sum_{j = 1}^n W_{ij} & \leq \sqrt{\log n \sum_{j = 1}^n \mathbb{E}[W_{ij}^2]} + \log n \max_{i, j \leq n} |W_{ij}| \\ & \leq \sqrt{\log n \sum_{j = 1}^n H_{ij}} + \log n \\
         & \leq C\sqrt{n\otheta\theta_i \log n} + \log n \\
         & \lesssim \sqrt{n\otheta\theta_i \log n}(1 + \sqrt{\frac{\log n}{n\otheta\theta_i}}) \lesssim \sqrt{n\otheta\theta_i \log n}
     \end{align*}
     holds with probability $1 - o(n^{-12})$ for some constant $C$,
     where the third line follows from plugging in the bound in \eqref{basicresultseq1}, and the fourth line from the fact that $\otheta\theta_i \geq \frac{\log n}{n}$, due to Assumption \ref{regconds2}(d). By union bound, inequalities \eqref{basicresultseq6} holds simultaneously for all $j \in [n]$ with probability $1 - o(n^{-11})$.

     By another application of Bernstein's inequality, we can bound:
      \begin{align*}
          \sum_{i = 1}^n \sum_{j = 1}^n W_{ij} & \leq \sqrt{\log n \sum_{i = 1}^n \sum_{j = 1}^n \mathbb{E}[W_{ij}^2]} + \log n \\ & \leq \sqrt{\log n \sum_{i = 1}^n \sum_{j = 1}^n H_{ij}} + \log n \\
         & \leq C\sqrt{n^2\otheta^2 \log n} + \log n \\
         & \lesssim n\otheta\sqrt{\log n}
     \end{align*}
     with probability $1 - o(n^{-11})$.

     Likewise, since $|W_{ij}^2| \leq 1$ for any $j \leq n$, Bernstein's inequality implies the following with probability $1 - o(n^{-11})$,
 \begin{align*}
     \norm{\bW_{i, \cdot}}^2 & = \sum_{j = 1}^n W_{ij}^2 \lesssim \sqrt{\log n \sum_{j = 1}^n \mathbb{E}[W_{ij}^4]} + \log n \\ 
     & \leq \sqrt{\log n \sum_{j = 1}^n H_{ij}} + \log n \\
     & \lesssim \sqrt{n\otheta\theta_i \log n} \numberthis \label{lemmaAeq1}
 \end{align*}
 where the last line follows from plugging in the bound in \eqref{basicresultseq1}.

 \item \textbf{Proof of \eqref{basicresultseq7}:} The bound on $\norm{\bH}$ follows from Lemma C.2 of \citet{jin2022mixed}; the bound on $\norm{\bW}$ from equation (D.27) of \citet{jin2022mixed}. 

\end{itemize}

\subsection{Proof of Lemma \ref{basicsumresults}}
Since each expression in Lemma \ref{basicsumresults} comprises of terms that are bounded in Lemmas \ref{basicexpansiondegreeresults} and \ref{incoherence}, our strategy will be to simply plug in the bounds in said lemmas.

    \underline{Equation \eqref{sum1}:} Plugging in the bounds in Lemma \ref{incoherence} and \ref{basicexpansiondegreeresults}, we can deduce
    \begin{align*}
        |\trueregdegree^{-\frac{3}{2}}(a, a)||\xi_k(a)| & \lesssim \sqrt{\frac{\mu_a}{n^4\otheta^3(\otheta \vee \theta_a)^3}}\\
& \lesssim \frac{1}{n^2\otheta^2(\otheta \vee \theta_a)}
    \end{align*}

    \underline{Equation \eqref{sum2}:} Likewise, plugging in the bounds in Lemma \ref{incoherence} and \ref{basicexpansiondegreeresults} yields
    \begin{align*}
        \sum_{a \neq i} \mu_a \theta_a \trueregdegree^{-1}(a, a) & \lesssim \sum_{a \neq i} \frac{\sqrt{\theta_a^3(\otheta \wedge \theta_a) \log n}}{n\otheta^2(\theta_a \vee \otheta)} \\
        & \leq \sqrt{\frac{\log n}{n^2 \otheta^3}} \sum_{a \neq i} \sqrt{\theta_a} \\
        & \leq \sqrt{\frac{\log n}{n^2 \otheta^3}} \cdot n\sqrt{\otheta} \leq \frac{\sqrt{\log n}}{\otheta}
    \end{align*}
    where the bound on $\sum_{a \neq i} \sqrt{\theta_a}$ in the last line follows from applying the Cauchy-Schrawz inequality.

     \underline{Equations \eqref{sum3}, \eqref{sum4}, \eqref{sum5}:} Proceeding in a similar manner for the remaining sums,
     \begin{align*}
         \sum_{a = 1}^n \xi_k(a) \trueregdegree^{-\frac{3}{2}}(a, a) & \lesssim \sum_{a = 1}^n \frac{1}{n^2\otheta^3} \leq \frac{1}{n\otheta^3} \\
         \sum_{a = 1}^n \sqrt{\frac{\mu_a \theta_a^2}{\minthetaa}} & \lesssim \sum_{a = 1}^n \sqrt{\frac{\theta_a^3}{\otheta^2}} \\
         & \leq n\sqrt{\frac{\theta_a^3}{\otheta^2}} \\
         \sum_{a = 1}^n \theta_a\sqrt{\frac{\mu_a}{\otheta \vee \theta_a}} & \leq  \sum_{a = 1}^n \sqrt{\frac{\theta_a^3(\otheta \wedge \theta_a)}{\otheta^2(\theta_a \vee \otheta)}} \\
         & \leq \sum_{a, \theta_a \geq \otheta} \sqrt{\frac{\theta_a^2}{\otheta}} + \sum_{a, \theta_a \leq \otheta} \sqrt{\otheta} \leq 2n\sqrt{\otheta}
     \end{align*}

\section{Proofs of Main Theorems and Propositions}

\subsection{Proof of Proposition \ref{recoverP}}\label{recoverPproof}
The proof is similar to that of Lemma 2.2 in Jin et al. (2022) \citet{jin2022mixed}, so we simply provide an outline below. By imitating the same reasoning as in Equation (C.19) of \citet{jin2022mixed}, we have the equality 
$$
\bG(\bB \bLambda \bB^T) \bG = \bG\bP\bG
$$
where $B$ is the non-singular matrix defined by $\bXi = \boldsymbol{\trueregdegreehalf} \bTheta \bPi \bB$ (which is shown to exist, for instance, by the proof of Lemma 2.1 of \cite{ke2022optimal}). By definition, $\bB = \diag(\bb_1)^T\bQ^T$, where $\boldeta_1$ and $\bQ$ are defined in Section 3.1. Since $\bG$ and $\bP$ are invertible, $\bP = \bB \bLambda \bB^T = \diag(\bb_1)^T\bQ^T\bLambda \bQ \diag(\bb_1)$. On the other hand, $\bXi = \bD_0^{-\frac{1}{2}}\bTheta \bPi \bB$ by construction, which immediately implies the desired expression for $\theta_i$.

\subsection{Proof of Theorem \ref{estimationP1}}\label{estimationP1proof}
By triangle inequality,
\begin{align*}
    \norm{\hat{P}_{t(a)t(b)} - P_{ab}} \leq & \norm{\hat{b}_1(t(a))- b_1(a)}\norm{\hat{\bQ}_{t(a), :}\hat{\bLambda}\hat{\bQ}_{t(b), :}^T\hat{b}_1(t(b))} \\
    & + \norm{b_1(a)}\norm{\hat{\bQ}_{t(a), :}\hat{\bLambda}\hat{\bQ}_{t(b), :}^T - \bQ_{a, :} \bLambda \bQ_{b, :}^T}\norm{\hat{b}_1(t(b))} \\
    & + \norm{b_1(a)}\norm{\bQ_{a, :} \bLambda \bQ_{b, :}^T} \norm{\hat{b}_1(t(b)) - b_1(b)} \\
    & := \mathcal{B}_1 + \mathcal{B}_2 + \mathcal{B}_3
\end{align*}
Our approach will be to bound each of the terms on the RHS separately. Before doing so, we first need to characterize the deviation of our estimated vertices $\hat{\bQ}_{t(a), :}$ from the true vertices $\bQ_{a, :}$. Define
$$
 \Delta_{\br} = \min_{k\in [K]}\min_{i\in [n]\backslash \mathcal{C}_k}\left\|\br_i-\bv_k\right\|_2
$$
\begin{thm} \label{VHthm}
   Assume that 
   \begin{align*}
       \Delta_{\br} > 2\phi> (1 + C_{SP}) \cdot \epsilon_{\text{max}} \numberthis \label{vhasmp} 
   \end{align*} for some fixed constant $C_{SP}$, where  $\epsilon_{\text{max}} := C \sqrt{\frac{K\log n}{n\bar{\theta} (\bar{\theta}\wedge \min_{1 \leq i \leq n} \theta_i)\nu_n^2}}$ for some constant $C$ (which is defined properly later on, in Corollary \ref{lem:hatR}). With probability $1 - o(n^{-10})$, there exists a permutation $t$ of $[K]$, such that the outputs of Algorithm \ref{algVH} satisfy $\hat{\mathcal{C}}_{t(k)} = \mathcal{C}_{k}$ for all $k\in[K]$. As a result, for all $k \in [K]$,
   \begin{align*}
          \Vert \wh \bv_{t(k)} \bO_1 - \bv_k \Vert_{2} & \lesssim \sqrt{\frac{K\log n}{n\bar{\theta}(\otheta \wedge \min_{i \in \mathcal{C}_k} \theta_i)\nu_n^2}} := \epsilon_{0, k} \numberthis \label{VHthmeq1}
   \end{align*}
   where $\bT$ is the permutation matrix corresponding to $t$, i.e. $T_{k, :} = \mathbf{e}_{t(k)}$ for all $k \in [K]$, and $\bO_1$ is the rotation matrix in Corollary \ref{lem:hatR}.
\end{thm}
$\Delta_{\br}$ measures the minimum separation between any pure node and a mixed one. Our Assumption \eqref{vhasmp} ensures that the vertex hunting algorithm correctly identifies all pure nodes, and it is crucial to our analysis for obtaining optimal rates for $\bP$ and $\bTheta$. We leave the setting of a small separation between pure and mixed nodes for future research.  

Taking Theorem \ref{VHthm} for granted, we now turn to bounding the estimation rate of $\wh \bP$. In order to do so, we first recollect the following basic facts.
\begin{itemize}
    \item $b_1(k) \asymp 1, |\wh b_1(t(k)) - b_1(k)| \lesssim \epsilon_{0, k}\sqrt{K} + \sqrt{\frac{K \log n}{n\otheta^2\lambda_K^2}}$, where $\epsilon_{0, k}$ is defined in Equation \eqref{VHthmeq1}.
    \begin{itemize}
        \item Proof: $b_1(k) \asymp 1$ follows from the proof of Lemma B.2 of \citet{ke2022optimal}. Likewise, a similar argument to the proof of equation (D.5) of \citet{ke2022optimal} shows
        \begin{align*}
            |\wh b_1(t(k)) - b_1(k)| & \lesssim K^{\frac 12} |\lambda_2| \Vert \wh \bv_{t(k)} \bO_1 - \bv_k \Vert_{2} + \big|\bv_k^T (\bO_1^T\hat{\bLambda}_{-1} \bO - \bLambda_1)\bv_k\big| \\
            & \lesssim  K^{\frac 12} |\lambda_2|\epsilon_{0, k} + \big|\bv_k^T (\bO_1^T\hat{\bLambda}_{-1} \bO - \bLambda_1)\bv_k\big| \\
            & \lesssim K^{\frac 12} |\lambda_2|\epsilon_{0, k} + \sqrt{\frac{K \log n}{n \otheta^2 \lambda_K^2}} \numberthis \label{estimationPeq2}
        \end{align*}
        where $\bT$ is the permutation matrix defined in Theorem \ref{VHthm}, and the last line follows from the proof of equation (D.5) in \cite{ke2022optimal} again.
    \end{itemize}
    \item For all rows $k$, $\norm{\bQ_{k, :}}_2 \leq C\sqrt{K}$. Furthermore, $\hat{\bQ}_{t(k), :} \lesssim \norm{\bQ_{k, :}}$, since
    \begin{align*}
    \norm{\hat{\bQ}_{t(k), :}} & = \norm{\hat{\bQ}_{t(k), :}\bO_1} \\
    & \leq \norm{\hat{\bQ}_{t(k), :}\bO_1 - \bQ_{k, :}} + \norm{\bQ_{k, :}} \lesssim \norm{\bQ_{k, :}}
\end{align*}
\end{itemize}
By the submultiplicity of the spectral norm, $\mathcal{B}_1$ can be bounded by 
\begin{align*}
    \mathcal{B}_1 & \leq \norm{\hat{b}_1(t(a)) - b_1(a)}\norm{\hat{\bQ}_{t(a), :}^T}\norm{\hat{\bLambda}}\norm{\hat{\bQ}_{t(b), :}}\norm{\hat{b}_1(t(b))} \\
    & \lesssim C K^{\frac 32} \lambda_1|\lambda_2|\epsilon_{0, a} + \sqrt{\frac{K^3 \lambda^2_1 \log n}{n \otheta^2 \lambda_K^2}}
\end{align*} In the same manner, the third RHS term can also be bounded by the same quantity. Regarding the second term, the following bound holds by a similar argument as to the proof of Equation (D.5) in \citet{ke2022optimal}:
\begin{align*}
   \norm{\hat{\bQ}_{t(a), :}\hat{\bLambda}\hat{\bQ}_{t(b), :}^T - \bQ_{a, :} \bLambda \bQ_{b, :}^T} \leq  K^{\frac 12} |\lambda_2|\Vert \bT\widehat{\bV}\bO - \bV\Vert_{2\to \infty} + \big|\bv_a^T (\bO_1^T\hat{\bLambda}_{-1} \bO - \bLambda_1)\bv_b\big|
\end{align*} The quantity on the RHS has already been bounded in \eqref{estimationPeq2}. Pulling our bounds for the $\mathcal{B}$ terms all together, we thus ultimately obtain a rate of
\begin{align*}
    \norm{\hat{P}_{t(a)t(b)} - P_{ab}} & \leq C K^{\frac 32} \lambda_1|\lambda_2|\max\{\epsilon_{0, a}, \epsilon_{0, b}\}  + \sqrt{\frac{K^3 \lambda^2_1 \log n}{n \otheta^2 \lambda_K^2}} \\
    & \lesssim \sqrt{\frac{K^4\lambda_1^2\lambda_2^2}{n\otheta(\otheta \wedge (\minthetaab))\nu_n^2}} + \sqrt{\frac{K^3 \lambda^2_1 \log n}{n \otheta^2 \lambda_K^2}}
\end{align*}

\subsection{Proof of Theorem \ref{VHthm}} \label{VHthmpf}
Proof: Our strategy will be to first bound the deviation $|\wh \br_i - \br_i|$ of each node, then to bound the deviation of our vertices through analyzing Algorithm \ref{algVH}. By Corollary 3.1 of \citet{ke2022optimal}, we know the following.
\begin{cor} (Modified Corollary 3.1 of \cite{ke2022optimal}) \label{lem:hatR} 
With probability $1-o(n^{-10})$, there exists an orthogonal matrix $\bO_1\in \mathbb{R}^{K-1, K-1}$ such that, simultaneously for $i \in [n]$:
\begin{equation} \label{bound-hatR}
\Vert \bO_1^T \hat{\br}_i - \br_i \Vert \leq C \sqrt{\frac{K\log n}{n\bar{\theta} (\bar{\theta}\wedge \theta_i)\nu_n^2}}
\end{equation}
Define $\epsilon_{\text{max}} = C \sqrt{\frac{K\log n}{n\bar{\theta} (\bar{\theta}\wedge \min_{1 \leq i \leq n} \theta_i)\nu_n^2}}$. Then $\Vert \bO_1^T \hat{\br}_i - \br_i \Vert \leq C\epsilon_{\text{max}}$ for all $i \in [n]$.
\end{cor}
The only difference between our version and that presented in \citet{ke2022optimal} is the latter requires the nodes to have degree exceeding a certain threshold. However, that assumption is not required in their proof, so we can freely drop it.

Now, we turn to analyzing the vertex hunting step. Due to our choice of SP as our vertex hunting algorithm (\cite{jin2022mixed}; Lemma 3.1), we know \begin{align}\label{est:vertex}
 \norm{\bT\widehat{\bV} \bO_1- \bV}_{2\to \infty}  \leq C_{SP} \max_{i\in [n]} \Vert \hat{\br}_i \bO_1 - \br_i \Vert.
\end{align}
for some $K\times K$ permutation matrix $\bT$ and fixed constant $C_{SP}$. Let $t$ be the permutation corresponding to $\bT$, i.e. for all $1 \leq k \leq K$, $t(k)$ is equal to the index of the nonempty column in row $k$. Recall that we defined $i_k$ to be the index of the node selected in the $k$th round step of the successive projection algorithm (Steps 3-5 of Algorithm \ref{algVH}). We first show that node $i_{t(k)}$ must belong to $\mathcal{C}_k$. Indeed, by triangle inequality,
\begin{align*}
    \norm{\br_{i_{t(k)}} - \bv_k}_2 & \leq |\br_{i_{t(k)}} - \hat{\br}_{i_{t(k)}}\bO_1| + |\hat{\br}_{i_{t(k)}}\bO_1 - \bv_k| \\
    & \leq  \max_{i\in [n]} \Vert \hat{\br}_i \bO_1 - \br_i \Vert + |\hat{\br}_{i_{t(k)}}\bO_1 - v_k| \\
    & \leq (1 + C_{SP})\sqrt{K}\epsilon_{\text{max}}
\end{align*}
    with probability $1 - o(n^{-10})$, where the bounds in the last inequality come from Corollary \ref{lem:hatR} and equation \eqref{est:vertex}. If $i_{t(k)}$ does not belong to $\mathcal{C}_k$, the LHS is at least $\Delta_r$, a contradiction. Thus, $i_{t(k)}$ must belong to $\mathcal{C}_k$. As this holds for all $1 \leq k \leq K$, it follows that the estimated vertices consist of $K$ pure nodes, one from each community.

    Next, we show that the estimated vertex set $\hat{\mathcal{C}}_{t(k)}$ contains all pure nodes $\mathcal{C}_k$ WHP. For any arbitrary node $i \in [n]$, triangle inequality implies
    \begin{align*}
        \norm{\hat{\br}_i - \hat{\br}_{i_{t(k)}}} & \leq  \norm{\hat{\br}_i - \br_i\bO_1^T} + \norm{\br_i\bO_1^T - \br_{i_{t(k)}}\bO_1^T} + \norm{\br_{i_{t(k)}}\bO_1^T - \hat{\br}_{i_{t(k)}}} \\
        & = \norm{\hat{\br}_i \bO_1 - \br_i} + \norm{\br_i - \br_{i_{t(k)}}} + \norm{\br_{i_{t(k)}}- \hat{\br}_{i_{t(k)}}\bO_1} \\
        & \lesssim \norm{\br_i - \br_{i_{t(k)}}} + (1 + C_{SP})\sqrt{K}\epsilon_{\text{max}}
    \end{align*}
    with probability $1 - o(n^{-10})$. If $i \in \mathcal{C}_k$, the RHS is therefore bounded by $(1 + C_{SP})\sqrt{K}\epsilon_{\text{max}} < \phi$, implying that $i \in \hat{\mathcal{C}}_{t(k)}$. In contrast, if $i \notin \mathcal{C}_k$, a similar triangle inequality argument implies
    \begin{align*}
         \norm{\hat{\br}_i - \hat{\br}_{i_{t(k)}}} & \geq \norm{\br_i - \br_{i_{t(k)}}} - \norm{\br_i - \br_{i_{t(k)}}} - \norm{\br_{i_{t(k)}}- \hat{\br}_{i_{t(k)}}\bO_1} \\
         & \geq \norm{\br_i - \br_{i_{t(k)}}} - (1 + C_{SP})\sqrt{K}\epsilon_{\text{max}}
    \end{align*}
    By the condition $\Delta_r > 2\phi$, however, the RHS is bounded below by $\phi$, meaning $i$ cannot lie in $\hat{\mathcal{C}}_{t(k)}$. In summary, $\hat{\mathcal{C}}_{t(k)} = \mathcal{C}_k$ for all communities, as desired.

\subsection{Proof of Theorem \ref{estimationP2}}\label{estimationP2proof}
\begin{proof}
    The key result needed for our proof are the following improved $l_{\infty}$ eigenvector bounds, which we will obtain using eigenvector expansions.
    \begin{thm} \label{improveddevbounds}
        Under Assumptions \ref{regconds1}, \ref{regconds2} and Assumption \ref{regconds3}, the following improved eigenvector bounds hold with probability $1 - o(n^{-10})$ for all $i \in [n], k \in [K]$:
        \begin{align*}
            |\hat{\xi}_k(i) - \xi_k(i)| & \lesssim \frac{\sqrt{\log n}}{n\otheta\lambda_k} \\
            |\hat{r}_i(k - 1) - \br_i(k - 1)| & \lesssim \sqrt{\frac{\log n}{n\theta_i(\minthetai)\lambda_K}} \numberthis \label{improvedeq1}
        \end{align*}
    \end{thm}
    Theorem \ref{improveddevbounds} improves over the eigenvector bounds in \citet{ke2022optimal} by i) providing deviation bounds for individual entries $|\wh \br_i(k - 1) - \br_i(k - 1)|$ of the rows, as opposed to a bound on the entire row and ii) depending more optimally on $\theta_i$, since the error on the RHS of \eqref{improvedeq1} continues to decrease even when $\theta_i > \otheta$. Equipped with these improved bounds, we can subsequently refine our bounds on the estimated vertices.
    \begin{thm} \label{VHthm2}
   Assume that the conditions of Theorem \ref{VHthm} are satisfied, and let $T$ be the permutation of vertices in Theorem \ref{VHthm}. If we additionally assume Assumption  \ref{regconds3}, the estimated vertices output by Algorithm \ref{algVH} satisfy
   \begin{align*}
    \displaystyle |\hat{v}_{t(a)}(k - 1) - v_a(k - 1)| \lesssim \sqrt{\frac{\log n}{n \lambda_K \min_{i \in \mathcal{C}_a}(\theta_i(\minthetai))}} \numberthis \label{VHthm2eq1}
   \end{align*}
   with probability $1 - o(n^{-10})$ for all $a \in [K]$.
\end{thm}

Taking these theorems for granted, we now establish an improved estimation rate for $\wh \bP$. First, we show the following basic facts, which are improvements over those established in Theorem \ref{estimationP1}. 
\begin{itemize}
    \item $|\wh b_1(t(a)) - b_1(a)| \lesssim \sqrt{\frac{K^2\log n}{n \lambda_K \min_{i \in \mathcal{C}_a}(\theta_i(\minthetai))}} + \sqrt{\frac{K^2\log n}{n \otheta^2}}$
    \begin{itemize}
        \item \begin{proof} Since $b_1(k) \asymp 1$, we first note that $|\wh b_1(t(a)) -b_1(a)| \asymp |\frac{1}{\wh b_1(t(a))^2} - \frac{1}{b_1(a)^2}|$. By triangle inequality,
        \begin{align*}
            |\frac{1}{\wh b_1(t(a))^2} - \frac{1}{b_1(a)^2}| & \leq |\wh \lambda_1 - \lambda_1| +  \big|\wh \bv_{t(a)}^T \hat{\bLambda}_{-1} \wh \bv_{t(a)} -   \bv_a^T {\bLambda}_{-1} \bv_a\big| \numberthis \label{estimationP2eq3}
        \end{align*}
        The first RHS term can be bounded via Weyl's inequality. Regarding the second term, further decomposition shows
        \begin{align*}
            \big|\wh \bv_{t(a)}^T \hat{\bLambda}_{-1} \wh \bv_{t(a)} -   \bv_a^T {\bLambda}_{-1} \bv_a\big| & = \sum_{k = 2}^{K} ( \wh v_{t(a)}(k - 1))( \wh v_{t(a)}(k - 1))\wh \lambda_{k} - v_a(k - 1)v_b(k - 1) \lambda_{k} \\
            & \leq \sum_{k = 2}^K | \wh v_{t(a)}(k - 1) - v_a(k - 1)||\wh v_{t(a)}(k - 1)||\wh \lambda_{k}| \\
            & + |v_a(k - 1)||\wh v_{t(b)}(k - 1) - v_b(k - 1)||\wh \lambda_{k}| \\
            & + |v_a(k - 1)||v_b(k - 1)||\wh \lambda_{k} - \lambda_k| \numberthis \label{estimationP2eq4}
        \end{align*}
        Plugging in our bounds from Theorem \ref{VHthm2}, we can bound each of the terms in the summand
        \begin{align*}
            | \wh v_{t(a)}(k - 1) - v_a(k - 1)||\wh v_{t(a)}(k - 1)||\wh \lambda_{k}| & \lesssim \sqrt{\frac{\log n}{n \lambda_K \min_{i \in \mathcal{C}_a}(\theta_i(\minthetai))}} \cdot 1 \cdot \lambda_k \\
            & \leq \sqrt{\frac{\log n}{n\lambda_K \min_{i \in \mathcal{C}_a}(\theta_i(\minthetai))}}\\
            |v_a(k - 1)||\wh v_{t(b)}(k - 1) - v_b(k - 1)||\wh \lambda_{k}| & \lesssim \sqrt{\frac{\log n}{n \lambda_K \min_{i \in \mathcal{C}_b}(\theta_i(\minthetai))}} \\
            |v_a(k - 1)||v_b(k - 1)||\wh \lambda_{k} - \lambda_k| & \lesssim 1 \cdot 1 \cdot \sqrt{\frac{\log n}{n \otheta^2}} \numberthis \label{estimationP2eq5}
        \end{align*}
        where $|\hat{\lambda}_k - \lambda_k| \leq \norm{\hat \bL_0 - \bL_0} \lesssim \sqrt{\frac{\log n}{n \otheta^2}}$ by Weyl's inequality. Note that the bounds in \eqref{estimationP2eq5} are independent of $k$. Combining our bounds and summing over all $k \in [K]$ in \eqref{estimationP2eq4}, we thus see
        \begin{align*}
            \displaystyle \big|\wh \bv_{t(a)}^T \hat{\bLambda}_{-1} \wh \bv_{t(a)} - \bv_a^T {\bLambda}_{-1} \bv_a\big| \lesssim \sqrt{\frac{K^2\log n}{n \lambda_K \min_{i \in \mathcal{C}_a}(\theta_i(\minthetai))}} + \sqrt{\frac{K^2\log n}{n \otheta^2}}
        \end{align*}
        which, in conjunction with \eqref{estimationP2eq3}, yields
        \begin{align*}
            |\wh b_1(t(a)) - b_1(a)| \lesssim \sqrt{\frac{K^2\log n}{n \lambda_K \min_{i \in \mathcal{C}_a}(\theta_i(\minthetai))}} + \sqrt{\frac{K^2\log n}{n \otheta^2}} \numberthis \label{estimationPeq6}
        \end{align*}
        \end{proof}
            \end{itemize}
\end{itemize}
Define the vectors $q_{a} = \bQ_{a, :}, \hat{q}_a = \hat{\bQ}_{a, :} \in \mathbb{R}^{1 \times K}$ for all $a \in [K]$. By triangle inequality,
\begin{align*}
    |\hat{P}_{t(a)t(b)} - P_{ab}| & = |\hat{b}_1(t(a)) \hat{q}_{t(a)} \hat{\bLambda} \hat{q}_{t(b)}^T \hat{b}_1(t(b)) - b_1(a)q_{a} \bLambda q_{b} b_1(b)| \\
    & = |\sum_{k = 1}^K \hat{b}_1(t(a))   \hat{b}_1(t(b))  \hat{q}_{t(a)}(k)  \hat{q}_{t(b)}(k)  \hat{\Lambda}_{k, k} - b_1(a) b_1(b) q_{a}(k) q_{b}(k) \Lambda_{k, k})| \\
    & \leq \sum_{k = 1}^K \big(| \hat{b}_1(t(a)) - b_1(a)||  \hat{b}_1(t(b))| |\hat{q}_{t(a)}(k)|| \hat{q}_{t(b)}(k)| \lambda_k \\
    & + |b_1(a)||\hat{b}_1(t(b)) - b_1(b)|| \hat{q}_{t(a)}(k)|| \hat{q}_{t(b)}(k)| \lambda_k \\
    & + |b_1(a)||b_1(b)|| \hat{q}_{t(a)}(k) - q_{a}(k)|| \hat{q}_{t(b)}(k)| \lambda_k \\
    & + |b_1(a)||b_1(b)||q_{a}(k)|| \hat{q}_{t(b)}(k) - q_{b}(k)| \lambda_k \\
    & + |b_1(a)||b_1(b)||q_{a}(k)||q_{b}(k)|| \hat{\lambda}_k - \lambda_k| \big) \numberthis \label{estimationPproofeq1}
\end{align*}
Plugging in our bounds from \eqref{estimationPeq6} and Theorem \ref{VHthm2} into the RHS of \eqref{estimationPproofeq1} yields a rate of 
$$
\displaystyle |\hat{P}_{t(a)t(b)} - P_{ab}| \lesssim\sqrt{\frac{K^4 \log n}{n \lambda_K \min_{i \in \mathcal{C}_a \cup \mathcal{C}_b}(\theta_i(\minthetai))}} + \sqrt{\frac{K^2\log n}{n \otheta^2}} + \sqrt{\frac{K^4 \lambda_1^2 \log n}{n\otheta^2}}
$$
\end{proof}

\subsection{Proof of Theorem \ref{estimationTheta}}
\begin{proof}
By triangle inequality,
    \begin{align*}
    |\hat{\Theta}(i, i) - \Theta(i, i)| & \leq |\hat{\xi}_1(i) - \xi_{1}(i)||\trueregdegreehalf(i, i)(\bpi_i^T\bb_1)^{-1}| \\
    & + |\hat{\xi}_{1}(i)||\obsregdegreehalf(i, i) - \trueregdegreehalf(i, i)||(\bpi_i^T\bb_1)^{-1}| \\
        & + |\hat{\xi}_1(i)||\obsregdegreehalf(i, i)||(\hat{\bpi}_i^T\hat{\bb}_1)^{-1} - (\bpi_i^T\bb_1)^{-1}| \\
    & = [1] + [2] + [3] \numberthis \label{estimationThetaeq0}
\end{align*}
To bound terms [1] and [3] in the above equation, we borrow the following bounds from \citet{ke2022optimal}.
\begin{thm}\label{crudexi1bound} (Theorem 3.1, Ke and Wang)
Under Assumptions \ref{regconds1} and \ref{regconds2}, there exists $\omega\in \{1, -1\}$ such that 
\begin{align}
|\omega\hat{\xi}_1(i)-\xi_1(i)| & \leq C \sqrt{\frac{K\theta_i\log(n)}{n^2\bar{\theta}^3\lambda_1^2 }}\Biggl(1+ \sqrt{\frac{\log(n)}{n\bar{\theta}\theta_i}}\, \Biggr),\label{result-1}
\end{align}
simultaneously for all $i \in [n]$, where $C$ is a fixed constant.
\end{thm}
\begin{thm} (Theorem 3.2, Ke and Wang) \label{crudepibound}
Let $\hat{\bPi}$ be the estimator output by the Mixed-SCORE-Laplacian (Ke and Wang, 2022), where the tuning parameters are such that $c>0$ and $0<\gamma<c_0$ (here $c_0$ is the same as in Assumption \eqref{regconds4}). Then, with probability $1-o(n^{-10})$, there exists a permutation $\bT$ on $\{1,2,\ldots,K\}$, such that 
$$
\norm{\bT\hat{\bpi}_i- \bpi_i}_1 \leq  C \min\left\{ \sqrt{\frac{K^3\log n}{n\bar{\theta} (\bar{\theta}\wedge \theta_i)\nu_n^2}},\;\; 1\right\}
$$
simultaneously for all $1\leq i\leq n$, where $C$ is a fixed constant.
\end{thm}

Both [1] and [2] in \eqref{estimationThetaeq0} consist of terms that are bounded in Theorems \ref{crudexi1bound},  \ref{crudepibound} and Lemmas \ref{lem:order}, \ref{basicexpansiondegreeresults}. Plugging in the corresponding bounds and using the fact that $\frac{\log n}{n\otheta \theta_i} \ll 1$ for all $i \in [n]$ (by Assumption \ref{regconds2}(d)) yields
\begin{align*}
    [1] & \leq \sqrt{\frac{K \theta_i^2 \log n}{n\otheta(\minthetai)}} \\
    [2] & \leq \sqrt{\frac{\theta_i^3\log n}{n\otheta(\maxthetai)^2}} \numberthis \label{ubthetaeq1}
\end{align*}
To bound (3), we apply triangle inequality again
\begin{align*}
    [3] & \leq |\hat{\xi}_1(i)||\hat{\trueregdegree}(i, i)||\bpi_i^T\bb_1|^{-2}|(\norm{\hat{\bpi}_i - \bpi_i}_{1}\norm{\bb_1}_{\infty} + |\hat{\bpi}_i|_{1}\norm{\hat{\bb}_1 - \bb_1}_{\infty})
\end{align*}
By equation (D.5) of \citet{ke2022optimal}, we can bound the $\norm{\hat{\bb}_1 - \bb_1}_{\infty}$ term in the above expression by $\sqrt{\frac{K\log n}{n\otheta^2 \lambda_k^2}}$. The rest of the terms can be bounded by again using Theorems \ref{crudexi1bound},  \ref{crudepibound} and Lemmas \ref{basicexpansiondegreeresults}. Plugging in the bounds yields $[3] \leq \sqrt{\frac{K\theta_i^2\log n}{n\bar{\theta} (\bar{\theta}\wedge \theta_i)\nu_n^2}} + \sqrt{\frac{K\theta_i^2\log n}{n\otheta^2 \lambda_k^2}}$. In conjunction with \eqref{ubthetaeq1}, we see
$$
|\hat{\Theta}(i, i) - \Theta(i, i)| \leq \sqrt{\frac{K\theta_i^2\log n}{n\bar{\theta} (\bar{\theta}\wedge \theta_i)\nu_n^2}} + \sqrt{\frac{K\theta_i^2\log n}{n\otheta^2 \lambda_k^2}}
$$
\end{proof}

\subsection{Proof of Theorem \ref{lowerboundsP}}\label{lowerboundsPproof}
\begin{proof}
To prove the desired lower bound, we restrict to a much smaller parameter space. Specifically, let $\mathcal{Q} = (\bTheta, \bPi, \bP)$ denote the parameters of a DCMM model. Define $\tilde{\mathcal{S}}_n$ to be a parameter space consisting of solely two parameter sets: $\mathcal{Q}^{\mu} = (\bTheta^{\mu}, \bPi^{\mu}, \bP^{\mu})$ and $\mathcal{Q}^{\alpha} = (\bTheta^{\alpha}, \bPi^\alpha, \bP^\alpha)$. If we can construct $\mathcal{Q}^{\mu}, \mathcal{Q}^{\alpha}$ such that $|P_{ab}^{\alpha} - P_{ab}^{\mu}| \geq C(n\otheta\tilde \theta)^{-\frac{1}{2}}$ and $KL(\mathcal{Q}^{\alpha}, \mathcal{Q}^{\mu}) \leq C'$ for constants $C, C'$, standard lower bound techniques [e.g. (Lemma 2.9 of \cite{Tsybakov_2009})] will imply a lower bound of $n^{-\frac{1}{2}}$ for the smaller space $\tilde{\mathcal{S}}_n$ and thus also the larger space $\tilde {\cal Q}_{n, ab, c_{10}}(\sigma)$. 

To this end, the main idea for the construction is to not work with the original parameters, but to instead reparameterize the DCMM in terms of the vectors $\bZ = \bTheta \bPi \bB \bLambda^{\frac{1}{2}}$ and $\bY = \bB\bLambda^{\frac{1}{2}}$, such that $\bZ\bZ^T = \bH$ and $\bY\bY^T = \bP$. In the following lemma, we prove that such a parameterization is indeed valid.
\begin{lem}\label{reparamlem}
    Let $\mathcal{Y}$ be the set of all nonnegative, non-singular matrices $\bY \in \mathbb{R}^{K \times K}$ such that $\bY\bY^T$ has unit diagonals. Furthermore, given a fixed $\bY \in \mathcal{Y}$, let $\mathcal{X}_Y$ be the set of all nonnegative matrices $\bZ \in \mathbb{R}^{n \times K}$ such that $\bZ = \bTheta' \bPi' \bY$, where
    \begin{itemize}
        \item $\bTheta' \in \mathbb{R}^{n \times n}$ is a nonnegative diagonal matrix.
        \item $\bPi' \in \mathbb{R}^{n \times K}$ is a nonnegative matrix that a) contains at least one row of the form $e_k = (0, 0, \cdots 1, \cdots, 0)$ for all $k \in [K]$ and b) has rows summing to one. In other words, $\bPi'$ is a valid membership matrix.
    \end{itemize}
    Then for any pair $(\bZ, \bY)$ in $\mathcal{X}_Y$ and $\mathcal{Y}$ respectively, there exists a unique DCMM model $(\bTheta, \bPi, \bP)$ satisfying Assumption \ref{regconds1}(d) and the relation
    \begin{align*}
        \bZ\bZ^T = \bTheta \bPi \bP \bPi^T \bTheta \numberthis \label{lowerboundsPproofeq1}
    \end{align*}
    In particular,
    \begin{align*}
        \bTheta = \bTheta', \bPi = \bPi', \bP = \bY\bY^T
    \end{align*}
\end{lem}
The proof of this Lemma is presented after the main proof, in Section \ref{reparamlempf}. Taking it to be true for the moment, we show the lower bound. For the sake of convenience, we 
\begin{itemize}
    \item Assume WLOG that $\theta_{a, \text{min}} \leq \theta_{b, \text{min}}$.
    \item Subsequently, we reorder the communities such that $\mathcal{C}_a$ is the first community.
    \item Assume WLOG that the pure node with the smallest degree in $\mathcal{C}_a$ appears first in the degree matrix $\bTheta$. 
\end{itemize} Consider the null hypothesis $\mathcal{Q}^{\mu} = (\bZ^{\mu}, \bY^{\mu})$, where \small
\begin{equation}
   \bY^{\mu} = \begin{bmatrix}
    \sqrt{1 - c_{12} (K - 1)\lambda_K} & \sqrt{\frac{c_{12} K\lambda_K}{2}} & \sqrt{\frac{c_{12} K\lambda_K}{6}} & \cdots &  \sqrt{\frac{c_{12} K \lambda_K}{(K - 1)^2 + (K - 1)}} \\
     \sqrt{1 - c_{12} (K - 1)\lambda_K} & -\sqrt{\frac{c_{12} K\lambda_K}{2}} & \sqrt{\frac{c_{12} K\lambda_K}{6}} & \cdots &  \sqrt{\frac{c_{12} K \lambda_K}{(K - 1)^2 + (K - 1)}} \\
      \sqrt{1 - c_{12} (K - 1)\lambda_K} &  0 & -2\sqrt{\frac{c_{12} K\lambda_K}{6}} & \cdots &  \sqrt{\frac{c_{12} K \lambda_K}{(K - 1)^2 + (K - 1)}} \\
      \vdots & \vdots & \vdots & \vdots & \vdots \\
       \sqrt{1 - c_{12} (K - 1)\lambda_K} & 0 & 0 & 0 & -(K - 1)\sqrt{\frac{c_{12} K \lambda_K}{(K - 1)^2 + (K - 1)}}
\end{bmatrix} \label{defYmu}
\end{equation}
\normalsize
and $c_{12}$ is a sufficiently small constant to be determined later. Also, define
$$
\bZ^{\mu} = \bTheta^{\mu} \begin{bmatrix}
\mathbf{e}_1 \\
\mathbf{1}_{\frac{2n}{K(K + 1)} - 1}^T \begin{bmatrix}
    \frac{3}{4} & \frac{1}{4} & 0 & \cdots & 0
\end{bmatrix} \\
    \mathbf{1}_{\frac{4n}{K(K + 1)}} \mathbf{e}_2^T \\
    \vdots \\
    \mathbf{1}_{\frac{2Kn}{K(K + 1)}} \mathbf{e}_K^T
\end{bmatrix} \bY^{\mu}
$$
where $\bTheta^{\mu}$ is defined such that
$$
\Theta^{\mu}(i, i) = \begin{cases}
    \tilde{\theta} & \hbox{ if $i = 1$}\\
    \check{\theta} & \hbox{ otherwise}
\end{cases}
$$
and $\sum_{i = 1}^n \Theta^{\mu}(i, i) = n\otheta$. In other words, every node has the same degree, except for the first node in community $a$ (which has been relabled to community one WLOG). We note that, by construction, the $k$th community in $\bPi^{\mu}$ for $k \geq 2$ has $n \cdot \frac{2k}{K(K + 1)}$ pure nodes in each community. 

In contrast, we define the alternative hypothesis $\mathcal{Q}^{\alpha} = (\bZ^{\alpha}, \bY^{\alpha})$ by \small
$$
\bY^{\alpha} = \begin{bmatrix}
    \; \sqrt{1 - c_{12} (K - 1)\lambda_K} - \Delta_Y & \sqrt{\frac{\lambda_k}{2} + 2 \Delta_Y \sqrt{1 - c_{12} (K - 1)\lambda_K} - \Delta_Y^2} & \sqrt{\frac{c_{12} K\lambda_K}{6}} & \cdots &  \sqrt{\frac{c_{12} K \lambda_K}{(K - 1)^2 + (K - 1)}} \; \\
     & \bY^{\mu}_{-1, :} & & 
\end{bmatrix}
$$ \normalsize
where
$$
\Delta_Y := c_0(n\otheta\tilde \theta)^{-\frac{1}{2}}
$$
for a properly small $c_0$ to be chosen later, and
$$
\bZ^{\alpha} = \begin{bmatrix}
     \theta_1 * \bY^{\alpha}_{1, :} \\
     \bZ^{\mu}_{-1, :}
\end{bmatrix}
$$
i.e. only the first row of $\bY^{\mu}$ and $\bZ^{\mu}$ are modified. To verify that $(\bZ^{\alpha}, \bY^{\alpha})$ comprise a valid pair, we note that for any row $\bZ^{\alpha}_{i, :}$ of $\bZ^{\alpha}$, the line between the origin and the point $\bZ^{\alpha}_{i, :}$ intersects the K-simplex formed by the vertices of $\bY^{\alpha}$, i.e. there exists a scalar multiple $r\bZ^{\alpha}_{i, :}$ of each row $\bZ^{\alpha}_{i, :}$ such that $r\bZ^{\alpha}_{i, :}$ lies within the K-simplex formed by $\bY^{\alpha}$. Thus, there exist matrices $\bTheta', \bPi'$ such that $\bZ^{\alpha} = \bTheta' \bPi' \bY^{\alpha}$. Moreover, $\bZ^{\alpha}$ contains a scalar multiple of each of the rows of $\bY^{\alpha}$, so $\bPi'$ contains an identity submatrix.

By Lemma \ref{reparamlem}, the null and alternative hypotheses correspond to DCMM models. Denote their corresponding models by $(\bTheta^{\mu}, \bPi^{\mu}, \bP^{\mu})$ and $(\bTheta^{\alpha}, \bPi^{\alpha}, \bP^{\alpha})$ respectively. First, we verify $\mathcal{H}^{\mu}, \mathcal{H}^{\alpha}$ satisfy Assumption \ref{regconds1}. Routine computation shows
\begin{align*}
    \bG^{\mu} & = K\frac{(\theta_1^{\mu})^2}{D_0^{\mu}(1, 1)} \mathbf{e_1}^T\mathbf{e_1} + K\Big(\sum_{i=2}^{\frac{2n}{K(K + 1)}} \frac{(\theta_i^{\mu})^2}{D_0^{\mu}(i,i)}\Big) (\begin{bmatrix}
        \frac{3}4 & \frac{1}{4} & 0 & \cdots & 0
    \end{bmatrix})^T\begin{bmatrix}
        \frac{3}4 & \frac{1}{4} & 0 & \cdots & 0
    \end{bmatrix} \\
    & + K \textbf{diag} \Big(0, \sum_{i\in \mathcal{C}_{2}}\frac{(\theta_i^{\mu})^2}{D_0^{\mu}(i,i)} \, , \cdots, \sum_{i\in \mathcal{C}_{K}}\frac{(\theta_i^{\mu})^2}{D_0^{\mu}(i,i)} \Big) 
\end{align*}
Since $\theta_j$ is the same for $j \geq 2$, it follows that for pure communities $k \geq 2$,
\begin{align*}
   \sum_{i\in \mathcal{C}_{k}}\frac{(\theta_i^{\mu})^2}{D_0^{\mu}(i,i)} \asymp \frac{|\mathcal{C}_k|}{n} = \frac{2k}{K(K + 1)} \numberthis \label{lowerboundthetaeq1}
\end{align*}
and
\begin{align*}
   \sum_{i=2}^{\frac{2n}{K(K + 1)}} \frac{\theta_i^2}{D_0^{\mu}(i,i)} \asymp \frac{2}{K(K + 1)} \numberthis \label{lowerboundthetaeq15}
\end{align*}
It is straightforward to deduce $\norm{\bG^{\mu}} \leq C$ for some fixed constant $C$, up to factors of $K$. To show $\norm{(\bG^{\mu})^{-1}}$ is also bounded above, we use the Gergoshonin circle theorem (see i.e. \cite{horn2012matrix}) to lower bound $\lambda_K(\bG^{\mu})$. Specifically, we know that every eigenvalue of $\bG^{\mu}$ lies in at least one of the discs $D(a_{kk}, R_{kk})$, where $a_{kk} = G^{\mu}_{kk}$ and $R_{kk} = \sum_{j \neq k} |G^{\mu}_{kj}|$. 
\begin{itemize}
    \item For $k \geq 3$, \eqref{lowerboundthetaeq1} implies $a_{kk} \asymp \frac{2k}{K + 1}$ and $R_{kk} = 0$. 
    \item For $k = 1$, 
    \begin{align*}
        a_{11} & = K \cdot \frac{\theta_1^2}{D_0(1, 1)} + \frac{9K}{16} \sum_{i=2}^{\frac{2n}{K(K + 1)}} \frac{\theta_i^2}{D_0^{\mu}(i,i)} \\
        R_{11} & = |G^{\mu}_{12}| = \frac{3K}{16} \sum_{j=2}^{\frac{2n}{K(K + 1)}} \frac{\theta_i^2}{D_0^{\mu}(i,i)}
    \end{align*}
    To bound $R_{11}$, we note that for sufficiently small $c_{12}$ (as defined in \eqref{defYmu}), $\max_{i \in \mathcal{C}_2} D_0^{\mu}(i,i) \leq 2\min_{i \in \mathcal{C}_2} D_0^{\mu}(i,i)$, since the probability matrix becomes arbitrarily close to the rank-one matrix $\mathbf{1}\mathbf{1}^T$ as $c_{12}$ decreases. In turn,
    $$
    \displaystyle R_{11} \leq  \frac{3K}{8} \sum_{i \in \mathcal{C}_{2}} \frac{\theta_i^2}{\max_{i \geq 2} D_0^{\mu}(i,i)} \leq \frac{3}{4} \cdot  a_{11}
    $$
    implying that any point within the disc $D(a_{11}, R_{11})$ is of order $a_{11}$.
    \item For $k = 2$, 
    \begin{align*}
        a_{22} & = K \sum_{i\in \mathcal{C}_{2}}\frac{\theta_i^2}{D_0^{\mu}(i,i)} + \frac{K}{16} \sum_{i=2}^{\frac{2n}{K(K + 1)}} \frac{\theta_i^2}{D_0^{\mu}(i,i)} \\
        R_{22} & = |G^{\mu}_{21}| = \frac{3K}{16} \sum_{i=2}^{\frac{2n}{K(K + 1)}} \frac{\theta_i^2}{D_0^{\mu}(i,i)}
    \end{align*}
    As in the $k = 1$ case, when $c_{12}$ is sufficiently small, we can show that any point within the disc $D(a_{22}, R_{22})$ is of order $a_{22}$.
\end{itemize}
In summary, any point within a Gergoshonin disk is of order $a_{kk}$ for some $k \in [K]$. Since $\min_{k \in [K]} a_{kk} \asymp \frac{2}{K + 1}$, in particular, $\lambda_{K}(\bG^{\mu})$ is of order $\frac{2}{K + 1}$. In turn, $\norm{(\bG^{\mu})^{-1}}$ is of order $K$, and since $K$ is fixed by assumption, it is in fact of constant order. 

To characterize the norm of $\bG^{\alpha}$, we use triangle inequality,
\small
\begin{align*}
    \norm{\bG^{\mu} - \bG^{\alpha}} & \leq \norm{\bG^{\mu} - K \cdot (\bPi^{\mu})^T \bTheta^{\alpha} (\bD^{\alpha})^{-1} \bTheta^{\alpha} \bPi^{\mu}} + \norm{K \cdot (\bPi^{\mu})^T \bTheta^{\alpha} (\bD^{\alpha})^{-1} \bTheta^{\alpha} \bPi^{\mu} - \bG^{\alpha}}  \numberthis \label{lowerboundPeq3}
\end{align*}
\normalsize

We bound each of the terms on the RHS of \eqref{lowerboundPeq3} separately. First off,
\small
\begin{align*}
    (\bPi^{\mu})^T \bTheta^{\alpha} (\bD^{\mu})^{-1} \bTheta^{\alpha} \bPi^{\mu} - \bG^{\mu} & =  K\diag\left(\frac{(\theta^{\alpha}_1)^2}{D_0^{\alpha}(1, 1)} - \frac{(\theta_1^{\mu})^2}{D_0^{\mu}(1,1)}, 0, \cdots, 0\right) \\
    & + K\Big(\sum_{i=2}^{\frac{2n}{K(K + 1)}} (\frac{(\theta^{\alpha}_i)^2}{D_0^{\alpha}(i, i)} - \frac{(\theta_i^{\mu})^2}{D_0^{\mu}(i,i)})\Big) (\begin{bmatrix}
        \frac{3}4 & \frac{1}{4} & 0 & \cdots & 0
    \end{bmatrix})^T\begin{bmatrix}
        \frac{3}4 & \frac{1}{4} & 0 & \cdots & 0
    \end{bmatrix} \\
    & +  K\diag\left(0, \sum_{i \in C_{2}} \frac{(\theta^{\alpha}_i)^2}{D_0^{\alpha}(i, i)} - \frac{(\theta^{\mu}_i)^2}{D_0^{\mu}(i, i)}, \cdots, \sum_{i \in C_{K}} \frac{(\theta^{\alpha}_i)^2}{D_0^{\alpha}(i, i)} - \frac{(\theta^{\mu}_i)^2}{D_0^{\mu}(i, i)} \right) \\
    & = K\diag\left(\frac{(\theta^{\alpha}_1)^2}{D_0^{\alpha}(1, 1)} - \frac{(\theta_1^{\mu})^2}{D_0^{\mu}(1,1)}, 0, \cdots, 0\right) \\
    & + K \sum_{i \in C_{2}} \left(\frac{(\theta^{\alpha}_i)^2}{D_0^{\alpha}(i, i)} - \frac{(\theta^{\mu}_i)^2}{D_0^{\mu}(i, i)}\right) \cdot \begin{bmatrix}
        0 & \frac{3}{16} & 0 & \cdots & 0 \\
        \frac{3}{16} & 0 & 0 & \cdots & 0 \\
        \vdots & \ddots & & & \vdots \\
        0 & \cdots & & & 0
    \end{bmatrix}\\
    & + K\diag\left(\frac{9}{16} \sum_{i=2}^{\frac{2n}{K(K + 1)}} (\frac{(\theta^{\alpha}_i)^2}{D_0^{\alpha}(i, i)} - \frac{(\theta_i^{\mu})^2}{D_0^{\mu}(i,i)}), \frac{1}{16} \sum_{i=2}^{\frac{2n}{K(K + 1)}} (\frac{(\theta^{\alpha}_i)^2}{D_0^{\alpha}(i, i)} - \frac{(\theta_i^{\mu})^2}{D_0^{\mu}(i,i)}), \cdots, 0\right) \\
    & + \small K\diag\left(0, \sum_{i \in C_{2}} \frac{(\theta^{\alpha}_i)^2}{D_0^{\alpha}(i, i)} - \frac{(\theta^{\mu}_i)^2}{D_0^{\mu}(i, i)}, \cdots, \sum_{i \in C_{K}} \frac{(\theta^{\alpha}_i)^2}{D_0^{\alpha}(i, i)} - \frac{(\theta^{\mu}_i)^2}{D_0^{\mu}(i, i)} \right) \numberthis \label{lowerboundPboundGalphaeq1}
\end{align*}
\normalsize

To bound the RHS of \eqref{lowerboundPboundGalphaeq1}, we characterize the differences $\theta^{\alpha}_i - \theta^{\mu}_i$ and $|D_0^{\alpha}(i, i) - D_0^{\mu}(i, i)|$ for $i \in [n]$. Since only the membership vectors of the first $\frac{2n}{K(K + 1)}$ nodes differ between $\bPi^{\alpha}$ and $\bPi^{\mu}$, only the first $\frac{2n}{K(K + 1)}$ nodes have $\Theta^{\alpha}(i, i) \neq \Theta^{\mu}(i, i)$. Furthermore, note that for any node $2 \leq i \leq \frac{2n}{K(K + 1)}$, its embedding $\bZ^{\alpha}_{i, :}$ intersects the two-simplex formed by the first two vertices $\bY^{\alpha}_{1, :}, \bY^{\alpha}_{2, :}$, so it belongs entirely to the first two communities, i.e. $\bpi^{\alpha}_i(k) \neq 0$ iff $k \leq 2$. In fact, $(\bpi^{\alpha}_i)^TY^{\alpha} = (\bpi^{\mu}_i)^TY^{\mu}$, so for all $2 \leq i \leq \frac{2n}{K(K + 1)}$,
\begin{align*}
    |\Theta^{\alpha}(i, i)| & = \norm{\bZ^{\alpha}_{i, :}} / \norm{(\bpi^{\alpha}_i)^T\bY^{\alpha}_{i, :}} \\
    & = \norm{\bZ^{\mu}_{i, :}} / \norm{(\bpi^{\mu}_i)^T\bY^{\mu}_{i, :}} \\
    & = |\Theta^{\mu}(i, i)|
\end{align*}
i.e. the degrees of all except the first node remain unchanged. For all $1 \leq i \leq \frac{2n}{K(K + 1)}$, routine computation reveals
\begin{align*}
|\bpi^{\alpha}_i(2) - \bpi^{\mu}_i(2)| \leq c\Delta_Y (c_{12} K \lambda_K)^{-1} 
\end{align*} 
for sufficiently small $c_0$ in the definition of $\Delta_Y$. From here, another routine computation shows
$$
\Theta^{\mu}(i, i) \leq \Theta^{\alpha}(i, i) \leq \Theta^{\mu}(i, i)(1 + \Delta_Y)
$$
In turn, $|\Theta^{\alpha}(i, i) - \Theta^{\mu}(i, i)| \ll \Theta^{\mu}(i, i)$. 
On the other hand, by construction, for $i \neq 1$, 
\begin{align*}
|D_0^{\alpha}(i, i) - D_0^{\mu}(i, i)| & = \sum_{j = 1}^n \left(\bZ^{\alpha}_{i, :}(\bZ^{\alpha}_{j, :})^T - \sum_{j = 1}^n \bZ^{\mu}_{i, :}(\bZ^{\alpha}_{j, :})^T\right) + \frac{\delta}{n} \sum_{i = 1}^n \sum_{j = 1}^n \left(\bZ^{\alpha}_{i, :}(\bZ^{\alpha}_{j, :})^T - \sum_{j = 1}^n \bZ^{\mu}_{i, :}(\bZ^{\alpha}_{j, :})^T\right)\\
& = \sum_{j = 1} \left(\bZ^{\alpha}_{i, :}(\bZ^{\alpha}_{j, :})^T - \sum_{j = 1}^n \bZ^{\mu}_{i, :}(\bZ^{\alpha}_{j, :})^T\right) + \frac{\delta}{n} \sum_{i = 1 \hbox{ or } j = 1} \left(\bZ^{\alpha}_{i, :}(\bZ^{\alpha}_{j, :})^T - \sum_{j = 1}^n \bZ^{\mu}_{i, :}(\bZ^{\alpha}_{j, :})^T\right) \\
& \leq c\Delta_Y \theta_1(\theta_i \wedge \otheta) \ll n \otheta (\theta_i \vee \otheta) = O(D(i, i))
\end{align*}
whereas for $i = 1$, similar computations show
\begin{align*}
    |D_0^{\alpha}(1, 1) - D_0^{\mu}(1, 1)| \leq cn\Delta_Y \otheta \theta_1
\end{align*}
As a result, for all $i \in [n]$,
\begin{align}
    |\frac{(\theta_i^{\alpha})^2}{D_0^{\alpha}(i, i)} - \frac{(\theta_i^{\mu})^2}{D^{\mu}_0(i, i)}| \lesssim \Delta_Y |\frac{(\theta_i^{\mu})^2}{D^{\mu}_0(i, i)}| \ll |\frac{(\theta_i^{\mu})^2}{D^{\mu}_0(i, i)}| \label{lowerboundPeq5}
\end{align}
implying each of the terms on the righthand side of \eqref{lowerboundPboundGalphaeq1} to be much smaller than $c \asymp |\bG|$, up to some factors of $K$. 

To bound the second RHS term of \eqref{lowerboundPeq3}, denote $\Gamma_{\Pi} = \bPi^{\alpha} - \bPi^{\mu}$. By triangle inequality,
\begin{align*}
    \norm{K \cdot (\bPi^{\alpha})^T \bTheta^{\alpha}(\bD^{\alpha})^{-1} \bTheta^{\alpha}\bPi^{\alpha} - \bG^{\alpha}} & \leq 2\norm{K(\Gamma_{\Pi})^T \bTheta^{\alpha}(\bD^{\alpha})^{-1} \bTheta^{\alpha}\bPi^{\mu}} + \norm{K(\Gamma_{\Pi})^T \bTheta^{\alpha}(\bD^{\alpha})^{-1} \bTheta^{\alpha}\Gamma_{\Pi}}
\end{align*} In order to bound the RHS of the above equation, we first note that it can be rewritten as
$$
2\norm{K(\Gamma_{\Pi})^T [\bTheta^{\alpha}(\bD^{\alpha})^{-1} \bTheta^{\alpha}]_{(-1, -1)}\Pi^{\mu}} + \norm{K(\Gamma_{\Pi})^T [\bTheta^{\alpha}(\bD^{\alpha})^{-1} \bTheta^{\alpha}]_{(-1, -1)}\Gamma_{\Pi}}
$$
where for a matrix $A$, $A_{(-1, -1)}$ denotes the matrix $A$ without the first entry. Removing the first entry is allowed because $\Gamma_{\Pi}$ has an empty first row. Also, we pick $c_{12}$ sufficiently small such that 
\begin{align*}
    |P^{\mu}|_{\min} \geq \frac{1}{2} \numberthis \label{Pboundedbelow}
\end{align*}
which in turn implies $D_0^{\alpha}(i, i) \geq \theta_i\sum_{j = 1}^n \theta_j|P^{\mu}|_{\min} = \frac{1}{2} \cdot n\theta_i\otheta$.
 With these simplifications, we now bound $(1)$ via the submultiplicty of norms:
\begin{align*}
    (1) & \leq 2K\norm{\Gamma_{\Pi}}_F\norm{[\bTheta^{\alpha}(\bD^{\alpha})^{-1} \bTheta^{\alpha}]_{(-1, -1)}}\norm{\bPi^{\mu}} + K\norm{\Gamma_{\Pi}}_F^2\norm{[\bTheta^{\alpha}(\bD^{\alpha})^{-1} \bTheta^{\alpha}]_{(-1, -1)}} \\
    & \leq 2K(\bpi^{\alpha}(2) \sqrt{2} \cdot \frac{2n}{K(K + 1)})(\max_{i \geq 2} \frac{(\theta_i^{\alpha})^2}{D_0^{\alpha}(i, i)}) + K(\bpi^{\alpha}(2) \sqrt{2} \cdot \frac{2n}{K(K + 1)})^2(\max_{i \geq 2} \frac{(\theta_i^{\alpha})^2}{D_0^{\alpha}(i, i)}) \\
    & \leq cK\Delta(c_{12} K \lambda_K)^{-1} + cK\Delta^2 (c_{12} K \lambda_K)^{-2}
\end{align*}
where $c$ is a constant. By taking $c_0$ in $\Delta$ small enough, moreover, we can bound $\Delta(c_{12} K \lambda_K)^{-1}$ and thus the RHS by an arbitrarily small constant. Thus, $\norm{\bG^{\alpha}}$ also satisfies the first two relations of Assumption \ref{regconds1}(a).

By the construction of $\bPi^{\mu}$ and $\bPi^{\alpha}$, the last relation of Assumption \ref{regconds1}(a) holds for both $\bH^{\mu}$ and $\bH^{\alpha}$. 

By our choice of $\bP^{\mu}$, it follows that $\lambda_1(\bP^{\mu}\bG^{\mu}) \asymp K$ and $\lambda_k(\bP^{\mu}\bG^{\mu}) \asymp K\lambda_K$ for all $2 \leq k \leq K$. Due to our choice of $\bPi^{\mu}$, furthermore, the eigenvalues are well-separated. Thus, there exist $c_{10}$ sufficiently small in the definition of $\tilde{\cal Q}_{n, ab, c_{10}}$ such that Assumptions \ref{regconds1}(b) and Assumption \ref{regconds2}(a) hold for $\bH^{\mu}$. As for $\bH^{\alpha}$, 

\setstretch{1.1}

\begin{align*}
    \lambda_1(\bP^{\alpha}\bG^{\alpha}) & = K \lambda_1(\bL_0^{\alpha}) \geq K(\lambda_1(\bL_0^{\mu}) - \norm{\bL_0^{\mu} - \bL_0^{\alpha}}) \\
    |\lambda_k(\bP^{\mu}\bG^{\mu})| + K\norm{\bL_0^{\mu} - \bL_0^{\alpha}} & \geq |\lambda_k(\bP^{\alpha}\bG^{\alpha})| \geq  |\lambda_k(\bP^{\mu}\bG^{\mu})| - K\norm{\bL_0^{\mu} - \bL_0^{\alpha}} \numberthis \label{lowerboundPeq6}
\end{align*}

\setstretch{1.0}

Since $\bH^{\alpha} = \bZ^{\alpha}(\bZ^{\alpha})^T$ and $\bH^{\mu} = \bZ^{\mu}(\bZ^{\mu})^T $ differ only in the first row $\bZ_{1, :}$, we expect their Laplacian matrices to be close. Indeed, some routine computation shows $\norm{\bL_0^{\mu} - \bL_0^{\alpha}} \leq \Delta_Y \sqrt{\frac{1}{c_{12} K \lambda_K n^2\otheta^4}}$. By choosing $c_0$ small enough in $\Delta_Y$, we can force $\Delta_Y \sqrt{\frac{1}{c_{12} K \lambda_K n^2\otheta^4}} \ll c_{12} K \lambda_K \leq K^{-1} \cdot |\lambda_k(\bP^{\mu}\bG^{\mu})|$, which, in conjunction with \eqref{lowerboundPeq6}, implies Assumptions \ref{regconds1}(b) and Assumption \ref{regconds2}(a) also hold for $\bH^{\alpha}$.

Next, we show Assumption \ref{regconds1}(c) holds. By Perron's theorem, the first singular vector, $\boldeta_1^{\mu}$, of $\bP^{\mu}\bG^{\mu}$ is positive. Since $\bP^{\mu}\bG^{\mu}$ has entries positive and all of the same order, $\boldeta_1^{\mu}$ satisfies
$$
\frac{\min_{k} \boldeta_1^{\mu}(k)}{\max_{K} \boldeta_1^{\mu}(k)} \geq \frac{\min_{i,j}\bP^{\mu}\bG^{\mu}(i,j) \sum_{k}  \boldeta_1^{\mu}(k) }{ \max_{i,j}\bP^{\mu}\bG^{\mu}(i,j)\sum_{k}  \boldeta_1^{\mu}(k) } >c
$$
For $\bP^{\alpha}\bG^{\alpha}$, sin-theta theorem implies
\begin{align*}
\Vert \boldeta^{\alpha}_1 -  \boldeta_1^{\mu}\Vert = \sqrt{2 |(\boldeta_1^{\mu})^T  \boldeta^{\alpha}_1  - 1|} \leq C \Vert \bP^{\mu}\bG^{\mu} - \bP^{\alpha}\bG^{\alpha} \Vert 
\end{align*}
By triangle inequality,
$$
\norm{ \bP^{\mu}\bG^{\mu} - \bP^{\alpha}\bG^{\alpha}} \leq \norm{\bP^{\mu}}\norm{\bG^{\mu} - \bG^{\alpha}} + \norm{\bP^{\mu} - \bP^{\alpha}}\norm{G^{\alpha}}
$$
Since only the first row of $Y^{\mu}$ and $Y^{\alpha}$ differ, we expect $\norm{\bP^{\mu}- \bP^{\alpha}}$ to be small. Indeed, some routine computation shows 
\begin{align*}
    \norm{\bP^{\mu}- \bP^{\alpha}} \lesssim \Delta_Y \numberthis \label{lowerboundPeqPdev}
\end{align*} Previously, we bounded $\norm{\bG^{\mu} - \bG^{\alpha}}$ by $\Delta(c_{12} K \lambda_K)^{-1}$ in \eqref{lowerboundPeq3}. As a result, 
\begin{align*}
    \norm{\bP^{\mu}\bG^{\mu} -  \bP^{\mu} \bG^{\alpha}} & \leq \norm{\bP^{\mu}} \cdot \Delta_Y(c_{12} K \lambda_K)^{-1} + \Delta_Y\sqrt{\frac{K}{\lambda_K}}\norm{G^{\alpha}} \\
    & \lesssim \Delta_Y \lambda_K^{-1} \numberthis \label{lowerboundsPeq15}
\end{align*} up to factors of $K$. In turn, for sufficiently small $c_0$ in $\Delta_Y$,
\begin{align*}
\frac{\min_{k}  \boldeta^{(j)}_1(k)} {\max_{k}  \boldeta^{(j)}_1(k)}> \frac{\min_{k}  \boldeta_1^{\mu}(k) - \Delta_Y \lambda_K^{-1}}{\max_{k}  \boldeta_1^{\mu}(k) + \Delta_Y \lambda_K^{-1}} > c 
\end{align*}
As for Assumption \ref{regconds1}(d), $\bP^{\mu} = \bY^{\mu}(Y^{\mu})^T = (1 - 2c_{12}K\lambda_K)\mathbf{1}\mathbf{1}^T + 2c_{12}K\lambda_K\mathbf{I}_{K \times K}$, implying $\bP^{\mu}$ has unit diagonals and is non-singular. As for $\bP^{\alpha}$ (which has unit diagonals by construction), taking $c_0$ in $\Delta_Y$ small enough ensures 
\begin{align*}
    \lambda_K(\bP^{\alpha}) & \geq \lambda_K(\bP^{\mu}) - \norm{\bP^{\alpha} - \bP^{\mu}} \\
    & \asymp K\lambda_K
\end{align*}
where the last line follows from \eqref{lowerboundPeqPdev}, i.e. $P^{\alpha}$ is also non-singular.

Now, we show the remaining points in Assumption \ref{regconds2}.
\begin{itemize}
    \item \underline{Assumption \ref{regconds2}(b):} $\lambda_K(\bL_0) \asymp \lambda_K \gg \sqrt{\frac{\log n}{n\otheta^2}}$, as $\otheta \geq \sqrt{\frac{\log n}{n}}$ by assumption. Moreover, 
    $$
    \lambda_K(\bL_0^{\alpha}) \geq \lambda_K(\bL_0) - \norm{\bL_0^{\alpha} - \bL_0^{\mu}} \geq \lambda_K(\bL_0) - \Delta_Y (K\lambda_K)^{-1}
$$
where the last inequality follows from the bound in equation \eqref{lowerboundsPeq15}. Using the fact that $\lambda_k \gg \sqrt{\frac{\log n}{n\otheta^2}}$ and $C\sqrt{\frac{\log n}{n}} \leq \theta_i \leq C$ for all $i \in [n]$, we can bound $\Delta_Y (K\lambda_K)^{-1} \leq C\sqrt{\frac{n}{\log n}}$. For sufficiently small $c_0$ in $\Delta_Y$, therefore, $\Delta_Y \lambda_K^{-1} \ll \lambda_K(\bL_0)$, implying Assumption \ref{regconds2}(b) to also hold for $\mathcal{Q}^{\alpha}$. 
    \item \underline{Assumption \ref{regconds2}(c):} By definition,
    \begin{align*}
        \bF^{\mu} & = K\diag(\sum_{i \in C_{1}} \frac{\theta_i^2}{\norm{\btheta}^2}, \cdots, \sum_{i \in C_{K}} \frac{\theta_i^2}{\norm{\btheta}^2})
    \end{align*}
    By a similar Gershgorin Circle argument as that used to show $\norm{\bG^{\mu}}, \norm{(\bG^{\mu})^{-1}} \leq C$, we can also show  $\norm{\bF^{\mu}}, \norm{(\bF^{\mu})^{-1}} \leq C$. Likewise, by a similar perturbation argument as to that used to show $\norm{\bG^{\alpha}}, \norm{(\bG^{\alpha})^{-1}} \leq C$, we can also show $\norm{\bF^{\alpha}}, \norm{(\bF^{\alpha})^{-1}} \leq C$, so we omit the argument.
\item \underline{Assumption \ref{regconds2}(d-f):} Assumption \ref{regconds2}(d) and Assumption \ref{regconds2}(f) hold by the assumptions of Theorem \ref{lowerboundsP}. Since $\bP^{\mu} = \bY^{\mu}(\bY^{\mu})^T$, $P_{ab}^{\mu} = \bY^{\mu}_{a, :}(\bY^{\mu}_{b, :})^T \leq \norm{\bY^{\mu}_{a, :}}\norm{\bY^{\mu}_{b, :}} = 1$ for all $a, b \in [n]$, implying the entries of $\bP^{\mu}$ are bounded. By a similar dot product argument, it follows that the entries of $\bP^{\alpha}$ are bounded by $1$ as well.
\end{itemize}
In summary, both $\bH^{\mu}$ and $\bH^{\alpha}$ belong to $\tilde {\cal Q}_{n, ab, c_{10}}(\sigma)$ for a sufficiently small constant $c_{10}$.

We now turn to showing the lower bound. First, we analyze the difference between the probability matrices $\bP^{\mu}$ and $\bP^{\alpha}$:
\begin{align*}
    P^{\mu}_{ab} - P^{\alpha}_{ab} & = \sum_{k = 1}^K (Y^{\alpha}_{1, k}Y^{\alpha}_{2, k} - Y^{0}_{1, k}Y^{0}_{2, k}) = \sum_{k = 1}^2 (Y^{\alpha}_{1, k}Y^{\alpha}_{2, k} - Y^{0}_{1, k}Y^{0}_{2, k})
\end{align*}
Each of the summands in the rightmost term are strictly positive, so we can lower bound by
$$
\sum_{k = 1}^2 (Y^{\alpha}_{1, k}Y^{\alpha}_{2, k} - Y^{0}_{1, k}Y^{0}_{2, k}) \geq \Delta_Y
$$
for sufficiently small $c_0$ in the definition of $\Delta_Y$.

Finally, we show the KL divergence between the null and alternative hypotheses are bounded by a constant. Since only the first row of $\bZ$ differs between the null and alternative, only terms in the first row and column of $\bH$ differ. In turn, the KL divergence can be expressed as
\begin{align*}
KL(\mathcal{Q}^\alpha,\mathcal{Q}^{\mu})& =2\sum_{i = 1, j \neq 1}H^{\alpha}_{ij}\log(\frac{H^{\alpha}_{ij}}{H^{\mu}_{ij}}) + \big(1- H^{\alpha}_{ij}\big )  \log\frac{1- H^{\alpha}_{ij}}{1- H^{\mu}_{ij}}
\numberthis \label{lowerboundsPeq9}
\end{align*} 
Thus, we attempt to bound the difference $|\bH^{\alpha}_{1j} - \bH^{\mu}_{1j}|$ between terms.
\begin{itemize}[topsep=1pt, itemsep=1pt]
    \item \underline{For $j$ among the first $\frac{2n}{K(K + 1)}$ nodes:} We first note that $Y_{12}^{\alpha} > Y_{12}^{\mu}$ and
$$
Y_{12}^{\alpha} \leq \sqrt{\frac{c_{12} K\lambda_K}{2} + 2\Delta_Y} \leq \sqrt{\frac{c_{12} K\lambda_K}{2}}(1 + \frac{4\Delta_Y}{c_{12} K \lambda_K})
$$
implying
\begin{align*}
    |Y_{12}^{\alpha} - Y_{12}^{\mu}| \leq \frac{\Delta_Y\sqrt{8}}{\sqrt{c_{12} K \lambda_K}} \numberthis \label{lowerboundsPeq11}
\end{align*}
This allows us to establish
\begin{align*}
 |H^{\alpha}_{1j} - H^{\mu}_{1j}| & = \theta_1(|Y_{11}^{\alpha}X_{11}^{\mu} - Y_{11}^{\mu}X_{11}^{\mu} + Y_{12}^{\alpha}\bZ_{12}^{\mu} - Y_{12}^{\mu}\bZ_{12}^{\mu}|) \\
 & \leq \theta_1(|Y_{11}^{\alpha} - Y_{11}^{\mu}|X_{11}^{\mu} + |Y_{12}^{\alpha} - Y_{12}^{\mu}|\bZ_{12}^{\mu}) \\
 & \leq \theta_1(2\theta_j\Delta_Y + 2\theta_j\Delta_Y) = 4\theta_1\theta_j \\
 & \leq 8H_{1j} \Delta_Y \numberthis \label{lowerboundsPeq10}
    \end{align*}
    where the third line follows from plugging in the bound in \eqref{lowerboundsPeq11}, and the fourth line by \eqref{Pboundedbelow}. 
    \item \underline{For $j$ in $C_2$:} The argument is similar, so we omit it. 
    \item \underline{For $j$ in communities $\{C_3, \cdots C_K\}$:} First off,
\begin{align*}
 H_{1j} & = \theta_1\theta_j(1 - \lambda_k - \Delta_Y \sqrt{1 - c_{12} (K - 1)\lambda_K}) \\
    & \geq H_{ij}(1 - \frac{\Delta_Y}{1 - \lambda_k}) \numberthis \label{lowerboundsPeq12}
\end{align*}
We can bound the KL divergence by the following bound on the summand in \eqref{lowerboundsPeq9}:
\begin{align*}
    H_{ij}^{\alpha} \log(\frac{H_{ij}^{\alpha}}{H_{ij}^{\mu}}) + (1 - H_{ij}^{\alpha})(\log(\frac{1 - H_{ij}^{\alpha}}{1 - H_{ij}^{\mu}}) & = H_{ij}^{\alpha} \log(1 + \frac{H_{ij}^{\alpha} - H_{ij}^{\mu}}{H_{ij}^{\mu}})\\
    & + (1 - H_{ij}^{\alpha}) \log(1 - \frac{H_{ij}^{\mu} - H_{ij}^{\alpha}}{1 - H_{ij}^{\mu}}) \\
    & \leq H_{ij}^{\alpha}(\frac{H_{ij}^{\alpha} - H_{ij}^{\mu}}{H_{ij}^{\mu}}) - (1 - H_{ij}^{\alpha})(\frac{H_{ij}^{\mu} - H_{ij}^{\alpha}}{1 - H_{ij}^{\mu}}) \\
    & = \frac{(H_{ij}^{\alpha} - H_{ij}^{\mu})^2}{H_{ij}^{\mu}(1 - H_{ij}^{\mu})} \numberthis \label{KLdvgbound}
\end{align*}
By our choice of $c_{0}$ in \eqref{Pboundedbelow}, $1 - H_{ij}^{\mu} \geq C$ for some constant $C$. 
\end{itemize} 
Alongside the bounds in equations \eqref{lowerboundsPeq10} and \eqref{lowerboundsPeq12}, we therefore obtain
\begin{align*}
    KL(\mathcal{Q}^\alpha,\mathcal{Q}^{\mu}) & \leq C\sum_{i = 1, 1 \leq j \leq n} \frac{(H_{ij}^{\alpha} - H_{ij}^{\mu})^2}{H_{ij}^{\mu}} \\
    & \leq C \sum_{i = 1, 1 \leq j \leq n} \frac{H_{ij}\Delta_Y^2}{(1 - c_{12} K \lambda_K)^2} \\
    & \leq C'\sum_{i = 1, 1 \leq j \leq n} H_{ij}\Delta_Y^2 \\
    & \leq C''n\theta_{1}\otheta \Delta_Y^2
\end{align*}
for some fixed constants $C, C', C''$, where the penultimate line follows from the assumption that $c_{12} K \lambda_K \leq (1 - c_0)K\lambda_1 \leq 1 - c_0$. Since $\Delta_Y \lesssim \frac{1}{\sqrt{n\theta_{1}\otheta}}$ by construction, the KL divergence is therefore bounded by a constant, thus establishing the lower bound.
\subsubsection{Proof of Lemma \ref{reparamlem}} \label{reparamlempf}
    Denote the DCMM model $\mathcal{Q} = (\bTheta', \bPi', \bY\bY^T)$. It is easy to check that $(\bTheta', \bPi', \bY\bY^T)$ satisfy Assumptions \ref{regconds1} \& \ref{regconds2} and that \eqref{lowerboundsPproofeq1} holds. Having shown that $\mathcal{Q}$ is a valid DCMM model, the fact that its probability matrix $\bP$ is non-singular and has unit diagonals shows it is unique (as proven in Proposition A.1 of \citet{jin2022mixed}).

\subsection{Proof of Theorem \ref{lowerboundsTheta}}\label{lowerboundsThetaproof}
In the first part of this proof, we show $\inf_{\hat{\bTheta}}\sup_{\bTheta \in \tilde {\cal R}_{n, ab, c_{11}}} |\hat{\Theta}(i, i) - \Theta(i, i)| \geq C\sqrt{\frac{\log n}{n\otheta^2}}(\frac{\theta_i}{\otheta})$ using one parameter set; in the second part, we show $\inf_{\hat{\bTheta}}\sup_{\bTheta \in \tilde {\cal R}_{n, ab, c_{11}}} |\hat{\Theta}(i, i) - \Theta(i, i)| \geq C\sqrt{\frac{\log n}{n\otheta^2}}(\sqrt{\frac{\theta_i}{\otheta}})$ using a different parameter set. Together, these two lower bounds will therefore imply the desired minimax rate.

As in the proof of Lemma \ref{lowerboundsP}, we seek to construct two parameter sets $\mathcal{Q}^{\mu} = (\bTheta^{\mu}, \bPi^{\mu}, \bP^{\mu})$ and $\mathcal{Q}^{\alpha} = (\bTheta^{\alpha}, \bPi^\alpha, \bP^\alpha)$ such that $|\theta_i^{\alpha} - \theta_i^{\mu}| \geq C\sqrt{\frac{\log n}{n\otheta^2}}(\frac{\theta_i}{\otheta})$ and $KL(\mathcal{Q}^{\alpha}, \mathcal{Q}^{\mu}) \leq C'$ for constants $C, C'$, from which standard lower bound techniques (e.g. Lemma 2.9 of \cite{Tsybakov_2009}) will imply the desired lower bound of $n^{-\frac{1}{2}}$. For the sake of convenience, we assume:
\begin{itemize}[itemsep=1pt, topsep=2pt]
    \item $\theta_i$ belongs to the first community and is the second node within the community.
\end{itemize}
As in Lemma \ref{lowerboundsP}, we use Lemma \ref{reparamlem} to define our null and alternative hypothesis. Consider the null hypothesis $\mathcal{Q}^{\mu} = (\bZ^{\mu}, \bY^{\mu})$, where
\begin{align*}
    \bY^{\mu} = \begin{bmatrix}
    \sqrt{1 - c_{13}(K - 1) } & \sqrt{\frac{c_{13}K}{2}} & \sqrt{\frac{c_{13}K}{6}} & \cdots &  \sqrt{\frac{c_{13}K}{(K - 1)^2 + (K - 1)}} \\
     \sqrt{1 - c_{13}(K - 1) } & -\sqrt{\frac{c_{13}K}{2}} & \sqrt{\frac{c_{13}K}{6}} & \cdots &  \sqrt{\frac{c_{13}K}{(K - 1)^2 + (K - 1)}} \\
      \sqrt{1 - c_{13}(K - 1) } &  0 & -2\sqrt{\frac{c_{13}K}{6}} & \cdots &  \sqrt{\frac{c_{13}K}{(K - 1)^2 + (K - 1)}} \\
      \vdots & \vdots & \vdots & \vdots & \vdots \\
       \sqrt{1 - c_{13}(K - 1) } & 0 & 0 & 0 & -(K - 1)\sqrt{\frac{c_{13}K}{(K - 1)^2 + (K - 1)}}
\end{bmatrix} \numberthis \label{thetaYmu}
\end{align*}
and $c_{13}$ is a sufficiently small constant to be determined later. Also, define
\begin{align*}
    \bZ^{\mu} = \bTheta^{\mu} \begin{bmatrix}
\mathbf{e}_1 \\
\mathbf{1}_{\frac{2n}{K(K + 1)} - 1}^T \begin{bmatrix}
    \frac{3}{4} & \frac{1}{4} & 0 & \cdots & 0
\end{bmatrix} \\
    \mathbf{1}_{\frac{4n}{K(K + 1)}} \mathbf{e}_2^T \\
    \vdots \\
    \mathbf{1}_{\frac{2Kn}{K(K + 1)}} \mathbf{e}_K^T
\end{bmatrix} \bY^{\mu} \numberthis \label{thetaXmu}
\end{align*}
where $\bTheta^{\mu}$ is defined such that
$$
\Theta^{\mu}(j, j) = \begin{cases}
    \theta_{i} & \hbox{ if $j = 2$}\\
    \check{\theta} & \hbox{ otherwise}
\end{cases}
$$
and $\sum_{i = 1}^n \Theta^{\mu}(i, i) = n\otheta$, i.e. all nodes except the second have the same degree $\check{\theta}$. Regarding the alternative hypothesis, we define
$$
\bY^{\alpha} = \begin{bmatrix}
    \sqrt{1 - c_{13}(K - 1) } - \Delta_Y & \sqrt{\frac{\lambda_K}{2} + 2 \Delta_Y \sqrt{1 - c_{13}(K - 1) } - \Delta_Y^2} & \sqrt{\frac{c_{13}K}{6}} & \cdots &  \sqrt{\frac{c_{13}K}{(K - 1)^2 + (K - 1)}} \\
     & \bY^{\mu}_{-1, :} & & 
\end{bmatrix}
$$
where $$
\Delta_Y := c_0(n\otheta^2)^{-\frac{1}{2}}
$$
for a properly small $c_0$ to be chosen later, and
$$
\bZ^{\alpha} = \begin{bmatrix}
     \check{\theta} * \bY^{\alpha}_{1, :} \\
     \bZ^{\mu}_{-1, :}
\end{bmatrix}
$$
which is in fact the same definition as in our lower bounds for $\bP$. (The only formal difference is in the definition of $\Delta_Y$, but since every pure node has degree of order $\otheta$, the two definitions are actually identical.) By similar arguments as to the proof of Theorem \ref{lowerboundsP}, then, we can show that $\mathcal{Q}^{\mu}$ and $\mathcal{Q}^{\alpha}$ both i) satisfy Assumptions \ref{regconds1} and \ref{regconds2} for $c_{13}$ sufficiently small and ii) 
lie in ${\cal R}_{n, i, c_{11}}(\tau)$ for $c_{11}$ sufficiently small. Moreover, one can show that $KL(\mathcal{Q}^{\alpha}, \mathcal{Q}^{\mu}) \leq C$ by a similar argument, so it suffices to show the lower bound for the difference $|\theta^{\alpha}_i - \theta^{\mu}_i|$.

To this end, we first note that $\theta^{\alpha}_i = \frac{\norm{\bZ^{\alpha}_{2, :}}}{\norm{(\bpi^{\alpha}_2)^T \bY^{\alpha}}}$ by definition. Denote $\bV = (\bpi^{\alpha}_2)^T \bY^{\alpha} \in \mathbb{R}^{1 \times K}$. For sufficiently small $c_0$ in $\Delta_Y$, furthermore, $\bpi^{\alpha}_2(1) \geq \frac{1}{2}$, meaning $V(1) \leq \frac{1}{2}Y_{11}^{\alpha} + \frac{1}{2}Y_{21}^{\alpha} \leq \sqrt{1 - c_{13}(K - 1)} - \frac{1}{2} \cdot \Delta_Y$. Since $\bV$ is a scalar multiple of $\bZ^{\alpha}_{2, :}$ by definition, it follows that
\begin{align*}
    \frac{\norm{\bZ^{\alpha}_{2, :}}}{\norm{(\bpi^{\alpha}_2)^T \bY^{\alpha}}} & = \frac{Z^{\alpha}_{21}}{V(1)} \\
    & \leq (\theta_2 \sqrt{1 - c_{13}(K - 1)}) \cdot (\sqrt{1 - c_{13}(K - 1)} - \frac{1}{2} \cdot \Delta_Y)^{-1} \\
    & \geq \theta_2(1 + \frac{1}{2} \cdot \Delta_Y)
\end{align*}
In turn, $\theta^{\alpha}_i - \theta^{\mu}_i \geq C\frac{\theta_2}{\Delta_Y} = C'\theta_i(n\otheta^2)^{-\frac{1}{2}}$. Furthermore, $\theta^{\alpha}_i > \theta^{\mu}_i$ by construction, meaning
$$
|\theta^{\alpha}_i - \theta^{\mu}_i| \geq C\sqrt{\frac{\log n}{n\otheta^2}}(\frac{\theta_i}{\otheta})
$$ as desired.

Moving onto the second half on the proof, we seek to construct two parameter sets $\mathcal{Q}^{\mu} = (\bTheta^{\mu}, \bPi^{\mu}, \bP^{\mu})$ and $\mathcal{Q}^{\alpha} = (\bTheta^{\alpha}, \bPi^\alpha, \bP^\alpha)$ such that $|\theta_i^{\alpha} - \theta_i^{\mu}| \geq C\sqrt{\frac{\log n}{n\otheta^2}}(\sqrt{\frac{\theta_i}{\otheta}})$ and $KL(\mathcal{Q}^{\alpha}, \mathcal{Q}^{\mu}) \leq C'$. Our construction in this case is far simpler: like before, we take $\bY^{\mu}, \bZ^{\mu}$ as in \eqref{thetaYmu}, \eqref{thetaXmu}. For the alternative hypothesis, however, we simply perturb the degree matrix
$$
\Theta^{\alpha}(j, j) = \begin{cases}
\theta_i(1 + c_{14}\sqrt{\frac{\theta_i}{\otheta}}) & \hbox{ if $j = 2$} \\
    \Theta^{\mu}(j, j) & \hbox{ otherwise}
\end{cases}
$$
and keep $\bPi^{\alpha} = \bPi^{\mu}, \bY^{\alpha} = \bY^{\mu}$ constant. By the same arguments as in the lower bound for $\bP$, we can show $\mathcal{Q}^{\mu}$ satisfies Assumptions \ref{regconds1} and \ref{regconds2}. Since all pure nodes have degree of order $\otheta$, furthermore, Assumption \ref{regconds4} is also satisfied. Lastly, one can show $\mathcal{Q}^{\alpha}$ also satisfies Assumptions \ref{regconds1} and \ref{regconds2} by straightforward pertubation arguments. 

Now, we turn to bounding the KL divergence. Since only $\theta_2$ was perturbed, the only terms that differ between $\bH^{\mu}$ and $\bH^{\alpha}$ lie in the second row and column. In turn, the KL divergence can be expressed as
\begin{align*}
KL(\mathcal{Q}^\alpha,\mathcal{Q}^{\mu})& =2\sum_{i = 2, j \neq 1}H^{\alpha}_{ij}\log(\frac{H^{\alpha}_{ij}}{H^{\mu}_{ij}}) + \big(1- H^{\alpha}_{ij}\big )  \log\frac{1- H^{\alpha}_{ij}}{1- H^{\mu}_{ij}}
\end{align*} 
By \eqref{KLdvgbound}, this expression can be bounded by 
\begin{align*}
   KL(\mathcal{Q}^\alpha,\mathcal{Q}^{\mu}) \leq 2\sum_{i = 2, j \neq 1} \frac{(H_{ij}^{\alpha} - H_{ij}^{\mu})^2}{H_{ij}^{\mu}(1 - H_{ij}^{\mu})} \numberthis \label{lowerboundthetaeq16}
\end{align*}
By construction, $|H_{ij}^{\alpha} - H_{ij}^{\mu}| = (c_{14}\sqrt{\frac{1}{n\otheta\theta_i}})H_{ij}^{\mu}$. Plugging this bound into \eqref{lowerboundthetaeq16}, and also the facts that $1 - H_{ij}^{\mu} > \frac{1}{2}$ (from choosing ${c_{11}}$ small enough) and $\sum_{i = 2, j \neq 1} H_{ij}^{\mu} \leq n\otheta\theta_i$ yield
\begin{align*}
   KL(\mathcal{Q}^\alpha,\mathcal{Q}^{\mu}) & \leq Cc_{14}^2(n\otheta\theta_i)^{-1}\sum_{i = 2, j \neq 1} H_{ij}^{\mu} \\
   & \leq Cc_{14}^2 \numberthis \label{lowerboundthetaeq17}
\end{align*}
for a constant $C$ independent of $c_{14}$. By picking $c_{14}$ properly small, the KL divergence can therefore be bounded by an arbitrarily small constant. Thus, $\inf_{\hat{\bTheta}}\sup_{\bTheta \in \tilde {\cal R}_{n, ab, c_{11}}} |\hat{\Theta}(i, i) - \Theta(i, i)| \geq C\sqrt{\frac{\log n}{n\otheta^2}}(\sqrt{\frac{\theta_i}{\otheta}})$.

In conjunction with the first half of the proof, we obtain the desired lower bound.

\section{Proof of Theorem \ref{improveddevbounds}} \label{improveddevboundsproof}
In order to improve the eigenvector deviation bounds in \citet{ke2022optimal}, we exploit Assumption \ref{regconds2}(a) (the well-separated eigenvalues condition). Specifically, we use contour integrals to derive expansions for the eigenvectors, inspired by \cite{fan2022asymptotic} and \cite{bhattacharya2023inferences}. First, we define some relevant notation, starting with second-order expansions of the degree and Laplacian matrices. Recall that we define $\mathcal{A}_1$ to be the event that Lemma \ref{basicresults} holds.
\subsubsection{Notation} \label{sectionexpansiondegrees}
\begin{lem}\label{expansiondegrees}
We have the following expansion
\begin{align}
    \boldsymbol{\obsregdegreehalf} = \boldsymbol{\trueregdegreehalf} + \bDelta_{D, 1} + \bDelta_{D, 2} + \phi_D \label{expansiondegreeseq0}
\end{align}
where $\bDelta_{D, 1}, \bDelta_{D, 2}, \phi_D$ are diagonal matrices obtained by Taylor expanding each entry of $\trueregdegree$. In particular,
\begin{align*}
    \Delta_{D, 1}(i, i) & := -\frac{1}{2} \cdot \frac{(\sum_{j = 1}^n W_{ij} + \delta / n \sum_{i = 1}^n \sum_{j = 1}^n W_{ij})}{D_0(i, i)^{\frac{3}{2}}} \numberthis \label{expansiondegreeseq1} \\ 
    \Delta_{D, 2}(i, i) & := \frac{3}{4} \cdot \frac{(\sum_{j = 1}^n W_{ij} + \delta / n \sum_{i = 1}^n \sum_{j = 1}^n W_{ij})^2}{D_0(i, i)^{\frac{5}{2}}} \quad \forall \; i \in [n] \numberthis \label{expansiondegreeseq2}
\end{align*}
\end{lem}
Denote an `order 0' term to be either of $\boldsymbol{\trueregdegreehalf}, \boldsymbol{\obsregdegreehalf}$; an `order 1' term to be either of $\bDelta_{D, 1}, \bR_{D, 1} := \boldsymbol{\obsregdegreehalf} - \boldsymbol{\trueregdegreehalf}$; an `order 2' term to be either of $\bDelta_{D, 2}, \bR_{D, 2} := \boldsymbol{\obsregdegreehalf} - \boldsymbol{\trueregdegreehalf} - \bDelta_{D, 1}$; an `order 3' term to be $\bphi_D$. Under event $\mathcal{A}_1$, we can readily bound the magnitude of the various degree orders.

\begin{lem}\label{basicexpansiondegreeresults}
    Define $A = \sum_{i \leq n} \theta_i^2 N(i, i)^2$, and let $i \leq n$ be arbitrary. Under event $\mathcal{A}_1$, the following bounds hold:
    \begin{center}
        \begin{tabular}{ c|c|c } 
N & $|N(i, i)|$ & $A$ \\ \hline
 $\hbox{order 0 term}$ & $(\min\{\frac{1}{n\otheta\theta_i}, \frac{1}{n\otheta^2}\} \})^{\frac{1}{2}}$& $1$\\
 $\hbox{order 1 term}$ & $\sqrt{\log n} \cdot \min\{\frac{1}{n\otheta\theta_i}, \frac{1}{n\otheta^2}\}$ & $\frac{\log n}{n\otheta^2}$ \\
 $\hbox{order 2 term}$ & $\log n \cdot (\min\{\frac{1}{n\otheta\theta_i}, \frac{1}{n\otheta^2}\})^{\frac{3}{2}}$ & $\frac{\log^2 n}{n^2\otheta^4}$ \\
  $\hbox{order 3 term}$ & $\sqrt{\log^3 n}(\min\{\frac{1}{n\otheta\theta_i}, \frac{1}{n\otheta^2}\})^{2}$ & $\frac{\log^3 n}{n^3\otheta^6}$
\end{tabular}
    \end{center}
\end{lem}

The proofs of Lemmas \ref{expansiondegrees} and \ref{basicexpansiondegreeresults} are deferred to the end of Section \ref{improveddevbounds}. Now, we can define the second-order expansion of the Laplacian.
\\~\\
\textbf{Definition 1.} \; We have the expansion $\hat{\bL} = \bL_0 + \bE_1 + \bE_2 + \phi_E$ where 
\begin{align}
    \bE_1 & := \boldsymbol{\trueregdegreehalf} \bW\boldsymbol{\trueregdegreehalf} + \bDelta_{D, 1} \bH \boldsymbol{\trueregdegreehalf} + \boldsymbol{\trueregdegreehalf} \bH \bDelta_{D, 1} \\
    \bE_2 & := \bDelta_{D, 2} \bH \boldsymbol{\trueregdegreehalf} + \bDelta_{D, 1}W \boldsymbol{\trueregdegreehalf} + \bDelta_{D, 1} \bH \bDelta_{D, 1} + \boldsymbol{\trueregdegreehalf} \bW \bDelta_{D, 1} + \boldsymbol{\trueregdegreehalf} \bH \bDelta_{D, 2}
\end{align}
For use in later proofs, we also define
\begin{align}
    \bDelta_{D, 1} & = \bDelta_{D, 1} + \bDelta_{D, 2} + \phi_D, \bDelta_{D, 2} = \bDelta_{D, 2} + \phi_D \\
    \bE' & = \bE_1 + \bE_2 \label{expansiondefinitions}
\end{align}
\subsection{Proof of Theorem \ref{improveddevbounds}}
Define $r := \frac{|\lambda_k|-|\lambda_{k + 1}|}{2}$ and let $\mathcal{C}_1$ be the circular contour around $\lambda_k$ with radius $r$. Then $\lambda_k$ is the only eigenvalue of $\bL_0$ that is inside $\mathcal{C}_1$. Assuming event $\mathcal{A}_1$ (as defined in Section \ref{basicresultssection} of the Appendix), the same reasoning as in the proof of Theorem 8 of \cite{bhattacharya2023inferences} shows that $\wh\lambda_k$ is the only eigenvalue of $\hat{\bL}$ that is inside $\mathcal{C}_1$. 
As in Theorem 8 of \cite{bhattacharya2023inferences}, we define
\begin{align*}
\bxi_k\bxi_k^{T} & := \bP_k \\
\wh\bxi_k\wh\bxi_k^\top & := \wh\bP_k \\
    \Delta\bP_k & = \frac{1}{2\bpi i}\oint_{\mathcal{C}_1}\left(\lambda\bI-\bL_0 \right)^{-1}\bE\left(\lambda\bI-\bL_0 \right)^{-1}d\lambda
\end{align*}
As in \cite{bhattacharya2023inferences}, we also define and analyze the second-order deviation $\delta_{i, k}$ of the eigenvectors. However, our definition of $\delta_{i, k}$ differs from that of \cite{bhattacharya2023inferences}, as the second-order deviation no longer comes solely from the term $(\wh\bP_k-\bP-\Delta\bP_k)\bxi_k$, but also from second-order terms within $\Delta\bP_k$. In particular, we denote the first-order terms within $\Delta\bP_k$ as
$$
\Delta\bP_k' = \sum_{1 \leq i \leq n, i \neq k}\frac{\bxi_i^{T}\bE_1 \bxi_k}{\lambda_k-\lambda_i}\bxi_i
$$
Then
\begin{align*}
      \delta_{i, k} & :=  \wh\bxi_k-\bxi_k-\Delta\bP_k'\bxi_k \nonumber \\
     & = \left[\left(\wh\bP_k-\bP-\Delta\bP_k\right)\bxi_k\right](i) + \left[(\Delta\bP_k - \Delta\bP_k')\bxi_k\right](i) + \left[(1-\wh\bxi_k^\top\bxi_k)\wh\bxi_k\right](i) \numberthis \label{theorem1eq1}
\end{align*}
In the first part of this proof, we will seek to derive upper bounds for $\delta_{i, k}$. We begin with bounding the first term on the RHS of \eqref{theorem1eq1}. By similar reasoning as to Theorem 8 of \cite{bhattacharya2023inferences},
\small \begin{align*}
    (\wh\bP_k-\bP-\Delta\bP_k)\bxi_k & = \left[\frac{1}{2\bpi i}\oint_{\mathcal{C}_1}\bDelta_1 d\lambda\right] \bxi_k + \left[\frac{1}{2\bpi i}\oint_{\mathcal{C}_1}\bDelta_2 d\lambda\right] \bxi_k \nonumber \\
    & = \bN_{1k}\bE\bN_{1k}\bE\bxi_k-\bN_{2k}\bE\bxi_k\bxi_k^{T}\bE\bxi_k-\bxi_k\bxi_k^{T}\bE\bN_{2k}\bE\bxi_k + \left[\frac{1}{2\bpi i}\oint_{\mathcal{C}_1}\bDelta_2 d\lambda\right] \bxi_k \numberthis \label{ointDelta1}
\end{align*}
\normalsize
where 
\begin{align*}
    \bDelta_1 & := \left(\lambda\bI-\bL_0 \right)^{-1} \bE\left(\lambda\bI-\bL_0 \right)^{-1}\bE\left(\lambda\bI-\bL_0 \right)^{-1}\\
    \bDelta_2 & := \left(\lambda\bI-\hat{\bL} \right)^{-1}\bE\left(\lambda\bI-\bL_0 \right)^{-1}\bE\left(\lambda\bI-\bL_0 \right)^{-1}\bE\left(\lambda\bI-\bL_0 \right)^{-1} \\
    \bN_{1k} &:= \sum_{1 \leq i \leq n, i \neq k} \frac{1}{\lambda_k-\lambda_i} \bxi_i\bxi_i^{T},\quad \bN_{2k} := \sum_{1 \leq i \leq n, i \neq k} \frac{1}{(\lambda_k-\lambda_i)^2}\bxi_i\bxi_i^{T} \numberthis \label{N1N2}
\end{align*}
Under the event $\mathcal{A}_1$,
\begin{align*}
    \left\|\left[\frac{1}{2\bpi i}\oint_{\mathcal{C}_1}\bDelta_2 d\lambda\right] \bxi_k\right\|_{\infty} &\leq \left\|\left[\frac{1}{2\bpi i}\oint_{\mathcal{C}_1}\bDelta_2 d\lambda\right] \bxi_k\right\|_{2} \nonumber \\
    &\leq \frac{1}{2\bpi}\oint_{\mathcal{C}_1}\left\|\bDelta_2\right\| d\lambda \lesssim\frac{\norm{E}^3 r}{\lambda_k^{3}} \numberthis \label{ointDelta2}
\end{align*} so it remains to bound $\left[\frac{1}{2\bpi i}\oint_{\mathcal{C}_1}\bDelta_1 d\lambda\right] \bxi_k$. To this end, we analyze each of its constituent terms separately.

\begin{itemize}[topsep=0pt,itemsep=1ex,partopsep=0ex,parsep=1ex] 
    \item[(i)] Control $\left\|\bN_{1k}\bE\bN_{1k}\bE\bxi_k\right\|_\infty$: For any $i \leq n$, we expand and write
    \begin{align*}
        [\bE\bN_{1k}\bE\bxi_k](i) &=  \left[(\boldsymbol{\trueregdegreehalf} \bW\boldsymbol{\trueregdegreehalf})\bN_{1k}\bE\bxi_k\right](i) + [(\bE_1 - \boldsymbol{\trueregdegreehalf}\bW\boldsymbol{\trueregdegreehalf})\bN_{1k}\bE\bxi_k](i) ] \\
        & + [(\bE - \bE_1)\bN_{1k}\bE\bxi_k](i) \numberthis \label{improveddevboundseq1}
    \end{align*} 
    To bound the first term on the RHS of \eqref{improveddevboundseq1}, we define a leave-one-out version of the first order error term. Specifically, we 
    \begin{enumerate}[topsep=0pt,itemsep=0ex,partopsep=1ex,parsep=1ex]
        \item Consider the leave-one-out noise matrix $\bW^{(i)}$ obtained by replacing all the elements in $i$-th row and $i$-the column of original $\bW$ with $0$ for $i\in [n]$.
        \item Define the leave-one-out version of the $\bDelta_{D, 1}$ matrix. Namely, define the diagonal $n \times n$ matrix $\bDelta^{(i)}_{D, 1}$ such that for all $j \leq n$,
        $$
\bDelta^{(i)}_{D, 1} (j, j) = \trueregdegree^{-\frac{3}{2}}(j, j) \left(\sum_{l = 1}^n W^{(i)}_{jl} + \frac{\delta}{n} \sum_{l = 1}^n\sum_{m = 1}^n W^{(i)}_{lm}\right)
$$
\item Then, we define
$$
    \bE_1^{(i)} := \boldsymbol{\trueregdegreehalf} \bW^{(i)} \boldsymbol{\trueregdegreehalf} + \bDelta^{(i)}_{D, 1} \bH \boldsymbol{\trueregdegreehalf} + \boldsymbol{\trueregdegreehalf} \bH \bDelta^{(i)}_{D, 1}
    $$
    \end{enumerate} As a result, $\bW^{(i)}$ is independent of $\bE_1^{(i)}$. To bound the first term on the RHS of \eqref{improveddevboundseq1}, we decompose further to apply triangle inequality:
    \begin{align*} 
\left[(\boldsymbol{\trueregdegreehalf} \bW\boldsymbol{\trueregdegreehalf})\bN_{1k}\bE\bxi_k\right](i) &\leq \left|[\boldsymbol{\trueregdegreehalf} \bW\boldsymbol{\trueregdegreehalf}\bN_{1k}\bE_1^{(i)}\bxi_k](i)\right| \\
    & + \left|[\boldsymbol{\trueregdegreehalf} \bW\boldsymbol{\trueregdegreehalf}\bN_{1k}(\bE_1 - \bE_1^{(i)})\bxi_k](i)\right| \\
    & + \left|[\boldsymbol{\trueregdegreehalf} \bW\boldsymbol{\trueregdegreehalf}\bN_{1k}(\bE - \bE_1)\bxi_k](i)\right| \numberthis \label{improveddevboundseq2} 
    \end{align*}
    We successively bound each of the terms on the RHS of \eqref{improveddevboundseq2}. By applying Lemma \ref{matrixbernstein1}, we obtain 
    \begin{align*}
        |[(\boldsymbol{\trueregdegreehalf} \bW\boldsymbol{\trueregdegreehalf})\bN_{1k}\bE_1^{(i)}\bxi_k](i)| & \lesssim \frac{\sqrt{\log n}}{n\otheta}  \norm{\bN_{1k}\bE_1^{(i)}\bxi_k}_{F} \\
       &  + \frac{\log n}{\sqrt{n\otheta^2}} \max_{l \leq n} \left\{ \frac{1}{\max\{n\otheta\theta_l, n\otheta^2\}}|[\bN_{1k}\bE_1^{(i)}\bxi_k](l)| \right\} \numberthis \label{u1leaveoneouteq2}
    \end{align*}
    The Frobenius norm term on the RHS of \eqref{u1leaveoneouteq2} can be bounded by a simple norm argument,
    \begin{align*}
        \norm{\bN_{1k}\bE_1^{(i)}\bxi_k}_{F} & \lesssim \norm{\bN_{1k}}\norm{\bE_1^{(i)}}\norm{\bxi_k}_{F} \lesssim \frac{\log n}{\sqrt{n\otheta^2} \lambda_k} \numberthis \label{theorem1eq3}
    \end{align*}
    whereas the second term on the RHS can be bounded by Lemmas \ref{matrixbernstein1}, \ref{Wconcentration2}, \ref{lemmaN1N2},
    \begin{align*}
        [\bN_{1k}\bE_1^{(i)}\bxi_k](i) & \lesssim \frac{1}{\lambda_k} |[\bE_1^{(i)}\bxi_k](i)| + \sqrt{\frac{(K - 1) \mu_i}{n\lambda_k^2}} \norm{\bE_1^{(i)}\bxi_k}_{2} \\
    & \lesssim  \frac{\log n}{\lambda_k} (\sqrt{\frac{\theta_i}{n^2\otheta^3}} + \sqrt{\frac{\mu_i}{n^3\otheta^4}}) + \frac{1}{\lambda_k} [\sqrt{\frac{(K - 1)\mu_i \log n}{n^2\otheta^2}}] \\
    & \lesssim \frac{\log n}{n\lambda_k\otheta} (\sqrt{\frac{\theta_i}{\otheta}} + 1) + \frac{1}{\lambda_k} [\sqrt{\frac{(K - 1)\mu_i \log n}{n^2\otheta^2}}] :=\rho_{1_{i, k}} \numberthis \label{u1leaveoneouteq3}
    \end{align*}

    with probability at least $1-O(n^{-11})$. Plugging \eqref{theorem1eq3} and \eqref{u1leaveoneouteq3} into the above equation tells us 
    \begin{align*}
        |[(\boldsymbol{\trueregdegreehalf} \bW\boldsymbol{\trueregdegreehalf})\bN_{1k}\bE_1^{(i)}\bxi_k](i)| &\lesssim \frac{\log n}{\sqrt{n^3\otheta^4} \lambda_k} + \frac{\log n}{\sqrt{n\otheta^2}} \max_{l \leq n} \left\{ \frac{\rho_{1_{l, k}} }{\max\{n\otheta\theta_l, n\otheta^2\}}\right\} \\
        &  := \chi_{\scriptscriptstyle 1_{l, k}} \numberthis \label{u1leaveoneouteq4}
    \end{align*}
    with probability at least $1-O(n^{-11})$.

    To bound the remaining terms, we use a simple spectral norm argument alongside the bounds in Lemmas \ref{e2concentrationlemma2}, \ref{spectralnormE}, \ref{Wconcentration2} and Corollary \ref{corN1N2}:
    \begin{align*}
        |[\boldsymbol{\trueregdegreehalf} \bW\boldsymbol{\trueregdegreehalf}\bN_{1k}(\bE_1 - \bE_1^{(i)}) \bxi_k](i)| & \leq \norm{[\boldsymbol{\trueregdegreehalf} \bW\boldsymbol{\trueregdegreehalf}]_{i, \cdot}} \norm{\bN_{1k}} \norm{\bE_1^{(i)} \bxi_k} \\
        & \lesssim [\frac{\sqrt{\log n}}{n\otheta^2}] \cdot [\lambda_k^{-1}] \cdot [\sqrt{\frac{1}{n^2\otheta^2}} \log n + \frac{1}{n\otheta^2} \sqrt{\frac{1}{n}} \log n] \numberthis \label{u1leaveoneouteq5} \\
        \norm{\boldsymbol{\trueregdegreehalf} \bW\boldsymbol{\trueregdegreehalf}\bN_{1k}(\bE - \bE_1) \bxi_k} & \leq \norm{[\boldsymbol{\trueregdegreehalf} \bW\boldsymbol{\trueregdegreehalf}]_{i, \cdot}} \norm{\bN_{1k}} \norm{\bE - \bE_1} \norm{\bxi_k} \\
         & \lesssim [\frac{\sqrt{\log n}}{n\otheta^2}] \cdot [\lambda_k^{-1}] \cdot \frac{\log n }{n\otheta^2} \\
         & \lesssim \frac{\log^{\frac{3}{2}} n}{n^{2}\otheta^4\lambda_k} \numberthis \label{u1leaveoneouteq6} \\
         |[(\bE_1 - \boldsymbol{\trueregdegreehalf} \bW\boldsymbol{\trueregdegreehalf})\bN_{1k} \bE \bxi_k](i)| & \leq \norm{[\bE_1 - \boldsymbol{\trueregdegreehalf} \bW\boldsymbol{\trueregdegreehalf}]_{i, \cdot}}\norm{\bN_{1k}} \norm{\bE} \norm{\bxi_k} \\
        & \lesssim (\frac{\log n}n \sqrt{\frac{\theta_i}{\otheta^3}} + \frac{\log n}{n\otheta}) \cdot \sqrt{\frac{\log n}{n \otheta^2 \lambda_k^2 }} \numberthis \label{u1leaveoneouteq7}
    \end{align*}
    where the last inequality follows from the triangle inequality bound $$\norm{[\bE_1 - \boldsymbol{\trueregdegreehalf} \bW\boldsymbol{\trueregdegreehalf}]_{i, \cdot}} \leq \norm{[\boldsymbol{\trueregdegreehalf}\bH\bDelta_{D, 1}]_{i, \cdot}} + \norm{[\bDelta_{D, 1}\bH\boldsymbol{\trueregdegreehalf}]_{i, \cdot}}$$ and plugging in the bounds of the corresponding lemmas for the RHS terms.
    
    Combining the bounds in \eqref{u1leaveoneouteq4}, \eqref{u1leaveoneouteq5}, \eqref{u1leaveoneouteq6}, and \eqref{u1leaveoneouteq7}, we obtain
    \begin{align}
        |[\bE_1 \bN_{1k} \bE \bxi_k](i)| & \lesssim \chi_{\scriptscriptstyle 1_{i, k}} + \frac{\log^{\frac{3}{2}} n}{n^{2}\otheta^4\lambda_k} + \sqrt{\frac{\max\{\theta_i, \otheta\}\log^3 n}{n^3 \otheta^5 \lambda_k^2}} 
        \label{u1leaveoneouteq9}
    \end{align}

    Lastly, we can bound the remaining, rightmost RHS term in \eqref{improveddevboundseq1} via a simple spectral norm argument, 
    \begin{align*}
        [(\bE - \bE_1)\bN_{1k}\bE\bxi_k](i) & = \norm{[\bE - \bE_1]_{i, \cdot}\bN_{1k}\bE\bxi_k} \\
        & \leq \norm{[\bE - \bE_1]_{i, \cdot}}\norm{\bN_{1k}}\norm{\bE}\norm{\bxi_k} \\
        & \lesssim \frac{\log n}{n\otheta^2} \cdot \frac{1}{\lambda_k} \cdot \sqrt{\frac{\log n}{n\otheta^2}} \cdot 1 \\
        & = \sqrt{\frac{\log^3 n}{n^3\otheta^6 \lambda_k^2}} \numberthis \label{u1leaveoneouteq10}
    \end{align*}

    Combining our bounds in equations \eqref{u1leaveoneouteq9} and \eqref{u1leaveoneouteq10}, we deduce
\begin{align*}
                [\bE\bN_{1k}\bE\bxi_k](i) & \lesssim \left(\chi_{\scriptscriptstyle 1_{i, k}} + \frac{\log^{\frac{3}{2}} n}{n^{2}\otheta^4\lambda_k} + \sqrt{\frac{\max\{\theta_i, \otheta\}\log^3 n}{n^3 \otheta^5 \lambda_k^2}} \right) + \sqrt{\frac{\log^3 n}{n^3\otheta^6 \lambda_k^2}} \\
                & 
 \lesssim\left(\chi_{\scriptscriptstyle 1_{i, k}} \sqrt{\frac{\max\{\theta_i, \otheta\}\log^3 n}{n^3 \otheta^5 \lambda_k^2}} \right) + \sqrt{\frac{\log^3 n}{n^3\otheta^6 \lambda_k^2}}
            \end{align*} 
            where the last line follows from using the first item in Assumption \ref{regconds3} to simplify.
            
           Now, by Corollary \ref{corN1N2}, we have
    \begin{align*}
\left\|\bN_{1k}\bE\bN_{1k}\bE\bxi_k\right\|_\infty\lesssim & \frac{1}{\lambda_k}\left\|\bE\bN_{1k}\bE\bxi_k\right\|_{\infty} +\sqrt{\frac{(K-1)\max_{i \leq n}\mu_i}{n\lambda_k^{2}}}\left\|\bE\bN_{1k}\bE\bxi_k\right\|_2 \\
        \lesssim &\frac{1}{\lambda_k}\left\|\bE\bN_{1k}\bE\bxi_k\right\|_{\infty}+\sqrt{\frac{(K-1)\max_{i \leq n}\mu_i}{n\lambda_k^{2}}}\left\|\bE\right\|^2\left\|\bN_{1k}\right\|\left\|\bxi_k\right\|_2\\
        \lesssim & \left( \max_{i \leq n}\left\{\frac{\chi_{\scriptscriptstyle 1_{i, k}}}{\lambda_k }\right\} + \sqrt{\frac{\theta_{\text{max}}\log^3 n}{n^3 \otheta^5 \lambda_k^4}} + \sqrt{\frac{\log^3 n}{n^3\otheta^6 \lambda_k^2}} \right) + \frac{\log n\sqrt{K\theta_{\text{max}}}}{n^{\frac{3}{2}}\otheta^{\frac{5}{2}} \lambda_k^2} \numberthis \label{u1leaveoneouteq8}
    \end{align*}
    with probability at least $1-O(n^{-11})$.
             
             \item[(ii)] Control $\left\|\bN_{2k}\bE_1\bxi_k\bxi_k^{T}\bE_1\bxi_k\right\|_\infty$: First, with probability at least $1-O(n^{-11})$, we have $|\bxi_k^{T}\bE\bxi_k|\lesssim \norm{\bxi_k}\norm{\bE}\norm{\bxi_k^T} \lesssim \norm{\bE}$. Using the definition of $\rho_{1_{i, k}}$ from \eqref{u1leaveoneouteq3}, 
    \begin{align*}
        &\left\|\bN_{2k}\bE_1\bxi_k\bxi_k^{T}\bE_1\bxi_k\right\|_\infty = \left|\bxi_k^{T}\bE_1\bxi_k\right|\left\|\bN_{2k}\bE_1\bxi_k\right\|_\infty\lesssim \sqrt{\frac{\log n}{n\otheta^2}}\left\|\bN_{2k}\bE_1\bxi_k\right\|_\infty \\
        &\lesssim\sqrt{\frac{\log n}{n\otheta^2}} \Big(\frac{1}{\lambda_k^{2}}\left\|\bE_1\bxi_k\right\|_\infty+\sqrt{\frac{(K-1)\max_{i \leq n}\mu_i}{n\lambda_k^{4}}}\left\|\bE_1\bxi_k\right\|_2\Big) \lesssim \frac{\sqrt{\frac{\log n}{n\otheta^2}}}{\lambda_k^{}} \max_{i \leq n} \rho_{1_{i, k}}
    \end{align*}
    with probability at least $1-O(n^{-11})$, where the second inequality uses Corollary \ref{corN1N2}, and the last inequality uses Lemma \ref{Wconcentration2}.
    \item[(iii)] Control $\left\|\bxi_k\bxi_k^{T}\bE_1\bN_{2k}\bE_1\bxi_k\right\|_\infty$: Using Lemma \ref{lemmaN1N2}, under the event $\mathcal{A}_1$,
    \begin{align*}
        \left\|\bxi_k\bxi_k^{T}\bE_1\bN_{2k}\bE_1\bxi_k\right\|_\infty &= \left\|\bxi_k\right\|_\infty\left|\bxi_k^{T}\bE_1\bN_{2k}\bE_1\bxi_k\right|\\
        &\leq \left\|\bxi_k\right\|_\infty\left|\bxi_k\right\|_2^2\left\|\bE_1\right\|^2\left\|\bN_{2k}\right\| \\
        & \leq \frac{\sqrt{\max_{i \leq n} \mu_i} \log n}{n^{\frac{3}{2}}\otheta^2\lambda_k^{2}} \lesssim \frac{\sqrt{\theta_{\text{max}}}\log n}{n^{\frac{3}{2}}\otheta^\frac{5}{2}\lambda_k^{2}}
    \end{align*}
    with probability at least $1-O(n^{-11})$, where the final inequality uses the bounds in Lemma \ref{incoherence}.
\end{itemize} 
Combine these three parts with \eqref{ointDelta1}, we get
\begin{align*}
\left\|\left[\frac{1}{2\bpi i}\oint_{\mathcal{C}_1}\bDelta_1 d\lambda\right] \bxi_k\right\|_\infty\lesssim& \left( ( \max_{i \leq n}\left\{\frac{\chi_{\scriptscriptstyle 1_{i, k}}}{\lambda_k }\right\} + \sqrt{\frac{\theta_{\text{max}}\log^3 n}{n^3 \otheta^5 \lambda_k^4}} + \sqrt{\frac{\log^3 n}{n^3\otheta^6 \lambda_k^2}} ) + \frac{\log n\sqrt{K\theta_{\text{max}}}}{n^{\frac{3}{2}}\otheta^{\frac{5}{2}} \lambda_k^2} \right)  \\
    &  + \frac{1}{\lambda_k}\sqrt{\frac{\log n}{n \otheta^2}} \rho_{1_{i, k}} + \frac{\sqrt{\theta_{\text{max}}}\log n}{n^{\frac{3}{2}}\otheta^\frac{5}{2}\lambda_k^{2}}  =: \chi_{\scriptscriptstyle 2_{i, k}}  \numberthis \label{ointDelta1bound}
\end{align*}
with probability at least $1-O(n^{-11})$. Combining \eqref{ointDelta1bound} and \eqref{ointDelta2}, we get
\begin{align}
    |[(\wh\bP_k-\bP_k-\Delta\bP_k)\bxi_k](i)|\lesssim&\frac{\log^{1.5}n}{\lambda_k^{3}n^{\frac{3}2}\otheta^3}+\chi_{\scriptscriptstyle 2_{i, k}} \label{extraeq1}
\end{align}

Now, to bound $(\Delta\bP_k - \Delta\bP_k')\bxi_k$, we apply Corollary \ref{corN1N2} and Lemmas \ref{spectralnormE} and \ref{Wconcentration2}, 
\begin{align*}
    |[(\Delta\bP_k - \Delta\bP_k')\bxi_k](i)| & = \sum_{1 \leq i \leq n, i \neq k}\frac{\bxi_i^{T}(\bE - \bE_1) \bxi_k}{\lambda_k-\lambda_i}\bxi_i = N_{1k}(\bE - \bE_1) \bxi_k \\
    & \lesssim \frac{1}{\lambda_k} |[(\bE - \bE_1) \bxi_k](i)| + \sqrt{\frac{(K - 1) \mu_i}{\lambda_k^2 n}} \norm{(\bE - \bE_1) \bxi_k}  \\
    & \lesssim \frac{1}{\lambda_k} |[(\bE - \bE_1) \bxi_k](i)| + \sqrt{\frac{(K - 1) \mu_i}{\lambda_k^2 n}} \norm{(\bE - \bE_1)}\norm{\bxi_k}_2 \\
    & \lesssim \frac{1}{\lambda_k} \left(\sqrt{\frac{\log^2(n)}{n^3\otheta^4}}(1 + \sqrt{\frac{1}{n\otheta\theta_i}}) + \sqrt{\frac{\log^3n}{n^3\otheta^6}} \right) + \sqrt{\frac{(K - 1) \mu_i \log^2 n}{\lambda_k^2 n^3\otheta^4}} \\
    & \lesssim \sqrt{\frac{\log^2(n)}{n^3\otheta^4\lambda_k^2}}(1 + \sqrt{K\mu_i}) + \sqrt{\frac{\log^3n}{n^3\otheta^6\lambda_k^2}}\numberthis \label{extraeq2}
\end{align*}
where in the last line, we use the fact that $\theta_i \geq C\frac{\log n}{n\otheta}$ to simplify.

It remains to bound $\|[(1-\wh\bxi_k^\top\bxi_k)\wh\bxi_k](i)\|$. To this end,
\begin{align*}
    \|(1-\wh\bxi_k^\top\bxi_k)\wh\bxi_k\|_\infty &= |1-\wh\bxi_k^\top\bxi_k|\|\wh\bxi_k\|_\infty\lesssim (\frac{\norm{E}}{\lambda_k})^2\left(\left\|\bxi_k\right\|_\infty+\left\|\hat{\bxi}_k-\bxi_k\right\|_\infty\right) \nonumber \\
    &\leq \frac{\log n}{n\otheta^2\lambda_k^{2}}\left(\sqrt{\frac{\max_{i \leq n} \mu_i}{n}}+\left\|\hat{\bxi}_k-\bxi_k\right\|_2\right)\lesssim \frac{\log n \sqrt{\theta_{\text{max}}}}{n^{\frac{3}{2}}\otheta^{\frac{5}{2}}\lambda_k^{2}}+ \frac{\log^{\frac{3}{2}}n}{n^{\frac{3}{2}}\otheta^3} \numberthis \label{extraeq3}
\end{align*}
where the first and last inequalities follow by the Davis-Kahan sin theta theorem. Plugging \eqref{extraeq1}, \eqref{extraeq2}, and \eqref{extraeq3} in \eqref{theorem1eq1}, and using the assumption that $\sqrt{\frac{\log n}{n\otheta\theta_i}} \ll 1$ to simplify, we obtain
\begin{align*}
    \left\|\delta\right\|_\infty & \leq \max_{i \leq n} \left(|[(\wh\bP_k-\bP_k-\Delta\bP_k)\bxi_k](i)|+ |[(\Delta\bP_k - \Delta\bP_k')\bxi_k](i)|\right) + \|(1-\wh\bxi_k^\top\bxi_k)\wh\bxi_k\|_\infty \\
    & \lesssim \max_{i \leq n} \left\{\chi_{\scriptscriptstyle 2_{i, k}}\right\} + \frac{\log^{1.5}n}{\lambda_k^{3}n^{\frac{3}2}\otheta^3} + (\frac{\log n \sqrt{\theta_{\text{max}}}}{n^{\frac{3}{2}}\otheta^{\frac{5}{2}}\lambda_k^{2}}+ \frac{\log^{\frac{3}{2}}n}{n^{\frac{3}{2}}\otheta^3})\numberthis \label{theorem9eq10}
\end{align*}
with probability at least $1-O(n^{-10})$.

In the second part of the proof, we turn to bounding the first-order error $\Delta\bP_k'(i) = [\bN_{1k}E_{1} \bxi_k](i)$. To this end, we decompose \begin{equation}
    \Delta\bP_k' = \bN_{1k}\boldsymbol{\trueregdegreehalf} \bW \boldsymbol{\trueregdegreehalf} \bxi_k + \bN_{1k}\bDelta_{D, 1} \bH \boldsymbol{\trueregdegreehalf} \bxi_k + \boldsymbol{\trueregdegreehalf} \bH \bDelta_{D, 1}\bxi_k \numberthis \label{theorem9begbound} 
\end{equation} and bound each of the terms separately. Since each of the terms are a weighted sum $\sum_{1 \leq a \leq b \leq n} c_{ab} W_{ab}$ of entries of $\bW$, our general strategy will be to obtain bounds on the coefficients $|c_{ab}|$ and variance of the sum, respectively, in order to apply Bernstein's inequality.
\begin{enumerate}[topsep=0pt,itemsep=0ex,partopsep=1ex,parsep=1ex]
    \item[(i)] \underline{Bound $[\bN_{1k}\boldsymbol{\trueregdegreehalf} \bW \boldsymbol{\trueregdegreehalf} \bxi_k](i)$:} By algebraic manipulations, 
    \begin{align*}
        [\bN_{1k}\boldsymbol{\trueregdegreehalf} \bW \boldsymbol{\trueregdegreehalf} \bxi_k](i) & = [\bN_{1k}]_{i, \cdot} [\boldsymbol{\trueregdegreehalf} \bW \boldsymbol{\trueregdegreehalf}] \bxi_k \\
        & = \sum_{a = 1}^n \sum_{b = 1}^n [N_{1k}]_{i, a} \xi_k(b) \trueregdegreehalf(a, a)\trueregdegreehalf(b, b) W_{ab}
    \end{align*}
    implying 
    $$
c_{ab} = \begin{cases}
[N_{1k}]_{i, a}\xi_k(b)\trueregdegreehalf(a, a)\trueregdegreehalf(b, b) + [N_{1k}]_{i, b}\xi_k(a)\trueregdegreehalf(b, b)\trueregdegreehalf(a, a) & \hbox{ if $a \neq b$} \\
[N_{1k}]_{i, a}\xi_k(a)\trueregdegree^{-1}(a, a) & \hbox{ if $a = b$}
\end{cases}
$$
In the above equation, we can bound $[N_{1k}]_{i, a}$ via Corollary \ref{corN1N2} and $|\trueregdegreehalf(a, a)| \lesssim \max(\{n\otheta\theta_a, n\otheta^2\})^{-\frac{1}{2}} $ via Lemma \ref{basicexpansiondegreeresults}. Plugging in these bounds, we derive upper bounds for the coefficients $c_{ab}$, 
\begin{align*}
    |c_{ab}| & \lesssim \trueregdegreehalf(a, a)\trueregdegreehalf(b, b) \cdot \begin{cases}
    \sqrt{\frac{K^2 \mu_i \mu_a \mu_b}{n^3\lambda_k^2}} & \hbox{if $a, b \neq i$} \\ 
    \sqrt{\frac{\mu_b}{n\lambda_k^2}} + \sqrt{\frac{K^2 \mu_i \mu_a \mu_b}{n^3\lambda_k^2}} & \hbox{if $a = i, b \neq i$} \\
     \sqrt{\frac{\mu_a}{n\lambda_k^2}} + \sqrt{\frac{K^2 \mu_i \mu_a \mu_b}{n^3\lambda_k^2}} & \hbox{if $a \neq i, b = i$} \\
     \sqrt{\frac{\mu_a}{n\lambda_k^2}} + \sqrt{\frac{K^2 \mu_i \mu_a \mu_b}{n^3\lambda_k^2}} & \hbox{if $a, b = i$}
\end{cases} \numberthis \label{corN1eq4}
\end{align*}
Using \eqref{corN1eq4}, we can bound the variance by 
\begin{align*}
    \hbox{var} & \lesssim \sum_{1 \leq a \leq b \leq n} c_{ab}^2 \theta_a\theta_b \\
    & = \sum_{a, b \neq i} c_{ab}^2 \theta_a\theta_b + \sum_{a = i, b \neq i} c_{ab}^2 \theta_a\theta_b + \sum_{a \neq i, b = i} c_{ab}^2 \theta_a\theta_b + \sum_{a, b = i} c_{ab}^2 \theta_a\theta_b
\end{align*}
Below, we bound each of the terms on the RHS. First of,
\begin{align*}
    \sum_{a, b \neq i} c_{ab}^2 \theta_a\theta_b & \leq \frac{K^2 \mu_i}{n^3 \lambda_k^2} \sum_{a, b \neq i} [\mu_a \theta_a \trueregdegree^{-1}(a, a)][\mu_b \theta_b \trueregdegree^{-1}(b, b)] \\
    & \lesssim \frac{K^2 \mu_i}{n^3 \lambda_k^2\otheta^2} \\
    \sum_{a = i, b > i} c_{ab}^2 \theta_a\theta_b & \leq \sum_{a = i, b > i} \left( \trueregdegreehalf(a, a)\trueregdegreehalf(b, b)  \sqrt{\frac{\mu_b}{n\lambda_k^2}}\right)^2 \theta_a\theta_b  \\
   & + \sum_{a = i, b > i} \left( \trueregdegreehalf(a, a)\trueregdegreehalf(b, b) \sqrt{\frac{K^2 \mu_i \mu_a \mu_b}{n^3\lambda_k^2}}\right)^2 \theta_a\theta_b\\
    & \lesssim \frac{\theta_i \trueregdegree^{-1}(i, i)}{n\lambda_k^2} \sum_{b \neq i} \mu_b \theta_b \trueregdegree^{-1}(b, b) + \frac{K^2 \theta_i \trueregdegree^{-1}(i, i) \mu_i^2}{n^3 \lambda_k^2}\sum_{b \neq i} \mu_b \theta_b \trueregdegree^{-1}(b, b)\\
    & \lesssim \frac{\theta_i}{n^2\otheta^2(\minthetai)\lambda_k^2} + \frac{K^2 \mu_i^2}{n^4 \lambda_k^2\otheta^2} \\
\end{align*}
where we use Lemma \ref{basicsumresults} to control the summands in the second-to-last line of each bound. Furthermore, by applying the same reasoning as in the bound for $\sum_{a = i, b > i} c_{ab}^2 \theta_a\theta_b$, we can also bound
$$
 \sum_{a < i, b = i} c_{ab}^2 \theta_a\theta_b \lesssim \frac{\theta_i}{n^2\otheta^2(\minthetai)\lambda_k^2} + \frac{K^2 \mu_i^2}{n^4 \lambda_k^2\otheta^2} 
$$ Lastly,
\begin{align*}
    \sum_{a, b = i} c_{ab}^2 \theta_a\theta_b & \leq  2 \theta_i^2 \left[ \sum_{a, b = i} \left(\sqrt{\frac{\mu_a}{n\lambda_k^2}}\trueregdegree^{-1}(i, i)\right)^2 + \left(\sqrt{\frac{K^2 \mu_i \mu_a \mu_b}{n^3\lambda_k^2}}\trueregdegree^{-1}(i, i)\right)^2 \right] \\
    & \lesssim \frac{\mu_i}{n^3\otheta^2\lambda_k^2} + \frac{K^2 \mu_i^3}{n^5 \lambda_k^2\otheta^2}
\end{align*}
Adding the bounds together, 
$$
\hbox{var} \lesssim \frac{\theta_i}{n^2\otheta^2(\minthetai)\lambda_k^2} + \frac{K^2 \mu_i}{n^3 \lambda_k^2\otheta^2} + \frac{K^2 \mu_i^2}{n^4 \lambda_k^2\otheta^2} + \frac{K^2 \mu_i^3}{n^5 \lambda_k^2\otheta^2}
$$
To obtain a uniform bound for all $|c_{ab}|$, we can plug in the bound for $\sqrt{\mu_a}$ from Lemma \ref{incoherence} into \eqref{corN1eq4} to obtain
$$
\max_{1 \leq a, b \leq n} |c_{ab}| \lesssim \frac{1}{n^{\frac{3}{2}}\otheta^2 d_k} + \frac{K\sqrt{\mu_i}}{n^{\frac{5}{2}}\otheta^2 d_k} 
$$
Applying Bernstein shows
\begin{align*}
    [\bN_{1k}\boldsymbol{\trueregdegreehalf} \bW \boldsymbol{\trueregdegreehalf} \bxi_k](i) & \lesssim (\frac{1}{n\otheta\lambda_k}\sqrt{\frac{\theta_i}{\minthetai}} + \frac{K\mu_i^{\frac{1}{2}}}{n^{\frac{3}{2}} \lambda_k\otheta} + \frac{K \mu_i}{n^2 \lambda_k\otheta} + \frac{K \mu_i^{\frac{3}{2}}}{n^{\frac{5}{2}} \lambda_k\otheta}) \sqrt{\log n} \\
    & + (\frac{1}{n^{\frac{3}{2}}\otheta^2 d_k} + \frac{K\sqrt{\mu_i}}{n^{\frac{5}{2}}\otheta^2 d_k}) \log n
\end{align*}

\item[(ii)]  \underline{Control $[\bN_{1k}\bDelta_{D, 1} \bH \boldsymbol{\trueregdegreehalf} \bxi_k](i)$:} By algebraic manipulations, 
\small \begin{align*}
[\bN_{1k}\bDelta_{D, 1} \bH \boldsymbol{\trueregdegreehalf} \bxi_k](i) & = [\bN_{1k}]_{i, \cdot} [\bDelta_{D, 1} \bH \boldsymbol{\trueregdegreehalf}] \bxi_k \\
& = \sum_{a = 1}^n \sum_{x = 1}^n [N_{1k}]_{i, a} \xi_k(x) \Delta_{D, 1}(a, a)  \trueregdegreehalf(x, x) H_{ax} \\
& = -\frac{1}{2} \sum_{a = 1}^n \sum_{x = 1}^n [N_{1k}]_{i, a} \xi_k(x) \left[ \frac{\sum_{y = 1}^n W_{ay} + \delta n^{-1} \sum_{x = 1}^n \sum_{y = 1}^n W_{xy}}{\trueregdegree^{\frac{3}{2}}(a, a)}\right] \trueregdegreehalf(x, x) H_{ax}   \\
& = -\frac{1}{2} \sum_{a = 1}^n \sum_{b = 1}^n W_{ab} \cdot \biggl[[N_{1k}]_{i, a} \trueregdegree^{-\frac{3}{2}}(a, a) \cdot \Bigl[\sum_{x = 1}^n \xi_k(x) \trueregdegreehalf(x, x) H_{ax} \Bigr] \\
& + \delta n^{-1} \sum_{a = 1}^n \sum_{x = 1}^n [N_{1k}]_{i, a} \xi_k(x) \trueregdegree^{-\frac{3}{2}}(a, a) \trueregdegreehalf(x, x) H_{ax}\biggr] \\
& := -\frac{1}{2} \sum_{a = 1}^n \sum_{b = 1}^n W_{ab} \cdot \biggl[d_{ab} \biggr] := -\frac{1}{2} \sum_{a = 1}^n \sum_{b = 1}^n W_{ab} \cdot \biggl[A_{1, ab} + A_{2, ab} \biggr] 
    \end{align*}

    \normalsize
    
     To bound $d_{ab}$, we bound each of its two constituent terms $A_{1, ab}, A_{2, ab}$ separately. First, we bound $A_{1, ab}$.  Plugging in the upper bounds for $H_{ax}, \trueregdegreehalf(x, x)$ from Lemma \ref{basicexpansiondegreeresults} yields
     \begin{align*}
         \sum_{x = 1}^n \xi_k(x) \trueregdegreehalf(x, x) H_{ax} & = \frac{\theta_a}{n\otheta^{\frac{1}{2}}} \sum_{x = 1}^n \sqrt{\frac{\mu_x \theta_x^2}{(\minthetax)}} \\
         & \lesssim \frac{\theta_a}{n\otheta^{\frac{1}{2}}} \cdot n \otheta^{\frac{1}{2}} = \theta_a \numberthis \label{theorem9extreq4}
     \end{align*}
     where the last line follows from Lemma \ref{basicsumresults}. As a result, $$|A_{1, ab}| \lesssim [N_{1k}]_{i, a} \trueregdegree^{-\frac{3}{2}}(a, a) \theta_a$$ Plugging in the bounds for $[N_{1k}]_{i, a}$ in Corollary \ref{corN1N2}, we obtain
     \begin{align*}
         [N_{1k}]_{i, a} \trueregdegree^{-\frac{3}{2}}(a, a) \theta_a \lesssim \begin{cases}
         \frac{K\sqrt{\mu_i}}{n^{\frac{5}{2}}\lambda_k \otheta^2} & \hbox{ if $i \neq a$} \\
    \frac{1}{n^{\frac{3}{2}}\lambda_k \otheta^{\frac{3}{2}}\sqrt{\minthetai}} + \frac{K\sqrt{\mu_i}}{n^{\frac{5}{2}}\lambda_k \otheta^2} & \hbox{ otherwise}
         \end{cases} \numberthis \label{theorem9extreq5}
     \end{align*}

     To bound $A_{2, ab}$,
     \begin{align*}
        |A_{2, ab}| & = \delta n^{-1} \sum_{a = 1}^n [N_{1k}]_{i, a} \trueregdegree^{-\frac{3}{2}}(a, a) \sum_{x = 1}^n \xi_k(x)  \trueregdegreehalf(x, x) H_{ax} \\
       & \lesssim \delta n^{-1} \sum_{a = 1}^n [N_{1k}]_{i, a} \trueregdegree^{-\frac{3}{2}}(a, a) \cdot (\theta_a) \\
       & =  \delta n^{-1} \cdot \left(\frac{1}{n^{\frac{3}{2}}\lambda_k \otheta^{\frac{3}{2}}\sqrt{\minthetai}} + \frac{K\sqrt{\mu_i}}{n^{\frac{5}{2}}\lambda_k \otheta^2} \right)
     \end{align*}
     The second line follows from \eqref{theorem9extreq4}, and the third line from summing all the bounds in \eqref{theorem9extreq5}. Combining our bounds for $|B_{1, ab}|$ and $|A_{2, ab}|$, we obtain a bound for $|d_{ab}|$ and thus $|c_{ab}|$:
     $$
|c_{ab}| \lesssim \begin{cases}
         \frac{1}{n^{\frac{5}{2}}\lambda_k \otheta^{\frac{3}{2}}\sqrt{\minthetai}} + \frac{K\sqrt{\mu_i}}{n^{\frac{5}{2}}\lambda_k \otheta^2} & \hbox{ if $i \neq a, b$} \\
     \frac{1}{n^{\frac{3}{2}}\lambda_k \otheta^{\frac{3}{2}}\sqrt{\minthetai}} + \frac{K\sqrt{\mu_i}}{n^{\frac{5}{2}}\lambda_k \otheta^2} & \hbox{ otherwise}
         \end{cases}
$$ Denote $L_2 = \frac{1}{n^{\frac{3}{2}}\lambda_k \otheta^{2}} + \frac{K\sqrt{\mu_i}}{n^{\frac{5}{2}}\lambda_k \otheta^2}$. Then $L_2$ is a global bound for $|c_{ab}|$. 

Now, we can bound the variance by 
\begin{align*}
    \hbox{var} & = \sum_{a, b \neq i} c_{ab}^2 \theta_a\theta_b + \sum_{a = i, b \neq i} c_{ab}^2 \theta_a\theta_b + \sum_{a \neq i, b = i} c_{ab}^2 \theta_a\theta_b + \sum_{a, b = i} c_{ab}^2 \theta_a\theta_b
\end{align*}
where each of the terms on the RHS can be bounded as
\begin{align*}
     \sum_{a, b \neq i} c_{ab}^2 \theta_a\theta_b & \lesssim \sum_{a, b \neq i} \left(\frac{1}{n^{\frac{5}{2}}\lambda_k \otheta^{\frac{3}{2}}\sqrt{\minthetai}}\right)^2 \theta_a\theta_b + \sum_{a, b \neq i} \left(\frac{K\sqrt{\mu_i}}{n^{\frac{5}{2}}\lambda_k \otheta^2}\right)^2 \theta_a\theta_b \\
     & \lesssim \frac{1}{n^3\lambda_k^2\otheta^2} + \frac{K^2 \mu_i}{n^3\lambda_k^2\otheta^2} \\
     \sum_{a = i, b \neq i} c_{ab}^2 \theta_a\theta_b & \lesssim \sum_{a = i, b \neq i} \left(\frac{1}{n^{\frac{3}{2}}\lambda_k \otheta^{\frac{3}{2}}\sqrt{\minthetai}}\right)^2 \theta_a\theta_b + \sum_{a = i, b \neq i} \left(\frac{K\sqrt{\mu_i}}{n^{\frac{5}{2}}\lambda_k \otheta^2}\right)^2 \theta_a\theta_b \\
     & \lesssim \frac{\theta_i}{n^2\lambda_k^2\otheta^2(\minthetai)} + \frac{K^2 \theta_i \mu_i}{n^4\lambda_k^2\otheta^3} \\
     \sum_{a \neq i, b = i} c_{ab}^2 \theta_a\theta_b & \lesssim \frac{\theta_i}{n^2\lambda_k^2\otheta^2(\minthetai)} + \frac{K^2 \theta_i \mu_i}{n^4\lambda_k^2\otheta^3} \\
     \sum_{a, b = i} c_{ab}^2 \theta_a\theta_b & \lesssim \frac{\theta_i}{n^3\lambda_k^2\otheta^3} + \frac{K^2 \mu_i \theta_i^2}{n^5 \lambda_k^2 \otheta^4}
\end{align*}
Therefore, $$\hbox{var} \lesssim \frac{\theta_i}{n^2\lambda_k^2\otheta^2(\minthetai)} + \frac{\minthetai}{n^3\lambda_k^2\otheta^3} + \frac{K^2 \theta_i \mu_i}{n^4\lambda_k^2\otheta^3} + \frac{K^2 \mu_i \theta_i^2}{n^5 \lambda_k^2 \otheta^4}$$ Applying Bernstein shows
\small \begin{align*}
     [\bN_{1k}\bDelta_{D, 1} \bH \boldsymbol{\trueregdegreehalf} \bxi_k](i) & \lesssim \left(\frac{1}{n\lambda_k\otheta}\sqrt{\frac{\theta_i}{\minthetai}} + \frac{1}{n^{\frac{3}{2}}\lambda_k\otheta}\sqrt{\frac{\minthetai}{\otheta}} + \frac{K \sqrt{\mu_i}}{n^2\lambda_k\otheta}\sqrt{\frac{\theta_i}{\otheta}} + \frac{K \theta_i \sqrt{\mu_i}}{n^{\frac{5}{2}} \lambda_k \otheta^2}\right) \sqrt{\log n} \\
     & + \left(\frac{1}{n^{\frac{3}{2}}\lambda_k \otheta^{2}} + \frac{K\sqrt{\mu_i}}{n^{\frac{5}{2}}\lambda_k \otheta^2} \right) \log n
\end{align*}
\normalsize

\item[(iii)]  \underline{Control $[\bN_{1k}\boldsymbol{\trueregdegreehalf} \bH \bDelta_{D, 1} \bxi_k](i)$:} By algebraic manipulations, 
\small \begin{align*}
        [N_{1k}(\boldsymbol{\trueregdegreehalf} \bH \bDelta_{D, 1} \bxi_k)](i) & = [\bN_{1k}]_{i, \cdot} [\boldsymbol{\trueregdegreehalf} \bH \bDelta_{D, 1}] \bxi_k \\
        & = \sum_{x = 1}^n \sum_{a = 1}^n [N_{1k}]_{i, x} \xi_k(a) \trueregdegreehalf(x, x)\Delta_{D, 1}(a, a) H_{xa} \\
        & = \sum_{x = 1}^n \sum_{a = 1}^n [N_{1k}]_{i, x} \xi_k(a) \trueregdegreehalf(x, x) \left[ \frac{\sum_{y = 1}^n W_{ay} + \delta n^{-1} \sum_{x = 1}^n \sum_{y = 1}^n W_{xy}}{\trueregdegree^{\frac{3}{2}}(a, a)}\right] H_{xa} \\
        & = -\frac{1}{2} \sum_{a = 1}^n \sum_{b = 1}^n W_{ab} \cdot \biggl[\trueregdegree^{-\frac{3}{2}}(a, a) \xi_k(a) \cdot \Bigl[\sum_{x = 1}^n [N_{1k}]_{i, x} H_{ax}\trueregdegreehalf(x, x)\Bigr] \\
        & + \delta n^{-1} \sum_{x = 1}^n \sum_{a = 1}^n [N_{1k}]_{i, x} H_{xa}\xi_k(a) \trueregdegreehalf(x, x)\trueregdegree^{-\frac{3}{2}}(a, a)\biggr] \\
        & := -\frac{1}{2} \sum_{a = 1}^n \sum_{b = 1}^n W_{ab} \cdot \biggl[d_{ab} \biggr] := -\frac{1}{2} \sum_{a = 1}^n \sum_{b = 1}^n W_{ab} \cdot \biggl[B_{1, ab} + B_{2, ab} \biggr]
    \end{align*}
    \normalsize To bound $d_{ab}$, we bound each of its two constituent terms $B_{1, ab}, B_{2, ab}$ separately. First, we bound $B_{1, ab}$.  Plugging in the upper bounds for $[N_{1k}]_{i, x}, H_{ax}, \trueregdegreehalf(x, x)$ from  Corollary \ref{corN1N2} and Lemmas \ref{basicresults} and \ref{basicexpansiondegreeresults} yields \small
    \begin{align*}
         \sum_{x = 1}^n [N_{1k}]_{i, x} H_{ax}\trueregdegreehalf(x, x) & \lesssim [N_{1k}]_{i, i} H_{ai}\trueregdegreehalf(i, i) + \sum_{x \neq i}^n [N_{1k}]_{i, x} H_{ax}\trueregdegreehalf(x, x)  \\
         & \lesssim  (\frac{1}{\lambda_k} + \sqrt{\frac{K^2 \mu_i^2}{n^2\lambda_k^2}} )(\theta_a \theta_i)(\frac{1}{\sqrt{n \otheta (\otheta \vee \theta_i)}}) \\
         & + \sum_{x \neq i}^n  \sqrt{\frac{K^2 \mu_i \mu_x}{n^2\lambda_k^2}} (\theta_a \theta_x)(\frac{1}{\sqrt{n \otheta (\otheta \vee \theta_x)}}) \\
         & \lesssim (\frac{1}{\lambda_k} + \frac{K \mu_i}{n\lambda_k})(\frac{\theta_a \theta_i}{\sqrt{n \otheta (\otheta \vee \theta_i)}})  + \sqrt{\frac{K^2 \mu_i \theta_a^2}{n\lambda_k^2}} \numberthis \label{theorem9extreq1}
    \end{align*} 
        \normalsize where the bound on the summation in the penultimate line follows from Lemma \ref{basicsumresults}. Lemma \ref{basicsumresults} implies $|\trueregdegree^{-\frac{3}{2}}(a, a) \xi_k(a)| \lesssim \frac{1}{n^{2}\otheta^2(\otheta \vee \theta_a)}$, so in conjunction with \eqref{theorem9extreq1}, we obtain
    \begin{align}
        B_{1, ab} \lesssim \frac{\theta_i}{n^{\frac{5}{2}}\otheta^2\sqrt{\otheta (\otheta \vee \theta_i)}}(\frac{1}{\lambda_k} + \frac{K \mu_i}{n\lambda_k}) + \frac{K \sqrt{\mu_i}}{n^{\frac{5}{2}}\lambda_k \otheta^2} \numberthis \label{theorem9extreq3}
    \end{align}
To bound $B_{2, ab}$, we apply Corollary \ref{corN1N2} and Lemmas \ref{basicresults}, \ref{basicexpansiondegreeresults} and \ref{basicsumresults} again,
\begin{align*}
    B_{2, ab} & \lesssim \delta n^{-1}  \sum_{a = 1}^n \xi_k(a) \trueregdegree^{-\frac{3}{2}}(a, a)\sum_{x = 1}^n [N_{1k}]_{i, x} H_{ax} \trueregdegreehalf(x, x) \\
    & = \delta n^{-1} \sum_{a = 1}^n \xi_k(a) \trueregdegree^{-\frac{3}{2}}(a, a) \left((\frac{1}{\lambda_k} + \frac{K \mu_i}{n\lambda_k})(\frac{\theta_a \theta_i}{\sqrt{n \otheta (\otheta \vee \theta_i)}})  + \sqrt{\frac{K^2 \mu_i \theta_a^2}{n\lambda_k^2}}\right) \\
    & \lesssim \frac{\theta_i}{n^{\frac{5}{2}}\otheta^2\sqrt{\otheta (\otheta \vee \theta_i)}}(\frac{1}{\lambda_k} + \frac{K \mu_i}{n\lambda_k}) + \frac{K \sqrt{\mu_i}}{n^{\frac{5}{2}}\lambda_k \otheta^2}  \end{align*}
where the bound in the second line follows from \eqref{theorem9extreq1}. Combining this bound with \eqref{theorem9extreq3} yields a bound for $|d_{ab}|$ and thus $|c_{ab}|$:
    $$
|c_{ab}| \lesssim \frac{\theta_i}{n^{\frac{5}{2}}\otheta^2\sqrt{\otheta (\otheta \vee \theta_i)}}(\frac{1}{\lambda_k} + \frac{K \mu_i}{n\lambda_k}) + \frac{K \sqrt{\mu_i}}{n^{\frac{5}{2}}\lambda_k \otheta^2}
$$
In turn, we can bound the variance by
\begin{align*}
    \hbox{var} & \lesssim \sum_{1 \leq a \leq b \leq n} c_{ab}^2 \theta_a\theta_b \leq \sum_{1 \leq a, b \leq n} c_{ab}^2 \theta_a\theta_b \\
    & \lesssim \sum_{1 \leq a, b \leq n} \left[\frac{\theta_i}{n^{\frac{5}{2}}\otheta^2\sqrt{\otheta (\otheta \vee \theta_i)}}(\frac{1}{\lambda_k} + \frac{K \mu_i}{n\lambda_k})\right]^2 \theta_a\theta_b + \sum_{1 \leq a, b \leq n} \left[\frac{K \sqrt{\mu_i}}{n^{\frac{5}{2}}\lambda_k \otheta^2}\right]^2 \theta_a\theta_b \\
    & \lesssim \frac{\theta_i^2}{n^{3}\otheta^3(\otheta \vee \theta_i)}(\frac{1}{\lambda_k} + \frac{K \mu_i}{n\lambda_k})^2 + \frac{K^2 \mu_i}{n^{3}\lambda_k^2 \otheta^2}
\end{align*}
Applying Bernstein and simplifying shows
\begin{align}
     [\bN_{1k}\boldsymbol{\trueregdegreehalf} \bH \bDelta_{D, 1} \bxi_k](i) & \lesssim \left(\frac{\theta_i}{n^{\frac{3}{2}}\otheta^{\frac{3}{2}}\sqrt{\otheta \vee \theta_i}}(\frac{1}{\lambda_k} + \frac{K \mu_i}{n\lambda_k}) + \frac{K \sqrt{\mu_i}}{n^{\frac{3}{2}}\lambda_k \otheta}\right) \sqrt{\log n} \label{theorem9bound3}
\end{align}
\end{enumerate}
Combining our bounds for the RHS terms in \eqref{theorem9begbound}, we finally obtain:
\begin{align*}
    & [\Delta\bP_k'](i) \lesssim (\frac{1}{n\otheta\lambda_k} + \frac{K\mu_i^{\frac{1}{2}}}{n^{\frac{3}{2}} \lambda_k\otheta} + \frac{K \mu_i}{n^2 \lambda_k\otheta} + \frac{K \mu_i^{\frac{3}{2}}}{n^{\frac{5}{2}} \lambda_k\otheta}) \sqrt{\log n} + (\frac{1}{n^{\frac{3}{2}}\otheta^2 d_k} + \frac{K\sqrt{\mu_i}}{n^{\frac{5}{2}}\otheta^2 d_k}) \log n \\
    & + \left(\frac{1}{n\lambda_k\otheta} + \frac{1}{n^{\frac{3}{2}}\lambda_k\otheta}\sqrt{\frac{\minthetai}{\otheta}} + \frac{K \sqrt{\mu_i}}{n^2\lambda_k\otheta}\sqrt{\frac{\theta_i}{\otheta}} + \frac{K \theta_i \sqrt{\mu_i}}{n^{\frac{5}{2}} \lambda_k \otheta^2}\right) \sqrt{\log n} \\
     & + \left(\frac{1}{n^{\frac{3}{2}}\lambda_k \otheta^{2}} + \frac{K\sqrt{\mu_i}}{n^{\frac{5}{2}}\lambda_k \otheta^2} \right) \log n + \left(\frac{\theta_i}{n^{\frac{3}{2}}\otheta^{\frac{3}{2}}\sqrt{\otheta \vee \theta_i}}(\frac{1}{\lambda_k} + \frac{K \mu_i}{n\lambda_k}) + \frac{K \sqrt{\mu_i}}{n^{\frac{3}{2}}\lambda_k \otheta}\right) \sqrt{\log n} \\
\end{align*}
Thanks to Assumption \ref{regconds2}(d), which states $C\sqrt{\frac{\log n}{n}} \lesssim \otheta$, we can simplify the expression greatly. In fact, when $K$ is fixed, the expression simply boils down to $[\Delta\bP_k'](i) \lesssim \frac{\sqrt{\log n}}{n\otheta\lambda_k}$. Alongside our bound for $\norm{\delta}_{\infty}$ from (\ref{theorem9eq10}), we can finally conclude 
\begin{align*}
    |\hat{\xi}_k(i) - \xi_k(i)| & \lesssim \frac{\sqrt{\log n}}{n\otheta\lambda_k} + (\max_{i \leq n} \left(|[(\wh\bP_k-\bP_k-\Delta\bP_k)\bxi_k](i)|+ |[(\Delta\bP_k - \Delta\bP_k')\bxi_k](i)|\right) \\
    & + \|(1-\wh\bxi_k^\top\bxi_k)\wh\bxi_k\|_\infty) \\
    & \lesssim \frac{\sqrt{\log n}}{n\otheta\lambda_k} + \max_{i \leq n} \{\chi_{\scriptscriptstyle 2_{i, k}}\} + \frac{\log^{1.5}n}{\lambda_k^{3}n^{\frac{3}2}\otheta^3} + (\frac{\log n \sqrt{\theta_{\text{max}}}}{n^{\frac{3}{2}}\otheta^{\frac{5}{2}}\lambda_k^{2}} + \frac{\log^{\frac{3}{2}}n}{n^{\frac{3}{2}}\otheta^3})
\end{align*}
The purpose of Assumption \ref{regconds3} is to ensure that the RHS of the above equation greatly simplifies -- specifically, to the point where the first order terms dominate the second order terms, meaning
\begin{equation}
   |\hat{\xi}_k(i) - \xi_k(i)| \lesssim \frac{\sqrt{\log n}}{n\otheta\lambda_k} \label{theorem9bound4}
\end{equation}
To show the deviation bounds on the normalized eigenvectors, we use triangle inequality
\begin{align*}
\Vert \hat{r}_i(k - 1) - r_i(k - 1) \Vert  = \Vert \frac{\hat{\xi}_k(i)}{\hat{\xi}_1(i)}-\frac{\xi_k(i)}{\xi_{1}(i)} \Vert & \leq \norm{\hat{\xi}_k(i) - \xi_k(i)}\norm{\frac{1}{\hat{\xi}_1(i)}} + \norm{\xi_k(i)}\norm{\frac{1}{\hat{\xi}_1(i)} - \frac{1}{\xi_1(i)}}
\end{align*}
To bound the first RHS term, we note that $\hat{\xi}_k(i) \asymp \xi_k(i)$, since $|\hat{\xi}_k(i) - \xi_k(i)| \leq \frac{\sqrt{\log n}}{n\otheta\lambda_k} \ll \xi_k(i)$ by \eqref{theorem9bound4}. In conjunction with Lemma \ref{lem:order}, it therefore follows that
\begin{align*}
\norm{\hat{\xi}_k(i) - \xi_k(i)}\norm{\frac{1}{\hat{\xi}_1(i)}} \leq \frac{\sqrt{C\log n}}{n\otheta\lambda_k} \cdot \sqrt{\frac{Cn\otheta^2}{\theta_i(\minthetai)}} = C'\sqrt{\frac{\log n}{n\theta_i(\minthetai)\lambda_K}}
\end{align*}
For the second RHS term,
\begin{align*}
\norm{\xi_k(i)}\norm{\frac{1}{\hat{\xi}_1(i)} - \frac{1}{\xi_1(i)}}  & \leq C\norm{\xi_k(i)}\cdot \frac{\norm{\hat{\xi}_1(i) - \xi_1(i)}}{\xi_1(i)^2} \\
& \leq C'\sqrt{\frac{\log n}{n\theta_i(\minthetai)\lambda_K}}
\end{align*}
under event $\mathcal{A}_1$, for all $i \in [n]$. Combining the above inequalities, we get the desired bounds on $\hat{r}_i(k - 1)$.

\subsection{Proof of Theorems \ref{expansiondegrees} and \ref{basicexpansiondegreeresults}}
\subsubsection{Proof of Theorem \ref{expansiondegrees}}\label{expansiondegreespf}
We obtain the desired expansion through Taylor expansion. Define $f(x) := (x)^{-\frac{1}{2}}$, and set $x_{0i} = D_0(i, i)$ and $x_{1i} = \hat{D}(i, i)$. By Taylor expansion, we know that there exists some $\tilde{x}_i$ between $x_{0i}$ and $x_{1i}$, such that
\begin{align*}
    f(x_{1i})-f(x_{0i}) & = f'(x_{0i})(x_{1i}-x_{0i})+\frac{f''(x_{0i})}{2}(x_{1i}-x_{0i})^2 + \frac{f'''(\tilde{x}_i)}{2}(x_{1i}-x_{0i})^3 \\
    & := \bDelta_{D, 1} + \bDelta_{D, 2} + \phi_D \numberthis \label{expansiondegreespfeq1}
\end{align*}

Since $f(x_{1i})-f(x_{0i}) = \hat{D}(i, i)^{\frac{-1}{2}} - D_0(i, i)^{\frac{-1}{2}}$ by construction, \eqref{expansiondegreespfeq1} is identical to \eqref{expansiondegreeseq0}.
Plugging in the relation $x_{1i} - x_{0i} := \hat{D}(i, i) - D_0(i, i)$ into \eqref{expansiondegreespf} yields the desired expressions for $\bDelta_{D, 1}, \bDelta_{D, 2}$ in equations \eqref{expansiondegreeseq1}, \eqref{expansiondegreeseq2}, so we are done.

\subsubsection{Proof of Theorem \ref{basicexpansiondegreeresults}}\label{expansiondegreesresultspf}
\underline{Order 0 terms:} The bound on $|N(i, i)|$ clearly follows from \eqref{basicresultseq2}. To bound $A$, we simply note 
\begin{align*}
    \sum_{i \leq n} \theta_i N(i, i)^2 & =  \sum_{\theta_i \leq \otheta} \theta_i N(i, i)^2 + \sum_{\theta_i > \otheta} \theta_i N(i, i)^2 \\
    & \leq \sum_{\theta_i \leq \otheta} \frac{\theta_i^2}{n\otheta^2} + \sum_{\theta_i \geq \otheta} \frac{\theta_i}{n\otheta} \\
    & \leq \sum_{\theta_i \leq \otheta}\frac{\otheta^2}{n\otheta^2} + \sum_{\theta_i} \frac{\theta_i}{n\otheta} \leq 2
\end{align*} holds, where the penultimate line follows from plugging in the bounds in \eqref{basicresultseq4} and \eqref{basicresultseq6}.

\underline{Order 1 terms:} 
First, we bound $\bDelta_{D, 1}$. 
\begin{align*}
    |\bDelta_{D, 1}(i, i)| & := |-\frac{1}{2} \cdot \frac{(\sum_{j = 1}^n W_{ij} + \delta / n \sum_{i = 1}^n \sum_{j = 1}^n W_{ij})}{D_0(i, i)^{\frac{3}{2}}}| \\
    & \leq |\frac{1}{2}\frac{(C\otheta\theta_i\sqrt{n \log n} + C' \delta \otheta \sqrt{\log n})}{D_0(i, i)^{\frac{3}{2}}}| \\
    & \lesssim \sqrt{\log n} \cdot \min\{\frac{1}{n\otheta\theta_i}, \frac{1}{n\otheta^2}\}| 
\end{align*}
where the penultimate line follows from plugging in the bounds in \eqref{basicresultseq4} and \eqref{basicresultseq6}; the last line from plugging in the bound on $|D_0(i, i)|$ proven in Lemma \ref{basicexpansiondegreeresults} and simplifying.

To bound $|R_{D, 1}(i, i)|$, we express it in terms of an order 1 Taylor expansion. In particular, define $f(x) := x^{-\frac{1}{2}}$, and set $x_{0i} = D_0(i, i)$ and $x_{1i} = \hat{D}(i, i)$, such that $f(x_{1i})-f(x_{0i}) = R_{D, 1}(i, i)$ by construction. By Taylor expansion, we know that there exists some $\tilde{x}_i$ between $x_{0i}$ and $x_{1i}$, such that
\begin{align*}
    |f(x_{1i})-f(x_{0i})| & = |f'(\tilde{x}_i)(x_{1i}-x_{0i})|
\end{align*}
We can bound the RHS of the above equation:
\begin{align*}
   |f'(\tilde{x}_i)(x_{1i}-x_{0i})| & = |-\frac{1}{2} \cdot \frac{1}{\tilde{x}_i^{\frac{3}{2}}}(x_{1i}-x_{0i})| \\
   & \leq |\frac{1}{2} \max\{D_0^{-\frac{1}{2}}(i, i), \hat{D}^{-\frac{1}{2}}(i, i)\}^3(x_{1i}-x_{0i})| \\
   & \lesssim |\frac{1}{2} (\min\{\frac{1}{n\otheta\theta_i}, \frac{1}{n\otheta^2}\})^3 (x_{1i}-x_{0i})| \\
   & \lesssim |\frac{1}{2} (\min\{\frac{1}{n\otheta\theta_i}, \frac{1}{n\otheta^2}\})^3 (\sqrt{n\otheta\theta_i \log n})| \\
   & \lesssim \sqrt{\log n} \cdot \min\{\frac{1}{n\otheta\theta_i}, \frac{1}{n\otheta^2}\}| 
\end{align*}
where the third line follows from the bound on $|D_0(i, i)|$ in Lemma \ref{basicexpansiondegreeresults}, and the fourth line follows from \eqref{basicresultseq8}.

To bound $A$, we plug in the bound for $|N(i, i)|$ we just established, then simply sum up across all indices $i \leq n$. Some routine computation yields the desired bound.

\underline{Order 2 terms:} 
The argument for $|\Delta_{D, 2}(i, i)|$ is similar to the order 1 case.
\begin{align*}
    |\bDelta_{D, 2}(i, i)| & := |\frac{3}{4} \cdot \frac{(\sum_{j = 1}^n W_{ij} + \delta / n \sum_{i = 1}^n \sum_{j = 1}^n W_{ij})^2}{D_0(i, i)^{\frac{5}{2}}}| \\
    & \leq |\frac{1}{2}\frac{(C\otheta\theta_i\sqrt{n \log n} + C' \delta \otheta \sqrt{\log n})^2}{D_0(i, i)^{\frac{5}{2}}}| \\
    & \lesssim \log n \cdot (\min\{\frac{1}{n\otheta\theta_i}, \frac{1}{n\otheta^2}\})^{\frac{3}{2}}
\end{align*}
where the last line follows from plugging in the bounds in \eqref{basicresultseq4}, \eqref{basicresultseq6}, and bound on $|D_0(i, i)|$ proven in Lemma \ref{basicexpansiondegreeresults}.

Likewise, the argument for $|R_{D, 2}(i, i)|$ is similar to the order 1 case. Defining $f(x), x_{0i}, x_{1i}$ as in the proof for the $|R_{D, 1}(i, i)|$ bound, Taylor expansion implies the existence of some $\tilde{x}_i$ between $x_{0i}$ and $x_{1i}$, such that
\begin{align*}
    |R_{D, 2}(i, i)| := |f(x_{1i})-f(x_{0i}) - f'(x_{0i})(x_{1i}-x_{0i})| & = |\frac{f''(\tilde{x}_i)}{2}(x_{1i}-x_{0i})^2|
\end{align*} Following a similar argument as to the proof for the $|R_{D, 1}(i, i)|$ bound, it is straightforward to deduce that $ |R_{D, 2}(i, i)| \lesssim \log n \cdot (\min\{\frac{1}{n\otheta\theta_i}, \frac{1}{n\otheta^2}\})^{\frac{3}{2}}$, as desired.

To bound $A$, we plug in the bound for $|N(i, i)|$ we just established, then simply sum up across all indices $i \leq n$. Some routine computation yields the desired bound.

\underline{Order 3 terms:} Defining $f(x), x_{0i}, x_{1i}$ as in the proof for the $|R_{D, 1}(i, i)|$ bound, Taylor expansion implies the existence of some $\tilde{x}_i$ between $x_{0i}$ and $x_{1i}$, such that
\begin{align*}
    |\phi_D| & := |f(x_{1i})-f(x_{0i}) - f'(x_{0i})(x_{1i}-x_{0i}) - f''(x_{0i})(x_{1i}-x_{0i})| \\
    & = |\frac{f'''(\tilde{x}_i)}{6}(x_{1i}-x_{0i})^3| \lesssim \sqrt{\log^3 n}(\min\{\frac{1}{n\otheta\theta_i}, \frac{1}{n\otheta^2}\})^{2}
\end{align*}
where the last inequality follows from the bound on $|D_0(i, i)|$ in Lemma \ref{basicexpansiondegreeresults}, alongside the bound in \eqref{basicresultseq8}. To obtain the desired bound for $A$, we plug in the just-shown bound for $|N(i, i)|$ into the expression for $A$, then simply sum up across all indices $i \leq n$.

\section{Supplementary Lemmas for Theorem \ref{improveddevboundsproof}}
 In order to prove Theorem \ref{improveddevboundsproof}, we make use of several supplementary lemmas, which we establish below. First, we start with Lemmas \ref{e2concentrationlemma2}-\ref{Wconcentration3}, which seek to bound the norm of the $i$th row of certain matrices.
\\
\begin{lem}\label{e2concentrationlemma2}
Let $(\bM, \bN)$ be a pair of matrices. Then with probability $1 - o(n^{-10})$, the following bounds on $\norm{[\bM \bH \bN]_{i, \cdot}}$ hold in the sense of $\lesssim$:
\begin{center}
    \begin{tabular}{ c|c|c } 
$(\bM, \bN)$ & $\norm{[\bM \bH \bN]_{i, \cdot}}$ \\ \hline
 $(\hbox{order 1 term}, \boldsymbol{\trueregdegreehalf})$ & $\sqrt{\frac{\log n}{n^2\otheta^2}}$ \\
 $(\boldsymbol{\trueregdegreehalf}, \hbox{order 1 term})$ & $\sqrt{\frac{\theta_i\log n}{n^2\otheta^3}}$  \\
 $(\hbox{order 1 term}, \hbox{order 1 term})$ & $\sqrt{\frac{\log^2 n}{n^3\otheta^4}}$  \\
 $(\boldsymbol{\trueregdegreehalf}, \hbox{order 2 term})$ & $\sqrt{\frac{\theta_i\log^2 n}{n^3\otheta^5}}$  \\
 $(\hbox{order 2 term}, \boldsymbol{\trueregdegreehalf})$ & $\sqrt{\frac{\log^2 n}{n^3\otheta^4}}$
\end{tabular}
\end{center}
where ``order 1" and ``order 2" terms are defined in Subsection \ref{sectionexpansiondegrees}.
\end{lem}

\begin{lem}\label{matrixbernstein1}
With probability $1 - o(n^{-10})$, the following bounds hold simultaneously for all $i \in [n]$:
\begin{itemize}
    \item $\norm{[\boldsymbol{\trueregdegreehalf} \bW \boldsymbol{\trueregdegreehalf} \bA]_{i, :}}_2 \lesssim \sqrt{\frac{\log^2(n)}{n\otheta^2}} \displaystyle \max_{l \leq n} \left\{ \frac{1}{\sqrt{\max\{n\otheta\theta_l, n\otheta^2\}}} \norm{\bA}_{2, \infty} \right\} + \frac{\sqrt{\log n}}{n\otheta} \norm{\bA}_{F}$, where $\bA$ is any fixed matrix independent of the entries in the ith row and column of the noise matrix $\bW$.
    \item $\norm{[\boldsymbol{\trueregdegreehalf} W \bDelta_{D, 1} \bxi_k]_{i, :}}_2 \lesssim \sqrt{\frac{\log^2 n}{n^3\otheta^4}}$, where $\bxi_k$ is any ground truth eigenvector.
    \item $\norm{[\bDelta_{D, 1} W \boldsymbol{\trueregdegreehalf} \bxi_k]_{i, :}}_2 \lesssim  \sqrt{\frac{\log^2(n)}{n^3\otheta^4}}(1 + \sqrt{\frac{1}{n\otheta\theta_i}})$.
\end{itemize}
\end{lem}

\begin{lem}\label{spectralnormE}
For any $i \leq n$, the following bounds hold with probability $1 - o(n^{-11})$, 
\begin{align*}
\norm{[\boldsymbol{\trueregdegreehalf} \bW \boldsymbol{\trueregdegreehalf}]_{i, \cdot}} & \lesssim \sqrt{\frac{\log n}{n\otheta^2}} \\
\norm{[\bE - \boldsymbol{\trueregdegreehalf} \bW \boldsymbol{\trueregdegreehalf}]_{i, \cdot}} & \lesssim \frac{\log n}{n\otheta^2}\\
\norm{[\bE - \bE_1]_{i, \cdot}} & \lesssim \frac{\log n}{n\otheta^2}\\
\norm{[\bE - \bE']_{i, \cdot}} & \lesssim \sqrt{\frac{\log^{3}n}{n^{3}\otheta^6}} \\
\norm{\bE - \bE_1} & \lesssim \frac{\log n}{n\otheta^2}
\end{align*}
\end{lem}

\begin{lem}\label{Wconcentration3}
    Suppose $A = \bDelta_{D, 1}^{(i)}\bxi_k$, where $\bxi_k$ is any ground truth eigenvector. For any $i\in [n]$, we have with probability at least $1-O(n^{-15}m)$,
    \begin{align*}
        \left\|\bW_{i,\cdot}\boldsymbol{A}\right\|_{2}\lesssim \sqrt{\log n}\theta_{\text{max}}\left\|\boldsymbol{A}\right\|_F +\log n\left\|\boldsymbol{A}\right\|_{2, \infty}
    \end{align*}
\end{lem}

Lemma \ref{Wconcentration2} provides several useful results for bounding the second order expansions in the proof of Theorem \ref{improveddevbounds}.
\\
\begin{lem}\label{Wconcentration2}
    For any $i\in [n]$ and ground truth eigenvector $\bxi_k$, we have with probability at least $1-o(n^{-10})$,
    \begin{align*}
        \norm{[\bE_1\bxi_k]_{i, \cdot}} & \lesssim \frac{1}{n\otheta}(1 + \sqrt{\frac{\theta_i}{\otheta}}) \\
        \norm{\bE_1\bxi_k}_2 & \lesssim \sqrt{\frac{\log n}{n\otheta^2}}\\ 
        \norm{[(\bE - \bE_1)\bxi_k]_{i, \cdot}} & \lesssim \sqrt{\frac{\log^2 n}{n^3\otheta^4}}(1 + \sqrt{\frac{1}{n\otheta\theta_i}}) + \sqrt{\frac{\log^3n}{n^3\otheta^6}}
    \end{align*}
\end{lem}

The following lemma and corollary bounds the effect of multiplying an arbitrary matrix $\bA$ by the $\bN_{1k}$ and $\bN_{2k}$ matrices, as defined in \eqref{N1N2}, on the norms of the rows of $\bA$.
\\
\begin{lem}\label{lemmaN1N2}
    For any $k \leq n$, let $\bN_{1k}$ and $\bN_{2k}$ be as defined in \eqref{N1N2}. Then,
    \begin{align*}
        & \norm{[\bN_{1k}-\frac{1}{\lambda_k}\bI]_{i, \cdot }}_2\lesssim \sqrt{\frac{(K-1)\mu_i}{n\lambda_k^{2}}}, \quad\; \left\|\bN_{1k}\right\|\lesssim \frac{1}{\lambda_k} \\
        &\norm{[\bN_{2k}-\frac{1}{\lambda_k^{2*}}\bI]_{i, \cdot }}_2 \lesssim \sqrt{\frac{(K-1)\mu_i}{n\lambda_k^{4}}},\quad \left\|\bN_{2k}\right\|\lesssim \frac{1}{\lambda_k^{2}}.
    \end{align*}
    Furthermore, for any $i, a \leq n$, the following entrywise bounds hold:
$$
[N_{1k}]_{j, a} \lesssim \begin{cases}
    \sqrt{\frac{K^2 \mu_j \mu_a}{n^2\lambda_k^2}} & \hbox{ if $j \neq a$} \\
    \frac{1}{\lambda_k} + \sqrt{\frac{K^2 \mu_j \mu_a}{n^2\lambda_k^2}} & \hbox{ otherwise}
\end{cases}
$$
\end{lem}
 From Lemma \ref{lemmaN1N2}, we immediately have the following corollary.
 \\
\begin{cor}\label{corN1N2}
    For $\bN_{1k}$ and $\bN_{2k}$ defined in \eqref{N1N2}, we have for all $\bx\in\mathbb{R}^n$
    \begin{align*}
        &\left\|\bN_{ik}\bx\right\|_{\infty}\lesssim \frac{1}{\lambda_1^{i}}\left\|\bx\right\|_{\infty}+\sqrt{\frac{(K-1)\max_{i \leq n}\mu_i}{n\lambda_1^{2i}}}\left\|\bx\right\|_{2}, \quad i=1,2. 
    \end{align*}
\end{cor}

\newpage
\subsection{Proofs of Supplementary Lemmas}

\subsubsection{Proof of Lemma \ref{e2concentrationlemma2}} \label{e2concentrationlemma2pf}
   Since $H_{ik} \lesssim \theta_i\theta_k$,
   \begin{align*}
       \norm{[\bM\bH \bN]_{i, \cdot}} & = \sqrt{\sum_{j \leq n} (M(i, i) H_{ij} N(j, j))^2} \\
       & \lesssim \sqrt{M(i, i)^2 \theta_i^2 \sum_{j \leq n} \theta_j^2 N(j, j)^2}
   \end{align*}
 Plugging in the bounds in Lemma \ref{basicexpansiondegreeresults} leads to the desired row norm bounds.

\subsubsection{Proof of Lemma \ref{matrixbernstein1}} \label{matrixbernstein1proof}
\underline{When $(\bM, \bN) = (\boldsymbol{\trueregdegreehalf}, \boldsymbol{\trueregdegreehalf})$:} $[\bM \bW \bN \bA]_{i, \cdot}$ is the sum of independent random matrices $\sum_{l = 1}^n W_{il} \cdot (N(l, l) M(i, i) [\bA]_{l, :})$, so we can bound the norm directly using the matrix Bernstein inequality (\cite{tropp2015introductionmatrixconcentrationinequalities}; Theorem 6.1.1). In particular, we can bound the parameters $\mathcal{L}$ and $\mathcal{V}$ in the matrix Bernstein inequality as
\begin{align*}
    L & := \max_{l \leq n} \norm{W_{il} \cdot (N(l, l) M(i, i) [\bA]_{l, })} \leq M(i, i) \max_{l \leq n} N(l, l) [\bA]_{l, } \\
    V & := \sum_{l = 1}^n E[\norm{W_{il} \cdot (N(l, l) M(i, i) [\bA]_{l, }]}^2] \leq \sum_{l = 1}^n \theta_i\theta_l (N(l, l) M(i, i) \norm{[\bA]_{l, }})^2\\
    & = \theta_i M(i, i)^2 \sum_{l = 1}^n \theta_l N(l, l)^2 \norm{[\bA]_{l, }}^2 \numberthis \label{matrixbernstein1eq1}
\end{align*}
where $\mathcal{L}$ is a uniform bound on every summand, and $\mathcal{V}$ a measure of the variance of the summation. By matrix Bernstein, it therefore follows that
$$[\bM \bW \bN \bA]_{i, \cdot} \lesssim \sqrt{\mathcal{V} \log n} +  \mathcal{L} \log n $$
The terms on the RHS of the expressions for $\mathcal{L}$ and $\mathcal{V}$ are the products of terms that are bounded in Lemma \ref{basicexpansiondegreeresults}, so plugging in their bounds yields the following bounds on $|\mathcal{L}|$ and $|\mathcal{V}|$:
\begin{align*}
    |\mathcal{L}| & \lesssim \frac{1}{\sqrt{n\otheta^2}} \max_{l \leq n} \left\{ \frac{1}{\max\{n\otheta\theta_l, n\otheta^2\}} \norm{\bA}_{2, \infty} \right\} \\
    \mathcal{V} & \lesssim \frac{1}{n^2\otheta^2} \norm{\bA}_{F}^2
\end{align*}
Applying Matrix Bernstein yields the desired bound.

\underline{When $M = \bDelta_{D, 1}, N = \boldsymbol{\trueregdegreehalf}$ and $\bA = \bxi_k$:} $[\bM\bW\bN\bA]_{i, \cdot}$ is again a sum of independent random vectors. As in the $(\bM, \bN) = (\boldsymbol{\trueregdegreehalf}, \boldsymbol{\trueregdegreehalf})$ case, we can therefore first bound $\bL_0$ and $V$ via equation \ref{matrixbernstein1eq1}, then plug in the bounds in Lemma \ref{basicexpansiondegreeresults} into equation \ref{matrixbernstein1eq1}. It follows that
\begin{align*}
    \mathcal{L} & \lesssim \frac{\sqrt{\log n}}{n\otheta^2} \max_{l \leq n} N(l, l) [\bA]_{l, } \\
    \mathcal{V} & \lesssim \frac{\log n}{n^3\otheta^4}\norm{\bA}_F^2
\end{align*}
Since $\bA = \bxi_k$, furthermore, we can refine our bound on $\mathcal{L}$ beyond a simple $2-\infty$ estimate. Specifically,
\begin{align*}
   \max_{l \leq n} N(l, l) [\bA]_{l, } & \lesssim \max_{l \leq n} \frac{1}{\sqrt{\max\{n\otheta\theta_l, n\otheta^2\}}} \cdot \sqrt{\frac{\mu_l}{n}} \\ 
   & = \frac{1}{\sqrt{n^2\otheta}}\sqrt{\frac{\mu_l}{\max{\theta_l, \otheta}}} \\
   & \lesssim \frac{1}{n\otheta}
\end{align*}
where the last line follows from Lemma \ref{incoherence}. As a result, $\mathcal{L} \lesssim \sqrt{\frac{\log n}{n^4\otheta^5}}$, implying
\begin{align*}
    \norm{[\bM\bW N\bxi_k]_{i, \cdot}}_2 & \lesssim \sqrt{\mathcal{V} \log n} + \mathcal{L} \log n \\
    & \lesssim \sqrt{\frac{\log^2 n}{n^3\otheta^4}}
\end{align*}
with probability $1 - o(n^{-10})$.

\underline{When $M = \boldsymbol{\trueregdegreehalf}, N = \bDelta_{D, 1}$ and $\bA = \bxi_k$:} We apply a leave one-out analysis. Define the leave-one-out matrix $\bDelta_{D, 1}^{(i)}$ by completely removing the effect of the $i$th row and column of $\bW$: \begin{align*}
&\bDelta_{D, 1}^{(i)}(j,j)=  
-\frac{\trueregdegree^{-\frac{3}{2}}(j, j)}2 \cdot \begin{cases}
\left(\sum_{s\neq i}W_{js}+ \frac 1n \cdot \sum_{s\neq i,l\neq i}W_{sl}\right), & \mbox{for }j\neq i,\cr
\frac 1n \cdot \sum_{s\neq i,l\neq i}W_{sl}, &\mbox{for }j=i. 
\end{cases}
\end{align*}
By triangle inequality,
\begin{align*}
    \norm{[\bM\bW\bDelta_{D, 1}\bxi_k]_{i, \cdot}} & \lesssim \norm{[\bM\bW\bDelta_{D, 1}^{(i)}\bxi_k]_{i, \cdot}} + \norm{[\bM\bW(\bDelta_{D, 1} - \bDelta_{D, 1}^{(i)})\bxi_k]_{i, \cdot}} \numberthis \label{Wconeq4}
\end{align*}
By Lemma \ref{Wconcentration3},
\begin{align*}
    \norm{[\bM\bW\bDelta_{D, 1}^{(i)}\bxi_k]_{i, \cdot}} & = M(i, i) \norm{[W\bDelta_{D, 1}^{(i)}\bA]_{i, \cdot}} \\
& \lesssim M(i, i) \sqrt{\frac{\log^2(n)\theta_i}{n^2\otheta^3}}(1 + \sqrt{\frac{1}{n\otheta\theta_i}}) \numberthis \label{matrixbernstein1eq2}
\end{align*}

To bound $\norm{[\bM\bW(\bDelta_{D, 1} - \bDelta_{D, 1}^{(i)})\bA]}$, we split $\bDelta_{D, 1} - \bDelta_{D, 1}^{(i)}$ into unregularized and regularized components. Denote $(\bDelta_{D, 1} - \bDelta_{D, 1}^{(i)})_{I}, (\bDelta_{D, 1} - \bDelta_{D, 1}^{(i)})_{R} \in \mathbb{R}^{n \times n}$, which are diagonal matrices defined by
\begin{align*}
    (\bDelta_{D, 1} - \bDelta_{D, 1}^{(i)})_{I}(l, l) & := -\frac{\trueregdegree^{-\frac{3}{2}}(l, l)}2 \cdot \begin{cases}
        W_{il} & \hbox{ if i $\neq$ l} \\
        \sum_{k = 1}^n W_{ik} & \hbox{ if i = l}
    \end{cases}\\
    (\bDelta_{D, 1} - \bDelta_{D, 1}^{(i)})_{R}(l, l) & := (-\frac{\delta}{n} \sum_{k = 1}^n W_{ik})\trueregdegree^{-\frac{3}2}(l, l)
\end{align*}
By plugging in the bounds in Lemmas \ref{basicresults} and \ref{basicexpansiondegreeresults}, we can deduce the following inequalities under event $A_1$:
\begin{align*}
    (\bDelta_{D, 1} - \bDelta_{D, 1}^{(i)})_{I}(l, l) & \lesssim \begin{cases}
        (n\otheta(\theta_l \vee \otheta))^{-\frac{3}{2}} & \hbox{ if i $\neq$ l} \\
        (n\otheta(\theta_i \vee \otheta))^{-1} \sqrt{\log n} & \hbox{ if i = l}
    \end{cases}\\
    (\bDelta_{D, 1} - \bDelta_{D, 1}^{(i)})_{R}(l, l) & \lesssim \sqrt{\frac{\theta_i \log n}{n^4\otheta^5}} \numberthis \label{Wconcentreq2}
\end{align*}
Now, 
\begin{align*}
     [\bM\bW(\bDelta_{D, 1} - \bDelta_{D, 1}^{(i)})_{I}\bA]_{i, :} & = M(i, i) \sum_{l \neq i} W_{ik}(\bDelta_{D, 1} - \bDelta_{D, 1}^{(i)})_{I}(l, l)\bA_{l, \cdot} \\
     & + M(i, i) \sum_{l = i} W_{ik}(\bDelta_{D, 1} - \bDelta_{D, 1}^{(i)})_{I}(l, l)\bA_{l, \cdot} \numberthis \label{Wconceq3} 
\end{align*}
The first term on the RHS is a sum of independent random vectors, so we can bound its norm via matrix Bernstein. In particular, under event $A_1$, 
\begin{align*}
    \mathcal{V} & \lesssim M(i, i)^2 \sum_{l \neq i} \theta_i\theta_l(\bDelta_{D, 1} - \bDelta_{D, 1}^{(i)})_{I}(l, l)^2\norm{A_{l, \cdot}}_2^2 \\
    |\mathcal{L}| & \lesssim M(i, i)(\bDelta_{D, 1} - \bDelta_{D, 1}^{(i)})_{I}(i, i)\norm{\bA}_{2, \infty}
\end{align*}
Each of the expressions on the RHS consists of terms that are bounded by equation \ref{Wconcentreq2} and Lemma \ref{basicexpansiondegreeresults}. Plugging in said bounds ultimately yields $$
|\mathcal{L}| \lesssim M(i, i)(n\otheta^2)^{-\frac{3}{2}}\norm{\bA}_{2, \infty}, V \lesssim M(i, i)^2(\frac{\theta_i}{n^3\otheta^5})\norm{\bA}_{F}^2
$$
The second RHS term in \ref{Wconceq3} can be bounded by $|M(i, i)W_{ik}||(\bDelta_{D, 1} - \bDelta_{D, 1}^{(i)})_{I}(i, i)|\norm{A_{l, \cdot}} \lesssim M(i, i)\frac{1}{n\otheta^2}\norm{\bA}_{2, \infty}$, so in all,
\begin{align*}
    \norm{[\bM\bW(\bDelta_{D, 1} - \bDelta_{D, 1}^{(i)})_{I}\bA]_{i, \cdot}} & \lesssim  M(i, i)\frac{1}{n\otheta^2}\norm{\bA}_{2, \infty} + M(i, i)\sqrt{\log n}(\frac{\theta_i}{n^3\otheta^5})^{\frac{1}{2}}\norm{\bA}_{F} \label{}
\end{align*}
To bound the unregularized term,
\begin{align*}
    \norm{[\bM\bW(\bDelta_{D, 1} - \bDelta_{D, 1}^{(i)})_{R}\bA]_{i, :}}_2 & \leq |M(i, i)|\norm{\bW_{i, \cdot}}\norm{(\bDelta_{D, 1} - \bDelta_{D, 1}^{(i)})_{R}}\norm{\bA} \\
& \lesssim (n\otheta\theta_i)^{-\frac{1}{2}} \cdot (n\otheta\theta_i)^{\frac{1}{2}} \cdot \sqrt{\frac{\log (n) \theta_i}{n^4\otheta^5}} \cdot \norm{\bA}_F \\
    & \lesssim \sqrt{\frac{\log (n)\theta_i}{n^4\otheta^5}} \norm{\bA}_F
\end{align*}
Combining our bounds for the regularized and unregularized terms, we therefore obtain
$$
\norm{[\bM\bW(\bDelta_{D, 1} - \bDelta_{D, 1}^{(i)})\bA]_{i, \cdot}} \lesssim \frac{M(i, i)}{n\otheta^2}\norm{\bA}_{2, \infty} + \sqrt{\frac{\log (n)\theta_i}{n^4\otheta^7}} \norm{\bA}_F
$$
Plugging the above and our bound in \eqref{matrixbernstein1eq2} into \eqref{Wconeq4}, we ultimately conclude
$$
\norm{[\bM\bW\bDelta_{D, 1} \bxi_k]_{i, \cdot}} \lesssim  \sqrt{\frac{\log^2(n)}{n^3\otheta^4}}(1 + \sqrt{\frac{1}{n\otheta\theta_i}})
$$

\subsubsection{Proof of Lemma \ref{spectralnormE}} \label{spectralnormEpf}
 First off,
 \begin{align*}
     \norm{[\boldsymbol{\trueregdegreehalf} \bW \boldsymbol{\trueregdegreehalf}]_{i, \cdot}} & \leq \trueregdegreehalf(i, i) \norm{\bW_{i, \cdot}}\norm{\boldsymbol{\trueregdegreehalf}} \numberthis \label{lemmaAeq2}
 \end{align*}
Each of the terms on the RHS are bounded in either Lemma \ref{basicresults} or \ref{basicexpansiondegreeresults}; plugging in their bounds, it follows that $\norm{[\boldsymbol{\trueregdegreehalf} \bW \boldsymbol{\trueregdegreehalf}]_{i, \cdot}} \lesssim \sqrt{\frac{\log n}{n\otheta^2}}$.
 
 Following the notation in Subsection \ref{sectionexpansiondegrees}, we note that $\bE - \boldsymbol{\trueregdegreehalf} \bW \boldsymbol{\trueregdegreehalf}$ can be expressed as 
 \begin{align*}
     \bE - \boldsymbol{\trueregdegreehalf} \bW \boldsymbol{\trueregdegreehalf} & = [\boldsymbol{\trueregdegreehalf}\bH\mathrm{\bR}_{D, 1} + \mathrm{\bR}_{D, 1} \bH \hat{D}^{-\frac{1}{2}}] + [\boldsymbol{\trueregdegreehalf} \bW \mathrm{\bR}_{D, 1} + \mathrm{\bR}_{D, 1} W \hat{D}^{-\frac{1}{2}}]
 \end{align*}
The norm of the ith row of all the terms in the first bracket can be bounded using Lemma \ref{e2concentrationlemma2}; plugging in the corresponding results yields an upper bound of $\frac{\log n}{n \otheta^2}$. To bound the second bracket, we use a simple spectral norm argument alongside the bounds in Lemmas \ref{basicresults}, \ref{expansiondegrees}, \ref{basicexpansiondegreeresults},
 \begin{align*}
     \norm{[\boldsymbol{\trueregdegreehalf} \bW \mathrm{\bR}_{D, 1} + \mathrm{\bR}_{D, 1} W \hat{D}^{-\frac{1}{2}}]_{i, \cdot}} & \lesssim \norm{\boldsymbol{\trueregdegreehalf}}\norm{\bW_{i, \cdot}}\norm{\mathrm{\bR}_{D, 1}} + \norm{\mathrm{\bR}_{D, 1}}\norm{\bW_{i, \cdot}}\norm{D_{\delta}^{-\frac{1}{2}}} \\
     & \lesssim \frac{\log n}{n \otheta^2} \numberthis \label{spectralnormEeq1}
 \end{align*}
 Combining \eqref{spectralnormEeq1} with our upper bound for the first bracket, we see
 $$\norm{[\bE - \boldsymbol{\trueregdegreehalf} \bW \boldsymbol{\trueregdegreehalf}]_{i, :}} \lesssim \frac{\log n}{n\otheta^2}$$
 
By a similar argument, one can also show the corresponding bounds for $\norm{[\bE-\bE_1]_{i, \cdot}}$ and $\norm{[\bE-\bE']_{i, \cdot}}$. 

Lastly, to show the spectral norm bound on $\bE - \bE_1$, we first decompose
\small \begin{align*}
    \norm{\bE - \bE_1} & \leq \norm{\boldsymbol{\trueregdegreehalf}\bH\mathrm{\bR}_{D, 2} + \boldsymbol{\trueregdegreehalf} \bW \mathrm{\bR}_{D, 1}} + \norm{\bDelta_{D, 1}\bH\mathrm{\bR}_{D, 1} + \bDelta_{D, 1}W\mathrm{\bR}_{D, 2}} + \norm{\mathrm{\bR}_{D, 2}\bX\boldsymbol{\obsregdegreehalf}} \\
    & := \mathcal{A}_1  + \mathcal{A}_2 + \mathcal{A}_3 \numberthis \label{W2conceq10}
\end{align*} \normalsize
We separately bound each of the terms on the RHS. 
\begin{itemize}
    \item \underline{Bound $\mathcal{A}_1$}: The first term can be further decomposed into
\begin{align*}
    \norm{\boldsymbol{\trueregdegreehalf}\bH\mathrm{\bR}_{D, 2} + \boldsymbol{\trueregdegreehalf} \bW \mathrm{\bR}_{D, 1}} & \leq \norm{\boldsymbol{\trueregdegreehalf}}\norm{\bH}\norm{\mathrm{\bR}_{D, 2}} + \norm{\boldsymbol{\trueregdegreehalf} \bW\bDelta_{D, 1}} \\
    & + \norm{\boldsymbol{\trueregdegreehalf}}\norm{\bW}\norm{\bR_{D, 2}}\\
    & \lesssim \frac{\log n}{n\otheta^2} + \norm{\boldsymbol{\trueregdegreehalf} \bW\bDelta_{D, 1}} + \sqrt{\frac{\theta_{\text{max}} \log n}{n^3\otheta^5}} \numberthis \label{Wconcentration2eq6}
\end{align*}
where the last line follows from plugging in the bounds in Lemmas \ref{basicresults} and \ref{basicexpansiondegreeresults}. To bound the second RHS term of (\ref{Wconcentration2eq6}), we first note 
\begin{align*}
    \norm{\boldsymbol{\trueregdegreehalf} \bW\bDelta_{D, 1}} & = \frac{1}{2} \cdot \norm{\boldsymbol{\trueregdegreehalf} \bW\bD^{-\frac{3}{2}}(\hat{\bD} - \bD)} \\
    & \leq \frac{1}{2}\norm{\boldsymbol{\trueregdegreehalf} \bW\bD^{-\frac{1}{2}}}\norm{\bD^{-1}(\hat{\bD} - \bD)} \numberthis \label{Wconcetr2eq7}
\end{align*}
Under event $A_1$, the quantity $\norm{\bD^{-1}(\hat{\bD} - \bD)}$ is bounded by $\sqrt{\frac{\log n}{n\otheta^2}}$ with probability $1 - o(n^{-10})$, as can be seen from plugging in the results in Lemma \ref{basicexpansiondegreeresults}. As for $\boldsymbol{\trueregdegreehalf} \bW \boldsymbol{\trueregdegreehalf}$, its norm is bounded by $\sqrt{\frac{\log n}{n\otheta^2}}$, as shown in the proof of Lemma B.3 of \citet{ke2022optimal}. Plugging these bounds into \eqref{Wconcetr2eq7}, we see $\norm{\boldsymbol{\trueregdegreehalf} \bW\bDelta_{D, 1}} \lesssim \frac{\log n}{\sqrt{n^3\otheta^6}}$. Returning to \eqref{Wconcentration2eq6} and simplifying, we obtain
\begin{align*}
    \norm{\boldsymbol{\trueregdegreehalf} \bW\bDelta_{D, 1}} \lesssim \sqrt{\frac{\log n}{n\otheta^2}}
\end{align*}

\item \underline{Bound $\mathcal{A}_2$ and $\mathcal{A}_3$}: By triangle inequality and submultiplicativty of the spectral norm,
\begin{align*}
     \mathcal{A}_2 & \leq \norm{\bDelta_{D, 1}}\norm{\bH}\norm{\mathrm{\bR}_{D, 1}} + \norm{\bDelta_{D, 1}}\norm{\bW}\norm{\mathrm{\bR}_{D, 2}} \numberthis \label{Wconcetr2eq9}
\end{align*}
Plugging in the bounds from Lemma \ref{basicexpansiondegreeresults} into the RHS yields an upper bound of $\frac{\log n}{n\otheta^2}$ for $\mathcal{A}_2$ (after simplification). Likewise, we can bound $\mathcal{A}_3$ by
\begin{align*}
    \mathcal{A}_3 & \leq \norm{\mathrm{\bR}_{D, 2}}(\norm{\bH} + \norm{\bW})\norm{\boldsymbol{\obsregdegreehalf}}
\end{align*}
Plugging in the bounds from Lemma \ref{basicexpansiondegreeresults} also results in an upper bound of $\frac{\log n}{n\otheta^2}$.
\end{itemize}
Combining our bounds for the various $\mathcal{A}$ terms, it follows that $\norm{\bE - \bE_1} \lesssim \frac{\log n}{n\otheta^2}$ with probability $1 - o(n^{-10})$.

\end{proof}

\subsection{Proof of Lemma \ref{Wconcentration3}} \label{Wconcentration3proof}
    $\bW_{i,\cdot}\boldsymbol{A} = \sum_{j = 1}^n W_{ij}\bA_{j, :}$ is the sum of independent random vectors, so we can bound its norm via matrix Bernstein inequality (\cite{tropp2015introductionmatrixconcentrationinequalities}; Theorem 6.1.1). Following the notation of \citep{tropp2015introductionmatrixconcentrationinequalities}, we know
    \begin{align*}
        \mathcal{L} & = \max_{1 \leq j \leq n} \norm{W_{ij}\bA_{j, :}}_2 \leq \max_{1 \leq j \leq n} \norm{\bA_{j, :}}_2 \leq \norm{\bA}_{2, \infty} \\
        V & = \sum_{j = 1}^n \norm{\bA_{j, :}}^2 H_{ij}(1 - H_{ij}) \leq\sum_{j = 1}^n C\theta_i\theta_j \norm{\bA_{j, :}}^2 \numberthis \label{Wconcentration3eq1}
    \end{align*}
    To bound $\norm{\bA_{j, :}}$ for $j \in [n]$, we can use the results in Lemmas \ref{basicexpansiondegreeresults} and \ref{incoherence}, thereby obtaining $\norm{\bA_{j, :}} \lesssim \sqrt{\frac{\log n}{n^3\otheta^4}}$. Plugging in our bound into \ref{Wconcentration3eq1}, it follows that $V \lesssim \frac{\log(n) \theta_i}{n^2\otheta^3}$. In turn, matrix Bernstein implies
    \begin{align*}
        \bW_{i,\cdot}\boldsymbol{A} & \lesssim \sqrt{V \log n} + \mathcal{L} \log n \\
        & \lesssim \sqrt{\frac{\log^2(n)\theta_i}{n^2\otheta^3}}(1 + \sqrt{\frac{1}{n\otheta\theta_i}})
    \end{align*}

\subsubsection{Proof of Lemma \ref{Wconcentration2}} \label{Wconcentration2proof}
We seek to bound the norm of
    \begin{align*}
        [\bE_1\bA]_{i, \cdot} & = [(\boldsymbol{\trueregdegreehalf} \bW \boldsymbol{\trueregdegreehalf})\bA]_{i, \cdot} + [(\bDelta_{D, 1} \bH \boldsymbol{\trueregdegreehalf})\bA]_{i, \cdot} + [(\boldsymbol{\trueregdegreehalf} \bH \bDelta_{D, 1})\bA]_{i, \cdot} 
    \end{align*}
    By Lemmas \ref{e2concentrationlemma2} and \ref{matrixbernstein1}, we can individually bound the norms of each of the terms on the RHS, yielding in all,
    \begin{align*}
        \norm{[\bE_1\bA]_{i, \cdot}} & \lesssim \frac{1}{\sqrt{n\otheta^2}} \max_{l \leq n} \left\{ \frac{1}{\max\{n\otheta\theta_i, n\otheta^2\}} \norm{\bA_{i, \cdot}} \right\} \log n + (\frac{1}{n\otheta} + \sqrt{\frac{\theta_i}{n^2\otheta^3}}) \sqrt{\log n} \norm{\bA}_{F} \numberthis \label{Wconcentration2proofeqA} 
    \end{align*}
    When $\bA = \bxi_k$ for some eigenvector, in particular, the RHS of \eqref{Wconcentration2proofeqA} can be bounded by $\frac{\sqrt{\log n}}{n\otheta}(1 + \sqrt{\frac{\theta_i}{\otheta}})$.

    To bound $\norm{\bE_1\bA}_2$, it suffices to bound $\norm{E}$, as $\norm{\bE_1\bA}_2 \leq \norm{\bE_1}_2\norm{\bA}_2$. To this end, we crudely bound $\norm{\bE_1}_2$ by its Frobenius norm. By summing across our row-bounds in \eqref{Wconcentration2proofeqA}, we can crudely bound the Frobenius norm, thereby obtaining
    \begin{align*}
        \norm{\bE_1}_F & = \sqrt{\sum_{i = 1}^n \norm{[\bE_1\bA]_{i, \cdot}}^2} \\
        & \lesssim \sqrt{\sum_{i = 1}^n \frac{\log n}{n^2\otheta^2}(1 + \frac{\theta_i}{\otheta})} \lesssim \sqrt{\frac{\log n}{n\otheta^2}}
    \end{align*}
    Since $A$ is an eigenvector, $\norm{\bE_1\bA}_2$ is therefore bounded from above by $\sqrt{\frac{\log n}{n\otheta^2}}$.

   Likewise, when $\bA = \bxi_k$ for some eigenvector $\bxi_k$, applying Lemmas \ref{e2concentrationlemma2}, \ref{matrixbernstein1} and \ref{spectralnormE} implies
    \begin{align*}
        \norm{[(\bE - \bE_1)\bA]_{i, \cdot}} & \leq \norm{[\bDelta_{D, 1} W \boldsymbol{\trueregdegreehalf} \bA]_{i, \cdot}} + \norm{[\boldsymbol{\trueregdegreehalf} \bW \bDelta_{D, 1} \bA]_{i, \cdot}} \\
        & + \norm{[\bDelta_{D, 1} \bH \bDelta_{D, 1} \bA]_{i, \cdot}} + \norm{[(\bE - \bE')\bA]_{i, \cdot}} \\
        & \lesssim  \sqrt{\frac{\log^2(n)}{n^3\otheta^4}}(1 + \sqrt{\frac{1}{n\otheta\theta_i}}) + \sqrt{\frac{\log^2 n}{n^3\otheta^4}} + \norm{\bE - \bE'}_{2. \infty} \norm{\bA}_{F} \\
        & \lesssim \sqrt{\frac{\log^2(n)}{n^3\otheta^4}}(1 + \sqrt{\frac{1}{n\otheta\theta_i}}) + \norm{\bE - \bE'}_{2. \infty} \norm{\bA}_{F} \\
        & \lesssim \sqrt{\frac{\log^2(n)}{n^3\otheta^4}}(1 + \sqrt{\frac{1}{n\otheta\theta_i}}) + \sqrt{\frac{\log^3n}{n^3\otheta^6}}
    \end{align*}

\subsubsection{Proof of Lemma \ref{lemmaN1N2}}
The spectral norm bounds follow directly from \eqref{N1N2}:
\begin{align*}
    \left\|\bN_{1k}\right\|  = \max_{2\leq i\leq n}\frac{1}{\lambda_k-\lambda_i}\lesssim \frac{1}{\lambda_k },\quad \left\|\bN_{2k}\right\|  = \max_{2\leq i\leq n}\frac{1}{(\lambda_k-\lambda_i)^2}\lesssim \frac{1}{\lambda_k^{2} }.
\end{align*}
By definition, for $\bN_{1k}$ we have
\begin{align}
    \bN_{1k}-\frac{1}{\lambda_k}\bI & = \sum_{1 \leq i \leq n, i \neq k}\frac{1}{\lambda_k-\lambda_i}\bxi_i\bxi_i^{T}  - \frac{1}{\lambda_k}\sum_{i=1}^n\bxi_i\bxi_i^{T}\\
    & = -\frac{1}{\lambda_k}\bxi_k\bxi_k^{T}+\sum_{1 \leq i \leq K, i \neq k}\frac{\lambda_i}{\lambda_k(\lambda_k-\lambda_i)}\bxi_i\bxi_i^{T} \label{corN1N2eq1}
\end{align}
On one hand, by Lemma \ref{incoherence} we have
\begin{align*}
    \left\|\frac{1}{\lambda_k}\bxi_k\bxi_k^{T}\right\|_{2,\infty} = \frac{1}{\lambda_k}\left\|\bxi_k\right\|_\infty\left\|\bxi_k\right\|_2\leq \sqrt{\frac{\mu_i}{n\lambda_k^{2}}}.
\end{align*}
On the other hand, defining $$\bC_1 := \textbf{diag}\left(\frac{\lambda_1}{\lambda_k(\lambda_k-\lambda_1)}, \frac{\lambda_2}{\lambda_k(\lambda_k-\lambda_2)},\dots, \frac{\lambda_K}{\lambda_k(\lambda_k-\lambda_K)}\right),$$ 
by Lemma \ref{incoherence} we have,
\begin{align*}
    \left\|\sum_{1 \leq i \leq K, i \neq k}\frac{\lambda_i}{\lambda_k(\lambda_k-\lambda_i)}\bxi_k\bxi_k^{T}\right\|_{2,\infty} &= \left\|\bXi_{-k}\bC_1\bXi_{-k}^\top\right\|_{2,\infty}\leq \left\|\bXi_{-k}\right\|_{2,\infty}\left\|\bC_1\bXi_{-k}^\top\right\| \lesssim \sqrt{\frac{(K-1)\mu_i}{n\lambda_k^{2}}}.
\end{align*}
As a result, we get $\left\|\bN_{1k}-\frac{1}{\lambda_k}\bI\right\|_{2,\infty}\lesssim\sqrt{\frac{(K-1)\mu_i}{n\lambda_k^{2}}}$. Similarly, for $\bN_{2k}$ we have 
\begin{align*}
    \bN_{2k}-\frac{1}{\lambda_k^{2}}\bI = -\frac{1}{\lambda_k^{2}}\bxi_k\bxi_k^{T}+\sum_{1 \leq i \leq K, i \neq k}\left(\frac{1}{(\lambda_k-\lambda_i)^2}-\frac{1}{\lambda_k^{2}}\right)\bxi_i\bxi_i^{T}.
\end{align*}
Again, we have $\|\frac{\bxi_k\bxi_k^{T}}{\lambda_k^{2}}\|_{2,\infty}\lesssim \sqrt{\frac{\mu_i}{n\lambda_k^{4}}}$. Defining
$$\bC_2 = \textbf{diag}\left(\frac{1}{(\lambda_k-\lambda_1)^2}-\frac{1}{\lambda_k^{2}},\frac{1}{(\lambda_k-\lambda_2)^2}-\frac{1}{\lambda_k^{2}}, \dots, \frac{1}{(\lambda_k-\lambda_K)^2}-\frac{1}{\lambda_K^{2}}\right)$$
where the $k$th term is excluded from the matrix. Thus, we have
\begin{align*}
    \left\|\sum_{1 \leq i \leq K, i \neq k}\left(\frac{1}{(\lambda_k-\lambda_i)^2}-\frac{1}{\lambda_k^{2}}\right)\bxi_k\bxi_k^{T}\right\|_{2,\infty} &= \left\|\bXi_{-k}\bC_2\bXi_{-k}^\top\right\|_{2,\infty}\leq \left\|\bXi_{-k}\right\|_{2,\infty}\left\|\bC_2\bXi_{-k}^\top\right\| \\
    & \lesssim \sqrt{\frac{(K-1)\mu_i}{n\lambda_k^{4}}}.
\end{align*}

Combining together, we get the desired conclusion
    $\|\bN_{2k}-\frac{1}{\lambda_k^{2}}\bI\|_{2,\infty}\lesssim\sqrt{\frac{(K-1)\mu_i}{n\lambda_k^{4}}}$.

To obtain entrywise bounds on $[N_{1k}]_{j, a}$, we split into cases. By \eqref{corN1N2eq1} and triangle inequality, 
\begin{align}
    |[N_{1k}]_{j, a}| & \leq |\frac{1}{\lambda_k}I_{j, a}| + |\frac{1}{\lambda_k}\xi_k(j)\xi_k(a)| + |[\bXi_{-k} \bC_1 \bXi_{-k}^T]_{j, a}| \label{corN1eq2}
\end{align}
where $\bXi_{-k}$ denotes $V$ with the kth column removed, i.e. the population eigenvector matrix with the kth eigenvector removed. When $a \neq j$, the first term on the RHS equals 0, while the third term can be bounded by
\begin{align*}
    [\bXi_{:, -k} \bC_1 \bXi_{:, -k}^T]_{j, a} = \Xi_{j, -k} \bC_1 \Xi_{a, -k} & \leq \norm{\Xi_{j, -k}}\norm{\bC_1}\norm{\Xi_{a, -k}} \\
    & \lesssim \sqrt{\frac{K\mu_j}{n}} \cdot \frac{1}{\lambda_k} \cdot \sqrt{\frac{K\mu_a}{n}} \lesssim \sqrt{\frac{K^2 \mu_j \mu_a}{n^2\lambda_k^2}} \numberthis \label{corN1eq3}
\end{align*}
where the bounds on $\norm{\Xi_{j, -k}}$ in the second line follow from Lemma \ref{lem:order}. Plugging in \eqref{corN1eq3} into \eqref{corN1eq2}, we therefore obtain the desired bound of $\sqrt{\frac{K^2 \mu_j \mu_a}{n^2\lambda_k^2}}$. 

When $a = j$, the bound in \eqref{corN1eq3} is still valid, so plugging it in again,
\begin{align*}
    |[N_{1k}]_{j, a}| & \lesssim |\frac{1}{\lambda_k}| + |\frac{1}{\lambda_k}\xi_k(j)\xi_k(a)| + \sqrt{\frac{K^2 \mu_j \mu_a}{n^2\lambda_k^2}} \\
    & \lesssim |\frac{1}{\lambda_k}| + \sqrt{\frac{K^2 \mu_j \mu_a}{n^2\lambda_k^2}}
\end{align*}

\subsubsection{Proof of Corollary \ref{corN1N2}} \label{corN1N2proof}
For all $i \in [n]$ and $j \in \{1, 2\}$, triangle inequality and Lemma \ref{N1N2} imply
\begin{align*}
    |[N_{jk}x](i)| & \leq |[(N_{jk} - \frac{1}{\lambda_k}I)x](i)| + |\frac{1}{\lambda_k} \cdot x(i)| \\
    & \leq \norm{N_{jk} - \frac{1}{\lambda_k}I}\norm{x} + |\frac{1}{\lambda_k} \cdot x(i)| \\
    & \lesssim \sqrt{\frac{(K - 1)\mu_i}{n\lambda_k^{2j}}}\norm{x} + |\frac{1}{\lambda_k} \cdot x(i)|
\end{align*}
Taking the maximum over all entries $i$ concludes.

\end{document}